\documentclass[aps,pra,letterpaper,notitlepage,longbibliography,reprint,onecolumn,superscriptaddress]{revtex4-1}
\usepackage[table]{xcolor}
\usepackage{graphicx,amsmath,amsfonts,amssymb,amsthm,bm,bbm}
\usepackage{tikz,pgfplots}
\pgfplotsset{compat=newest}
\pgfplotsset{plot coordinates/math parser=false}
\usetikzlibrary{arrows,decorations.pathmorphing,backgrounds,positioning,fit,calc,patterns,shapes,decorations.markings,shadings,snakes}
\usetikzlibrary{external}
\usepackage{ifthen}
\usepackage{multirow}
\usepackage{diagbox}
\usepackage{stmaryrd}
\usepackage{tabularx}
\usepackage{booktabs}
\newcommand{\greyhline}{\arrayrulecolor[gray]{.8}\hline\arrayrulecolor{black}}
\newcommand{\greycline}[1]{\arrayrulecolor[gray]{.8}\cline{#1}\arrayrulecolor{black}}

\newtheoremstyle{exampleA}
{\topsep}
{10pt}
{}
{}
{\bfseries}
{\\\\}
{ 10mm }
{\thmname{#1}\thmnumber{ #2}\thmnote{ (#3)}}

\theoremstyle{exampleA}
\newtheorem{exmp}{Example}[section]

\usepackage[pdfpagelabels,pdftex,bookmarks,breaklinks]{hyperref}
\definecolor{darkblue}{RGB}{0,0,127} 
\definecolor{darkgreen}{RGB}{0,180,0}
\definecolor{darkred}{RGB}{180,0,0}
\hypersetup{
	colorlinks, 
	linkcolor=darkblue, 
	citecolor=darkgreen, 
	filecolor=red, 
	urlcolor=blue,
	pdftitle={Fusing Binary Interface Defects in Topological Phases:The Z/pZ case}, 
	pdfauthor={Daniel Barter, Jacob C Bridgeman, and Corey Jones}
}

\definecolor{tctwistcolor}{RGB}{255,0,0}
\definecolor{tcmcolor}{RGB}{0,0,180}
\definecolor{tcecolor}{RGB}{0,180,0}
\definecolor{tcTTppcolor}{RGB}{184,134,11}
\definecolor{tcTTpmcolor}{RGB}{252,15,192}

\newcommand{\linf}{\,\,\rotatebox[origin=c]{-90}{$\curlyveedownarrow$}\,\,}

\theoremstyle{plain}

\theoremstyle{definition}

\newcommand{\ZZ}[1]{\mathbb{Z}/#1\mathbb{Z}}

\let\oldonlinecite\onlinecite
\renewcommand{\onlinecite}[1]{Ref.~[\oldonlinecite{#1}]}

\newcommand{\onlinecites}[1]{Refs.~[\oldonlinecite{#1}]}

\newcommand{\ket}[1]{|{#1}\rangle}

\newcommand{\D}[1]{\mathcal{D}{} }

\newcommand{\restrict}[1]{\raise-.2ex\hbox{\ensuremath|}_{#1}}

\definecolor{tensorblue}{rgb}{0.8,0.8,1}
\definecolor{tensorred}{rgb}{1,0.5,0.5}
\definecolor{tensorgreen}{rgb}{0.6,1,0.6}
\definecolor{tensorpurp}{rgb}{1,0.5,1}

\tikzset{ten/.style={fill=tensorblue}}
\tikzset{tenred/.style={fill=tensorred}}
\tikzset{tengreen/.style={fill=tensorgreen}}
\tikzset{tenpurp/.style={fill=tensorpurp}}

\makeatletter
\newcommand{\vast}{\bBigg@{4}}
\newcommand{\Vast}{\bBigg@{9}}
\makeatother

\def\Put(#1,#2)#3{\leavevmode\makebox(0,0){\put(#1,#2){#3}}}

\makeatletter
\def\pgf@plot@curveto@handler@finish{%
  \ifpgf@plot@started%
    \pgfpathcurvebetweentimecontinue{0}{0.995}{\pgf@plot@curveto@first}{\pgf@plot@curveto@first@support}{\pgf@plot@curveto@second}{\pgf@plot@curveto@second}%
  \fi%
}
\makeatother

\newlength\figureheight 
\newlength\figurewidth 
  
\newcommand{\includeTikz}[2]
{
	\tikzifexternalizing
	{
		\includeTikzrm{#1}{#2}
	}
	{
		\IfFileExists{figures/#1.pdf}{
			\includegraphics{figures/#1}
		}
		{
			\includeTikzrm{#1}{#2}
		}
	}
}

\newcommand{\includeTikzrm}[2]{
	\tikzset{external/remake next}
	\tikzsetnextfilename{#1}
	#2
}
  
  \definecolor{nicegreena}{RGB}{1,115,16}
  \definecolor{nicegreenb}{RGB}{1,240,16}
  
  \definecolor{nicegreen}{RGB}{60,183,82}
  
    \colorlet{ccred}{red!20}
    \colorlet{ccgreen}{green!50}
    \colorlet{ccblue}{blue!20}

  \tikzset{hexr/.style= {shape=regular polygon,regular polygon sides=6,minimum size=1cm, draw,inner sep=0,anchor=center,fill=red!50}}
  \tikzset{hexg/.style= {shape=regular polygon,regular polygon sides=6,minimum size=1cm, draw,inner sep=0,anchor=center,fill=green!50}}
  \tikzset{hexb/.style= {shape=regular polygon,regular polygon sides=6,minimum size=1cm, draw,inner sep=0,anchor=center,fill=blue!50}}
  
  \definecolor{tensorblue}{rgb}{0.8,0.8,1}
  \definecolor{tensorred}{rgb}{1,0.5,0.5}
  \definecolor{tensorpurp}{rgb}{1,0.5,1}
  
  \tikzset{nonesty/.style={fill=none,draw=none}}
  \tikzset{ten/.style={fill=tensorblue}}
  \tikzset{tenred/.style={fill=tensorred}}
  \tikzset{tengreen/.style={fill=green!50!black!50}}
  \tikzset{tenpurp/.style={fill=tensorpurp}}
  \tikzset{tengrey/.style={fill=black!20}}
  \tikzset{tenorange/.style={fill=orange!30}}
  \tikzset{u/.style={fill=blue!20,draw=black}}
  \tikzset{w/.style={fill=green!50!black!50,draw=black}}

\tikzset{external/system call={pdflatex \tikzexternalcheckshellescape -halt-on-error -interaction=batchmode -jobname "\image" "\texsource"; ps2pdf13 -dEmbedAllFonts=true -r100000  "\image".pdf "\image-13".pdf && cp "\image-13".pdf "\image".pdf && rm "\image-13".pdf && rm "\image".log && rm "\image".dpth && rm "\image"Notes.bib}}

\usetikzlibrary{calc}
\tikzstyle{inline text}=[text height=1.1ex, text depth=0.1ex, yshift=-.1ex]

\newcommand{\inflationalhs}[3]{
	\ifthenelse{\equal{#1}{}}{}
	{
		\draw[blue](-.75,1)--(0,.5);
		\node[above,inline text,blue] at (-.75,1) {#1};
	}
	\ifthenelse{\equal{#3}{}}{}
	{
		\draw[blue](.75,1)--(0,-.5);
		\node[above,inline text,blue] at (.75,1) {#3};
	}
	\draw[red,ultra thick] (0,-1)--(0,1);
	\node[below,inline text,red] at (0,-1) {#2};
}

\newcommand{\inflationarhs}[5]{
	\ifthenelse{\equal{#1}{}}{}
	{
		\draw[blue](-1.25,1)--(-.5,.5);
		\node[above,inline text,blue] at (-1.25,1) {#1};
	}
	\ifthenelse{\equal{#3}{}}{}
	{
		\draw[blue](-.5,-.75)--(.5,-.25);
		\node[above,inline text,blue] at (0,-.5) {#3};
	}
	\ifthenelse{\equal{#5}{}}{}
	{
		\draw[blue](1.25,1)--(.5,0);
		\node[above,inline text,blue] at (1.25,1) {#5};
	}
	\draw[red,ultra thick] (-.5,-1)--(-.5,1);
	\draw[orange,ultra thick] (.5,-1)--(.5,1);
	\node[below,inline text,red] at (-.5,-1) {#2};
	\node[below,inline text,orange] at (.5,-1) {#4};
}

\newcommand{\inflationarhss}[5]{
	\ifthenelse{\equal{#1}{}}{}
	{
		\draw[blue](-1.25,1)--(-.5,.5);
		\node[above,inline text,blue] at (-1.25,1) {#1};
	}
	\ifthenelse{\equal{#3}{}}{}
	{
		\draw[blue](-.5,-.75)--(.5,-.25);
		\node[above,inline text,blue] at (0,-.5) {#3};
	}
	\ifthenelse{\equal{#5}{}}{}
	{
		\draw[blue](1.25,1)--(.5,0);
		\node[above,inline text,blue] at (1.25,1) {#5};
	}
	\draw[red,ultra thick] (-.5,-1)--(-.5,1);
	\draw[red,ultra thick] (.5,-1)--(.5,1);
	\node[below,inline text,red] at (-.5,-1) {#2};
	\node[below,inline text,red] at (.5,-1) {#4};
}

\newcommand{\inflationblhs}[3]{
	\ifthenelse{\equal{#1}{}}{}
	{
		\draw[blue](-.75,-1)--(0,-.5);
		\node[below,inline text,blue] at (-.75,-1) {#1};
	}
	\ifthenelse{\equal{#3}{}}{}
	{
		\draw[blue](.75,-1)--(0,.5);
		\node[below,inline text,blue] at (.75,-1) {#3};
	}
	\draw[red,ultra thick] (0,-1)--(0,1);
	\node[above,inline text,red] at (0,1) {#2};
}

\newcommand{\inflationbrhs}[5]{
	\ifthenelse{\equal{#1}{}}{}
	{
		\draw[blue](-1.25,-1)--(-.5,-.5);
		\node[below,inline text,blue] at (-1.25,-1) {#1};
	}
	\ifthenelse{\equal{#3}{}}{}
	{
		\draw[blue](-.5,-.25)--(.5,.25);
		\node[above,inline text,blue] at (0,0) {#3};
	}
	\ifthenelse{\equal{#5}{}}{}
	{
		\draw[blue](1.25,-1)--(.5,.5);
		\node[below,inline text,blue] at (1.25,-1) {#5};
	}
	\draw[red,ultra thick] (-.5,-1)--(-.5,1);
	\draw[orange,ultra thick] (.5,-1)--(.5,1);
	\node[above,inline text,red] at (-.5,1) {#2};
	\node[above,inline text,orange] at (.5,1) {#4};
}

\newcommand{\inflationbrhss}[5]{
	\ifthenelse{\equal{#1}{}}{}
	{
		\draw[blue](-1.25,-1)--(-.5,-.5);
		\node[below,inline text,blue] at (-1.25,-1) {#1};
	}
	\ifthenelse{\equal{#3}{}}{}
	{
		\draw[blue](-.5,-.25)--(.5,.25);
		\node[above,inline text,blue] at (0,0) {#3};
	}
	\ifthenelse{\equal{#5}{}}{}
	{
		\draw[blue](1.25,-1)--(.5,.5);
		\node[below,inline text,blue] at (1.25,-1) {#5};
	}
	\draw[red,ultra thick] (-.5,-1)--(-.5,1);
	\draw[red,ultra thick] (.5,-1)--(.5,1);
	\node[above,inline text,red] at (-.5,1) {#2};
	\node[above,inline text,red] at (.5,1) {#4};
}

\newcommand{\annparamss}[4]
{
	\def\ta{#1};
	\def\tb{#2};
	\def\tap{#3};
	\def\tbp{#4};
}
\newcommand{\annss}[2]{
	\def\ra{.5};
	\def\rb{1.5};
	\ifthenelse{\equal{#1}{}}{}
		{
			\node[left,inline text,blue] at ($(0, 0) + (180:.75 and .75)$() {\footnotesize#1};
			\draw[blue,xscale=2.5] (0,-.75)to[out=90+45,in=270-45] (0,.75);
		}
	\ifthenelse{\equal{#2}{}}{}
	{
		\node[right,inline text,blue] at ($(0, 0) + (0:.75 and .75)$() {\footnotesize#2};
		\draw[blue,xscale=1.5] (0,-1.25)to[out=90-45,in=270+45] (0,1.25);
	}
	\draw[red,ultra thick] (0,\ra)--(0,\rb);
	\draw[red,ultra thick] (0,-\ra)--(0,-\rb);
	\node[below,inline text,red] at ($(0, 0) + (270:1.5 and 1.5)$() {\footnotesize\ta};
	\node[above,inline text,red] at ($(0, 0) + (90:1.5 and 1.5)$() {\footnotesize\tb};
	\node[above,inline text,red] at ($(0, 0) + (270:.5 and .5)$() {\footnotesize\tap};
	\node[below,inline text,red] at ($(0, 0) + (90:.5 and .5)$() {\footnotesize\tbp};
	\draw (0,0) circle (\ra);
	\draw (0,0) circle (\rb);
}
\newcommand{\annparamst}[4]
{
	\def\ta{#1};
	\def\tb{#2};
	\def\tap{#3};
	\def\tbp{#4};
}
\newcommand{\annst}[2]{
	\def\ra{.5};
	\def\rb{1.5};
	\ifthenelse{\equal{#1}{}}{}
	{
		\node[left,inline text,blue] at ($(0, 0) + (180:.75 and .75)$() {\footnotesize#1};
		\draw[blue,xscale=2.5] (0,-.75)to[out=90+45,in=270-45] (0,.75);
	}
	\ifthenelse{\equal{#2}{}}{}
	{
		\node[right,inline text,blue] at ($(0, 0) + (0:.75 and .75)$() {\footnotesize#2};
		\draw[blue,xscale=1.5] (0,-1.25)to[out=90-45,in=270+45] (0,1.25);
	}
	\draw[orange,ultra thick] (0,\ra)--(0,\rb);
	\draw[red,ultra thick] (0,-\ra)--(0,-\rb);
	\node[below,inline text,red] at ($(0, 0) + (270:1.5 and 1.5)$() {\footnotesize\ta};
	\node[above,inline text,orange] at ($(0, 0) + (90:1.5 and 1.5)$() {\footnotesize\tb};
	\node[above,inline text,red] at ($(0, 0) + (270:.5 and .5)$() {\footnotesize\tap};
	\node[below,inline text,orange] at ($(0, 0) + (90:.5 and .5)$() {\footnotesize\tbp};
	\draw (0,0) circle (\ra);
	\draw (0,0) circle (\rb);
}
\newcommand{\annparamstpq}[8]
{
	\def\ta{#1};
	\def\tb{#2};
	\def\tc{#3};
	\def\td{#4};
	\def\tap{#5};
	\def\tbp{#6};
	\def\tcp{#7};
	\def\tdp{#8};
}
\newcommand{\annstpq}[4]{
	\ifthenelse{\equal{#1}{}}{}
	{
		\node[left,inline text,blue] at ($(180:.6 and .75)$() {\footnotesize#1};
		\coordinate (A) at ({.5*cos(-140)},{1*sin(140)});
		\draw[blue] ({.5*cos(-140)},{1*sin(-140)}) to[out=90+45,in=270-45] (A);
	}
	\ifthenelse{\equal{#2}{}}{}
	{
		\node[above,inline text,blue] at ($(90:.5 and .9)$() {\footnotesize#2};
		\coordinate (B) at ({.5*cos(40)},{1.7*sin(40)});
		\draw[blue] ({.5*cos(140)},{1.2*sin(140)}) to[out=45,in=270-45] (B);
	}
	\ifthenelse{\equal{#3}{}}{}
	{
		\node[right,inline text,blue] at ($(0,.3)+(0:.75 and .75)$() {\footnotesize#3};
		\coordinate (C) at ({.5*cos(40)},{2*sin(40)});
		\draw[blue] ({.5*cos(40)},{1*sin(-40)}) to[out=90-45,in=270+45] (C);
	}
	\ifthenelse{\equal{#4}{}}{}
	{
		\node[above,inline text,blue] at ($(270:.5 and 1.1)$() {\footnotesize#4};
		\coordinate (A) at ({.5*cos(-40)},{1.5*sin(-40)});
		\draw[blue] ({.5*cos(-140)},{2*sin(-140)}) to[out=45,in=270-45] (A);
	}
	\draw[red,ultra thick] ({.5*cos(-140)},{.5*sin(-140)+.1})--({.5*cos(-140)},{-1.5*sin(acos(1/3*cos(-140)))});
	\draw[orange,ultra thick] ({.5*cos(-40)},{.5*sin(-40)+.1})--({.5*cos(-40)},{-1.5*sin(acos(1/3*cos(-40)))});
	\draw[darkgreen,ultra thick] ({.5*cos(140)},{.5*sin(140)-.1})--({.5*cos(140)},{1.5*sin(acos(1/3*cos(140)))});
	\draw[darkred,ultra thick] ({.5*cos(40)},{.5*sin(40)-.1})--({.5*cos(40)},{1.5*sin(acos(1/3*cos(40)))});
\draw[fill=white] (0,0) ellipse (.5 and .5);
\draw (0,0) circle (1.5 and 1.5);
	\node[below,inline text,red] at ({.5*cos(-140)},{-1.5*sin(acos(1/3*cos(-140)))}) {\footnotesize\ta};
	\node[below,inline text,orange] at ({.5*cos(-40)},{-1.5*sin(acos(1/3*cos(-40)))}) {\footnotesize\tb};
	\node[above,inline text,darkgreen] at ({.5*cos(140)},{1.5*sin(acos(1/3*cos(140)))}) {\footnotesize\tc};
	\node[above,inline text,darkred] at ({.5*cos(40)},{1.5*sin(acos(1/3*cos(40)))}) {\footnotesize\td};
}
\newcommand{\annsspq}[4]{
	\ifthenelse{\equal{#1}{}}{}
	{
		\node[left,inline text,blue] at ($(180:.6 and .75)$() {\footnotesize#1};
		\coordinate (A) at ({.5*cos(-140)},{1*sin(140)});
		\draw[blue] ({.5*cos(-140)},{1*sin(-140)}) to[out=90+45,in=270-45] (A);
	}
	\ifthenelse{\equal{#2}{}}{}
	{
		\node[above,inline text,blue] at ($(90:.5 and .9)$() {\footnotesize#2};
		\coordinate (B) at ({.5*cos(40)},{1.7*sin(40)});
		\draw[blue] ({.5*cos(140)},{1.2*sin(140)}) to[out=45,in=270-45] (B);
	}
	\ifthenelse{\equal{#3}{}}{}
	{
		\node[right,inline text,blue] at ($(0,.3)+(0:.75 and .75)$() {\footnotesize#3};
		\coordinate (C) at ({.5*cos(40)},{2*sin(40)});
		\draw[blue] ({.5*cos(40)},{1*sin(-40)}) to[out=90-45,in=270+45] (C);
	}
	\ifthenelse{\equal{#4}{}}{}
	{
		\node[above,inline text,blue] at ($(270:.5 and 1.1)$() {\footnotesize#4};
		\coordinate (A) at ({.5*cos(-40)},{1.5*sin(-40)});
		\draw[blue] ({.5*cos(-140)},{2*sin(-140)}) to[out=45,in=270-45] (A);
	}
	\draw[red,ultra thick] ({.5*cos(-140)},{.5*sin(-140)+.1})--({.5*cos(-140)},{-1.5*sin(acos(1/3*cos(-140)))});
	\draw[red,ultra thick] ({.5*cos(-40)},{.5*sin(-40)+.1})--({.5*cos(-40)},{-1.5*sin(acos(1/3*cos(-40)))});
	\draw[darkgreen,ultra thick] ({.5*cos(140)},{.5*sin(140)-.1})--({.5*cos(140)},{1.5*sin(acos(1/3*cos(140)))});
	\draw[darkred,ultra thick] ({.5*cos(40)},{.5*sin(40)-.1})--({.5*cos(40)},{1.5*sin(acos(1/3*cos(40)))});
	\draw[fill=white] (0,0) ellipse (.5 and .5);
	\draw (0,0) circle (1.5 and 1.5);
	\node[below,inline text,red] at ({.5*cos(-140)},{-1.5*sin(acos(1/3*cos(-140)))}) {\footnotesize\ta};
	\node[below,inline text,red] at ({.5*cos(-40)},{-1.5*sin(acos(1/3*cos(-40)))}) {\footnotesize\tb};
	\node[above,inline text,darkgreen] at ({.5*cos(140)},{1.5*sin(acos(1/3*cos(140)))}) {\footnotesize\tc};
	\node[above,inline text,darkred] at ({.5*cos(40)},{1.5*sin(acos(1/3*cos(40)))}) {\footnotesize\td};
}
\newcommand{\annstpp}[4]{
	\ifthenelse{\equal{#1}{}}{}
	{
		\node[left,inline text,blue] at ($(180:.6 and .75)$() {\footnotesize#1};
		\coordinate (A) at ({.5*cos(-140)},{1*sin(140)});
		\draw[blue] ({.5*cos(-140)},{1*sin(-140)}) to[out=90+45,in=270-45] (A);
	}
	\ifthenelse{\equal{#2}{}}{}
	{
		\node[above,inline text,blue] at ($(90:.5 and .9)$() {\footnotesize#2};
		\coordinate (B) at ({.5*cos(40)},{1.7*sin(40)});
		\draw[blue] ({.5*cos(140)},{1.2*sin(140)}) to[out=45,in=270-45] (B);
	}
	\ifthenelse{\equal{#3}{}}{}
	{
		\node[right,inline text,blue] at ($(0,.3)+(0:.75 and .75)$() {\footnotesize#3};
		\coordinate (C) at ({.5*cos(40)},{2*sin(40)});
		\draw[blue] ({.5*cos(40)},{1*sin(-40)}) to[out=90-45,in=270+45] (C);
	}
	\ifthenelse{\equal{#4}{}}{}
	{
		\node[above,inline text,blue] at ($(270:.5 and 1.1)$() {\footnotesize#4};
		\coordinate (A) at ({.5*cos(-40)},{1.5*sin(-40)});
		\draw[blue] ({.5*cos(-140)},{2*sin(-140)}) to[out=45,in=270-45] (A);
	}
	\draw[red,ultra thick] ({.5*cos(-140)},{.5*sin(-140)+.1})--({.5*cos(-140)},{-1.5*sin(acos(1/3*cos(-140)))});
	\draw[orange,ultra thick] ({.5*cos(-40)},{.5*sin(-40)+.1})--({.5*cos(-40)},{-1.5*sin(acos(1/3*cos(-40)))});
	\draw[darkgreen,ultra thick] ({.5*cos(140)},{.5*sin(140)-.1})--({.5*cos(140)},{1.5*sin(acos(1/3*cos(140)))});
	\draw[darkgreen,ultra thick] ({.5*cos(40)},{.5*sin(40)-.1})--({.5*cos(40)},{1.5*sin(acos(1/3*cos(40)))});
	\draw[fill=white] (0,0) ellipse (.5 and .5);
	\draw (0,0) circle (1.5 and 1.5);
	\node[below,inline text,red] at ({.5*cos(-140)},{-1.5*sin(acos(1/3*cos(-140)))}) {\footnotesize\ta};
	\node[below,inline text,orange] at ({.5*cos(-40)},{-1.5*sin(acos(1/3*cos(-40)))}) {\footnotesize\tb};
	\node[above,inline text,darkgreen] at ({.5*cos(140)},{1.5*sin(acos(1/3*cos(140)))}) {\footnotesize\tc};
	\node[above,inline text,darkgreen] at ({.5*cos(40)},{1.5*sin(acos(1/3*cos(40)))}) {\footnotesize\td};
}
\newcommand{\annsspp}[4]{
	\ifthenelse{\equal{#1}{}}{}
	{
		\node[left,inline text,blue] at ($(180:.6 and .75)$() {\footnotesize#1};
		\coordinate (A) at ({.5*cos(-140)},{1*sin(140)});
		\draw[blue] ({.5*cos(-140)},{1*sin(-140)}) to[out=90+45,in=270-45] (A);
	}
	\ifthenelse{\equal{#2}{}}{}
	{
		\node[above,inline text,blue] at ($(90:.5 and .9)$() {\footnotesize#2};
		\coordinate (B) at ({.5*cos(40)},{1.7*sin(40)});
		\draw[blue] ({.5*cos(140)},{1.2*sin(140)}) to[out=45,in=270-45] (B);
	}
	\ifthenelse{\equal{#3}{}}{}
	{
		\node[right,inline text,blue] at ($(0,.3)+(0:.75 and .75)$() {\footnotesize#3};
		\coordinate (C) at ({.5*cos(40)},{2*sin(40)});
		\draw[blue] ({.5*cos(40)},{1*sin(-40)}) to[out=90-45,in=270+45] (C);
	}
	\ifthenelse{\equal{#4}{}}{}
	{
		\node[above,inline text,blue] at ($(270:.5 and 1.1)$() {\footnotesize#4};
		\coordinate (A) at ({.5*cos(-40)},{1.5*sin(-40)});
		\draw[blue] ({.5*cos(-140)},{2*sin(-140)}) to[out=45,in=270-45] (A);
	}
	\draw[red,ultra thick] ({.5*cos(-140)},{.5*sin(-140)+.1})--({.5*cos(-140)},{-1.5*sin(acos(1/3*cos(-140)))});
	\draw[red,ultra thick] ({.5*cos(-40)},{.5*sin(-40)+.1})--({.5*cos(-40)},{-1.5*sin(acos(1/3*cos(-40)))});
	\draw[darkgreen,ultra thick] ({.5*cos(140)},{.5*sin(140)-.1})--({.5*cos(140)},{1.5*sin(acos(1/3*cos(140)))});
	\draw[darkgreen,ultra thick] ({.5*cos(40)},{.5*sin(40)-.1})--({.5*cos(40)},{1.5*sin(acos(1/3*cos(40)))});
	\draw[fill=white] (0,0) ellipse (.5 and .5);
	\draw (0,0) circle (1.5 and 1.5);
	\node[below,inline text,red] at ({.5*cos(-140)},{-1.5*sin(acos(1/3*cos(-140)))}) {\footnotesize\ta};
	\node[below,inline text,red] at ({.5*cos(-40)},{-1.5*sin(acos(1/3*cos(-40)))}) {\footnotesize\tb};
	\node[above,inline text,darkgreen] at ({.5*cos(140)},{1.5*sin(acos(1/3*cos(140)))}) {\footnotesize\tc};
	\node[above,inline text,darkgreen] at ({.5*cos(40)},{1.5*sin(acos(1/3*cos(40)))}) {\footnotesize\td};
}
\newcommand{\annssss}[4]{
	\ifthenelse{\equal{#1}{}}{}
	{
		\node[left,inline text,blue] at ($(180:.6 and .75)$() {\footnotesize#1};
		\coordinate (A) at ({.5*cos(-140)},{1*sin(140)});
		\draw[blue] ({.5*cos(-140)},{1*sin(-140)}) to[out=90+45,in=270-45] (A);
	}
	\ifthenelse{\equal{#2}{}}{}
	{
		\node[above,inline text,blue] at ($(90:.5 and .9)$() {\footnotesize#2};
		\coordinate (B) at ({.5*cos(40)},{1.7*sin(40)});
		\draw[blue] ({.5*cos(140)},{1.2*sin(140)}) to[out=45,in=270-45] (B);
	}
	\ifthenelse{\equal{#3}{}}{}
	{
		\node[right,inline text,blue] at ($(0,.3)+(0:.75 and .75)$() {\footnotesize#3};
		\coordinate (C) at ({.5*cos(40)},{2*sin(40)});
		\draw[blue] ({.5*cos(40)},{1*sin(-40)}) to[out=90-45,in=270+45] (C);
	}
	\ifthenelse{\equal{#4}{}}{}
	{
		\node[above,inline text,blue] at ($(270:.5 and 1.1)$() {\footnotesize#4};
		\coordinate (A) at ({.5*cos(-40)},{1.5*sin(-40)});
		\draw[blue] ({.5*cos(-140)},{2*sin(-140)}) to[out=45,in=270-45] (A);
	}
	\draw[red,ultra thick] ({.5*cos(-140)},{.5*sin(-140)+.1})--({.5*cos(-140)},{-1.5*sin(acos(1/3*cos(-140)))});
	\draw[red,ultra thick] ({.5*cos(-40)},{.5*sin(-40)+.1})--({.5*cos(-40)},{-1.5*sin(acos(1/3*cos(-40)))});
	\draw[red,ultra thick] ({.5*cos(140)},{.5*sin(140)-.1})--({.5*cos(140)},{1.5*sin(acos(1/3*cos(140)))});
	\draw[red,ultra thick] ({.5*cos(40)},{.5*sin(40)-.1})--({.5*cos(40)},{1.5*sin(acos(1/3*cos(40)))});
	\draw[fill=white] (0,0) ellipse (.5 and .5);
	\draw (0,0) circle (1.5 and 1.5);
	\node[below,inline text,red] at ({.5*cos(-140)},{-1.5*sin(acos(1/3*cos(-140)))}) {\footnotesize\ta};
	\node[below,inline text,red] at ({.5*cos(-40)},{-1.5*sin(acos(1/3*cos(-40)))}) {\footnotesize\tb};
	\node[above,inline text,red] at ({.5*cos(140)},{1.5*sin(acos(1/3*cos(140)))}) {\footnotesize\tc};
	\node[above,inline text,red] at ({.5*cos(40)},{1.5*sin(acos(1/3*cos(40)))}) {\footnotesize\td};
}
\newcommand{\pantsstpq}[5]
{
	\generalpantsstpq{#1}{#2}{#3}{#4}{}{#5}{}{};
}
\newcommand{\generalpantsstpq}[8]
{
	\pgfmathsetmacro{\ra}{.25};
	\pgfmathsetmacro{\sa}{.75};
	\pgfmathsetmacro{\rb}{2};
	\pgfmathsetmacro{\rc}{1.5};
	\coordinate (A) at ({-\sa},{-\rc*sin(acos(-\sa/\rb))});
	\coordinate (B) at ({\sa},{-\rc*sin(acos(-\sa/\rb))});
	\coordinate (C) at ({-\sa},{\rc*sin(acos(\sa/\rb))});
	\coordinate (D) at ({\sa},{\rc*sin(acos(\sa/\rb))});
		\ifthenelse{\equal{#1}{}}{}
		{
			\coordinate (A1) at ($(A)!.9!(-\sa,-\ra)$);
			\coordinate (B1) at ({-\sa-1.5*\ra},{0});
			\coordinate (C1) at ($(C)!.9!(-\sa,\ra)$);
			\node[left,inline text,blue] at (A1) {\footnotesize#1};
			\draw[blue] (A1) to[out=90+45,in=270] (B1) to[out=90,in=180+45] (C1);
		}
		\ifthenelse{\equal{#2}{}}{}
		{
			\coordinate (A1) at ($(A)!.6!(-\sa,-\ra)$);
			\coordinate (B1) at ({-\sa+2*\ra},{0});
			\coordinate (C1) at ($(C)!.75!(-\sa,\ra)$);
			\node[right,inline text,blue] at (A1) {\footnotesize#2};
			\draw[blue] (A1) to[out=90-45,in=270] (B1) to[out=90,in=0-45] (C1);
		}
		\ifthenelse{\equal{#3}{}}{}
		{
			\coordinate (A1) at ($(B)!.9!(\sa,-\ra)$);
			\coordinate (B1) at ({\sa-2*\ra},{0});
			\coordinate (C1) at ($(D)!.9!(\sa,\ra)$);
			\node[left,inline text,blue] at (A1) {\footnotesize#3};
			\draw[blue] (A1) to[out=90+45,in=270] (B1) to[out=90,in=180+45] (C1);
		}
		\ifthenelse{\equal{#4}{}}{}
		{
			\coordinate (A1) at ($(B)!.75!(\sa,-\ra)$);
			\coordinate (B1) at ({\sa+1.5*\ra},{0});
			\coordinate (C1) at ($(D)!.75!(\sa,\ra)$);
			\node[right,inline text,blue] at (A1) {\footnotesize#4};
			\draw[blue] (A1) to[out=90-45,in=270] (B1) to[out=90,in=0-45] (C1);
		}
		\ifthenelse{\equal{#5}{}}{}
		{
			\coordinate (A1) at ($(A)!.3!(-\sa,-\ra)$);
			\coordinate (B1) at ({-\sa-3*\ra},{0});
			\coordinate (C1) at ($(C)!.5!(-\sa,\ra)$);
			\node[left,inline text,blue] at (A1) {\footnotesize#5};
			\draw[blue] (A1) to[out=90+45,in=270] (B1) to[out=90,in=180+45] (C1);
		}
		\ifthenelse{\equal{#6}{}}{}
		{
			\coordinate (A1) at ($(C)!.25!(-\sa,-\ra)$);
			\coordinate (C1) at ($(D)!.25!(\sa,\ra)$);
			\coordinate (B1) at ($(A1)!.5!(C1)$);
			\node[below,inline text,blue] at (B1) {\footnotesize#6};
			\draw[blue] (A1) to[out=45,in=180+45] (C1);
		}
		\ifthenelse{\equal{#7}{}}{}
		{
			\coordinate (A1) at ($(B)!.3!(\sa,-\ra)$);
			\coordinate (B1) at ({\sa+3*\ra},{0});
			\coordinate (C1) at ($(D)!.1!(\sa,\ra)$);
			\node[right,inline text,blue] at (A1) {\footnotesize#7};
			\draw[blue] (A1) to[out=90-45,in=270] (B1) to[out=90,in=-45] (C1);
		}
		\ifthenelse{\equal{#8}{}}{}
		{
			\coordinate (A1) at ($(A)!.15!(-\sa,-\ra)$);
			\coordinate (C1) at ($(B)!.15!(\sa,\ra)$);
			\coordinate (B1) at ($(A1)!.5!(C1)$);
			\node[above,inline text,blue] at (B1) {\footnotesize#8};
			\draw[blue] (A1) to[out=45,in=180+45] (C1);
		}
\begin{scope}[even odd rule]
	\clip (0,0) ellipse [x radius=\rb,y radius=\rc] (-\sa,0) circle (\ra) (\sa,0) circle (\ra);
	\draw [red,ultra thick] ({-\sa},{0})--({-\sa},{-2*\rb});
	\draw [orange,ultra thick] ({\sa},{0})--({\sa},{-2*\rb});
	\draw [darkgreen,ultra thick] ({-\sa},{0})--({-\sa},{2*\rb});
	\draw [darkred,ultra thick] ({\sa},{0})--({\sa},{2*\rb});
\end{scope}
	\draw (-\sa,0) circle (\ra);
	\draw (\sa,0) circle (\ra);
	\draw (0,0) ellipse [x radius=\rb,y radius=\rc];
	\node[red,inline text,below] at (A) {\footnotesize\ta};
	\node[orange,inline text,below] at (B) {\footnotesize\tb};
	\node[darkgreen,inline text,above] at (C) {\footnotesize\tc};
	\node[darkred,inline text,above] at (D) {\footnotesize\td};
}

\newcommand{\tube}[3]{
	\ifthenelse{\equal{#2}{}}{}
	{
		\draw[domain=-90:0,smooth,variable=\x,blue] plot ({cos(\x)},{.1*sin(\x)+(\x)/360-.25});
		\draw[domain=0:180,smooth,variable=\x,blue,dashed] plot ({cos(\x)},{.1*sin(\x)+(\x)/360-.25});
		\draw[domain=180:270,smooth,variable=\x,blue] plot ({cos(\x)},{.1*sin(\x)+(\x)/360-.25});
	}
	\node[below,inline text,red] at ($(0, -1) + (270:1 and .1)$() {\footnotesize#1};
	\node[above,inline text,red] at ($(0, 1) + (-270:1 and .1)$() {\footnotesize#3};
	\draw (0,1) ellipse (1 and .1);
	\draw (-1,1)--(-1,-1) (1,1)--(1,-1);
	\draw[dashed] ($(0, -1) + (180:1 and .1)$() arc (180:0:1 and .1);
	\draw ($(0, -1) + (180:1 and .1)$() arc (180:360:1 and .1);
	\ifthenelse{\equal{#2}{}}{}
	{
		\node[inline text,blue] at (.5,-.3) {\footnotesize#2};
	}
	\draw [red,thick] ($(0, -1) + (270:1 and .1)$()--($(0, 1) + (-90:1 and .1)$();
}

\newcommand{\ttube}[6]{
	\draw [red,thick,dashed] ($(0, -1) + (70:1 and .1)$()--($(0, 1) + (70:1 and .1)$();
	\ifthenelse{\equal{#4}{}}{}
	{
		\draw[domain=-110:0,smooth,variable=\x,blue] plot ({cos(\x)},{.1*sin(\x)+(\x)/360-.25});
		\draw[domain=0:70,smooth,variable=\x,blue,dashed] plot ({cos(\x)},{.1*sin(\x)+(\x)/360-.25});
	}
	\ifthenelse{\equal{#3}{}}{}
	{
		\draw[domain=70:180,smooth,variable=\x,blue,dashed] plot ({cos(\x)},{.1*sin(\x)+(\x)/360});
		\draw[domain=180:250,smooth,variable=\x,blue] plot ({cos(\x)},{.1*sin(\x)+(\x)/360});
	}
	\node[below,inline text,red] at ($(0, -1) + (250:1 and .1)$() {\footnotesize#1};
	\node[below,inline text,red] at ($(0, -1) + (-70:1 and .1)$() {\footnotesize#2};
	\node[above,inline text,red] at ($(0, 1) + (-250:1 and .1)$() {\footnotesize#5};
	\node[above,inline text,red] at ($(0, 1) + (70:1 and .1)$() {\footnotesize#6};
	\draw (0,1) ellipse (1 and .1);
	\draw (-1,1)--(-1,-1) (1,1)--(1,-1);
	\draw[dashed] ($(0, -1) + (180:1 and .1)$() arc (180:0:1 and .1);
	\draw ($(0, -1) + (180:1 and .1)$() arc (180:360:1 and .1);
	\ifthenelse{\equal{#4}{}}{}
	{
		\node[right,inline text,blue] at ({cos(0)},{.1*sin(0)+(0)/360-.25}) {\footnotesize#4};
	}
	\ifthenelse{\equal{#3}{}}{}
	{
		\node[left,inline text,blue] at ({cos(180)},{.1*sin(180)+(180)/360}) {\footnotesize#3};
	}
	\draw [red,thick] ($(0, -1) + (250:1 and .1)$()--($(0, 1) + (250:1 and .1)$();
}

\newcommand{\ttttube}[4]{
	\draw [red,thick,dashed] ($(0, -1) + (100:1 and .1)$()--($(0, 1) + (100:1 and .1)$();
	\draw [red,thick,dashed] ($(0, -1) + (40:1 and .1)$()--($(0, 1) + (40:1 and .1)$();
	\draw (0,1) ellipse (1 and .1);
	\draw (-1,1)--(-1,-1) (1,1)--(1,-1);
	\draw[dashed] ($(0, -1) + (180:1 and .1)$() arc (180:0:1 and .1);
	\draw ($(0, -1) + (180:1 and .1)$() arc (180:360:1 and .1);
	\ifthenelse{\equal{#1}{}}{}
	{
		\draw[domain=230-360:280-360,smooth,variable=\x,blue] plot ({cos(\x)},{.1*sin(\x)+(\x)/360-.25});
	}
	\ifthenelse{\equal{#2}{}}{}
	{
		\draw[domain=280-360:0,smooth,variable=\x,blue] plot ({cos(\x)},{.1*sin(\x)+(\x)/360-.15});
		\draw[domain=0:40,smooth,variable=\x,blue,dashed] plot ({cos(\x)},{.1*sin(\x)+(\x)/360-.15});
	}
	\ifthenelse{\equal{#3}{}}{}
	{
		\draw[domain=40:100,smooth,variable=\x,blue,dashed] plot ({cos(\x)},{.1*sin(\x)+(\x)/360-.05});
	}
	\ifthenelse{\equal{#4}{}}{}
	{
		\draw[domain=100:180,smooth,variable=\x,blue,dashed] plot ({cos(\x)},{.1*sin(\x)+(\x)/360+.05});
		\draw[domain=180:230,smooth,variable=\x,blue] plot ({cos(\x)},{.1*sin(\x)+(\x)/360+.05});
	}
	\draw [red,thick] ($(0, -1) + (230:1 and .1)$()--($(0, 1) + (230:1 and .1)$();
	\draw [red,thick] ($(0, -1) + (280:1 and .1)$()--($(0, 1) + (280:1 and .1)$();
	\node[below,inline text,red] at ($(0, -1) + (230:1 and .1)$() {\footnotesize\ta};
	\node[below,inline text,red] at ($(0, -1) + (280:1 and .1)$() {\footnotesize\tb};
	\node[below,inline text,red] at ($(0, -1) + (-100:1 and .1)$() {\footnotesize\tc};
	\node[below,inline text,red] at ($(0, -1) + (-40:1 and .1)$() {\footnotesize\td};
	\node[above,inline text,red] at ($(0, 1) + (-230:1 and .1)$() {\footnotesize\tap};
	\node[above,inline text,red] at ($(0, 1) + (-280:1 and .1)$() {\footnotesize\tbp};
	\node[above,inline text,red] at ($(0, 1) + (100:1 and .1)$() {\footnotesize\tcp};
	\node[above,inline text,red] at ($(0, 1) + (40:1 and .1)$() {\footnotesize\tdp};
	\ifthenelse{\equal{#1}{}}{}
	{
		\node[above,inline text,blue] at ({cos(-90)},{.1*sin(-90)+(-90)/360-.25}) {\footnotesize#1};
	}
	\ifthenelse{\equal{#2}{}}{}
	{
		\node[right,inline text,blue] at ({cos(0)},{.1*sin(0)+(0)/360-.15}) {\footnotesize#2};
	}
	\ifthenelse{\equal{#3}{}}{}
	{
		\node[above,inline text,blue] at ({cos(60)},{.1*sin(60)+(60)/360-.05}) {\footnotesize#3};
	}
	\ifthenelse{\equal{#4}{}}{}
	{
		\node[left,inline text,blue] at ({cos(180)},{.1*sin(180)+(180)/360+.05}) {\footnotesize#4};
	}
}

\usetikzlibrary{intersections}

\tikzset{
	dot/.style={circle,inner sep=1pt,fill},
}

\newcommand{\pantsparams}[8]
{
	\def\ta{#1};
	\def\tb{#2};
	\def\tc{#3};
	\def\td{#4};
	\def\tap{#5};
	\def\tbp{#6};
	\def\tcp{#7};
	\def\tdp{#8};
}

\newcommand{\pants}[8]{
	\draw (0,2) ellipse (1.5 and .1);
	\draw[name path=left edge] (-2.5,-2)to[out=90,in=-90](-1.5,2);
	\draw[dashed] ($(-1.5, -2) + (180:1 and .1)$() arc (180:0:1 and .1);
	\draw ($(-1.5, -2) + (180:1 and .1)$() arc (180:360:1 and .1);
	\draw[name path=right edge] (2.5,-2)to[out=90,in=-90](1.5,2);
	\draw[dashed] ($(1.5, -2) + (180:1 and .1)$() arc (180:0:1 and .1);
	\draw ($(1.5, -2) + (180:1 and .1)$() arc (180:360:1 and .1);
	\draw[name path=centre] (-.5,-2)to[out=90,in=180](0,0)to[out=0,in=90](.5,-2);
	\draw[thick,red,dashed,name path=C] ($(-1.5, -2) + (70:1 and .1)$()to[out=90,in=-90]($(0, 2) + (100:1.5 and .1)$();
	\draw[thick,red,dashed,name path=D] ($(1.5, -2) + (70:1 and .1)$()to[out=90,in=-90]($(0, 2) + (40:1.5 and .1)$();
	\path[thick,red,name path=A] ($(-1.5, -2) + (250:1 and .1)$()to[out=90,in=-90]($(0, 2) + (230:1.5 and .1)$();
	\path[thick,red,name path=B] ($(1.5, -2) + (250:1 and .1)$()to[out=90,in=-90]($(0, 2) + (280:1.5 and .1)$();
	\ifthenelse{\equal{#1}{}}{}
	{
		\begin{scope}
			\path[name path=X] (-.5,-.25)--(-2.5,0);
			\path [name intersections={of=C and X,by={c1}}];
			\path [name intersections={of=X and left edge,by={c2}}];
			\draw[blue,dashed] (c1)to[out=140,in=70] (c2);
			\path[name path=Y] (c2) -- (0,-2);
			\path [name intersections={of=Y and A,by={d1}}];
			\draw[blue] (c2) to[out=-110,in=220]  (d1);
			\node[left,inline text,blue] at (c2) {\footnotesize#1};
		\end{scope}
	}
	\ifthenelse{\equal{#2}{}}{}
	{
		\begin{scope}
			\path[name path=X] (-2.5,-2)--(0,-.5);
			\path [name intersections={of=A and X,by={c1}}];
			\path [name intersections={of=X and centre,by={c2}}];
			\draw[blue] (c1) to[out=45,in=260]  (c2);
			\path[name path=Y] (c2) -- (-2.5,2);
			\path [name intersections={of=Y and C,by={d1}}];
			\draw[blue,dashed] (c2) to[out=95,in=-40]  (d1);
			\node[right,inline text,blue] at (c2) {\footnotesize#2};
		\end{scope}
	}
	\ifthenelse{\equal{#3}{}}{}
	{
		\begin{scope}
			\path[name path=X] (0,-.5)--(2.5,-.5);
			\path[name intersections={of=D and X,by={c1}}];
			\path[name intersections={of=X and centre,by={c2}}];
			\draw[blue,dashed] (c1)to[out=140,in=100] (c2);
			\path[name path=Y] (c2) -- (2,-2);
			\path [name intersections={of=Y and B,by={d1}}];
			\draw[blue] (c2) to[out=-80,in=220]  (d1);
			\node[left,inline text,blue] at (c2) {\footnotesize#3};
		\end{scope}
	}
	\ifthenelse{\equal{#4}{}}{}
	{
		\begin{scope}
			\path[name path=X] (0,-2.5)--(2.5,-1);
			\path [name intersections={of=B and X,by={c1}}];
			\path [name intersections={of=X and right edge,by={c2}}];
			\draw[blue] (c1) to[out=45,in=280]  (c2);
			\path[name path=Y] (c2) -- (0,.5);
			\path [name intersections={of=Y and D,by={d1}}];
			\draw[blue,dashed] (c2) to[out=110,in=-40]  (d1);
			\node[right,inline text,blue] at (c2) {\footnotesize#4};
		\end{scope}
	}
	\ifthenelse{\equal{#5}{}}{}
	{
		\begin{scope}
			\path[name path=X] (-2.5,2)--(0,.5);
			\path [name intersections={of=left edge and X,by={c1}}];
			\path [name intersections={of=X and A,by={c2}}];
			\draw[blue] (c1) to[out=-110,in=250]  (c2);
			\path[name path=Y] (c1) -- (0,1.5);
			\path [name intersections={of=Y and C,by={d1}}];
			\draw[blue,dashed] (c1) to[out=80,in=110]  (d1);
			\node[left,inline text,blue] at (c1) {\footnotesize#5};
		\end{scope}
	}
	\ifthenelse{\equal{#6}{}}{}
	{
		\begin{scope}
			\path[name path=X] (-2.5,.5)--(2.5,.5);
			\path [name intersections={of=A and X,by={c1}}];
			\path [name intersections={of=X and B,by={c2}}];
			\draw[blue] (c1) to[out=70,in=250]  (c2);
			\node[left,inline text,blue] at ($ (c2) + (0,.1) $) {\footnotesize#6};
		\end{scope}
	}
	\ifthenelse{\equal{#7}{}}{}
	{
		\begin{scope}
			\path[name path=X] (0,.5)--(2.5,1);
			\path [name intersections={of=B and X,by={c1}}];
			\path [name intersections={of=X and right edge,by={c2}}];
			\draw[blue] (c1) to[out=70,in=280]  (c2);
			\path[name path=Y] (c2) -- (0,2.2);
			\path [name intersections={of=Y and D,by={d1}}];
			\draw[blue,dashed] (c2) to[out=110,in=280]  (d1);
			\node[right,inline text,blue] at (c2) {\footnotesize#7};
		\end{scope}
	}
	\ifthenelse{\equal{#8}{}}{}
	{
		\begin{scope}
			\path[name path=X] (-2.5,1.3)--(2.5,1.3);
			\path [name intersections={of=D and X,by={c1}}];
			\path [name intersections={of=X and C,by={c2}}];
			\draw[blue,dashed] (c1) to[out=120,in=300]  (c2);
			\node[left,inline text,blue] at (c1) {\footnotesize#8};
		\end{scope}
	}
	\draw[thick,red] ($(-1.5, -2) + (250:1 and .1)$()to[out=90,in=-90]($(0, 2) + (230:1.5 and .1)$();
	\draw[thick,red] ($(1.5, -2) + (250:1 and .1)$()to[out=90,in=-90]($(0, 2) + (280:1.5 and .1)$();
	\node[below,inline text,red] at ($(-1.5, -2) + (250:1 and .1)$() {\footnotesize\ta};
	\node[below,inline text,red] at ($(1.5, -2) + (250:1 and .1)$() {\footnotesize\tb};
	\node[below,inline text,red] at ($(-1.5, -2) + (-70:1 and .1)$() {\footnotesize\tc};
	\node[below,inline text,red] at ($(1.5, -2) + (-70:1 and .1)$() {\footnotesize\td};
	\node[above,inline text,red] at ($(0, 2) + (-230:1.5 and .1)$() {\footnotesize\tap};
	\node[above,inline text,red] at ($(0, 2) + (-280:1.5 and .1)$() {\footnotesize\tbp};
	\node[above,inline text,red] at ($(0, 2) + (100:1.5 and .1)$() {\footnotesize\tcp};
	\node[above,inline text,red] at ($(0, 2) + (40:1.5 and .1)$() {\footnotesize\tdp};
}

\def\ladder[#1][#2][#3][#4][#5][#6]{
	\path(#1);
	\pgfgetlastxy{\XCoord}{\YCoord};
	\begin{scope}[xshift=\XCoord,yshift=\YCoord]
		\draw (-.6,-.5)--(-.6,.5);
		\draw (.6,-.5)--(.6,.5);
		\ifthenelse{\equal{#4}{}}{}
			{
				\draw (-.6,-.25)--(.6,.25);
			};
		\node[below,inline text] at (-.6,-.5) {\footnotesize#2};
		\node[below,inline text] at (.6,-.5) {\footnotesize#3};
		\node[below,inline text] at (0,0) {\footnotesize#4};
		\node[above,inline text] at (-.6,.5) {\footnotesize#5};
		\node[above,inline text] at (.6,.5) {\footnotesize#6};
	\end{scope}
}

\def\Lladder[#1][#2][#3][#4][#5][#6][#7]{
	\path(#1);
	\pgfgetlastxy{\XCoord}{\YCoord};
	\begin{scope}[xshift=\XCoord,yshift=\YCoord]
		\draw (-1.8,-.5)--(-1.8,.5);
		\draw (-.6,-.5)--(-.6,.5);
		\draw (.6,-.5)--(.6,.5);
		\ifthenelse{\equal{#4}{}}{}
		{
			\draw (-.6,-.25)--(.6,.25);
		};
		\node[below,inline text] at (-.6,-.5) {\footnotesize#2};
		\node[below,inline text] at (.6,-.5) {\footnotesize#3};
		\node[below,inline text] at (0,0) {\footnotesize#4};
		\node[above,inline text] at (-.6,.5) {\footnotesize#5};
		\node[above,inline text] at (.6,.5) {\footnotesize#6};
		\node[below,inline text] at (-1.8,-.5) {\footnotesize#7};
	\end{scope}
}

\def\Lmorphism[#1][#2][#3][#4][#5][#6][#7]{
	\path(#1);
	\pgfgetlastxy{\XCoord}{\YCoord};
	\begin{scope}[xshift=\XCoord,yshift=\YCoord]
		\draw (-1.8,-.5)to[out=90,in=190](-.6,.25);
		\draw (-.6,-.5)--(-.6,.5);
		\draw (.6,-.5)--(.6,.5);
		\ifthenelse{\equal{#4}{}}{}
		{
			\draw (-.6,-.25)--(.6,.25);
		};
		\node[below,inline text] at (-.6,-.5) {\footnotesize#2};
		\node[below,inline text] at (.6,-.5) {\footnotesize#3};
		\node[below,inline text] at (0,0) {\footnotesize#4};
		\node[above,inline text] at (-.6,.5) {\footnotesize#5};
		\node[above,inline text] at (.6,.5) {\footnotesize#6};
		\node[below,inline text] at (-1.8,-.5) {\footnotesize#7};
	\end{scope}
}

\def\Rladder[#1][#2][#3][#4][#5][#6][#7]{
	\path(#1);
	\pgfgetlastxy{\XCoord}{\YCoord};
	\begin{scope}[xshift=\XCoord,yshift=\YCoord]
		\draw (1.8,-.5)--(1.8,.5);
		\draw (-.6,-.5)--(-.6,.5);
		\draw (.6,-.5)--(.6,.5);
		\ifthenelse{\equal{#4}{}}{}
		{
			\draw (-.6,-.25)--(.6,.25);
		};
		\node[below,inline text] at (-.6,-.5) {\footnotesize#2};
		\node[below,inline text] at (.6,-.5) {\footnotesize#3};
		\node[below,inline text] at (0,0) {\footnotesize#4};
		\node[above,inline text] at (-.6,.5) {\footnotesize#5};
		\node[above,inline text] at (.6,.5) {\footnotesize#6};
		\node[below,inline text] at (1.8,-.5) {\footnotesize#7};
	\end{scope}
}

\def\Rmorphism[#1][#2][#3][#4][#5][#6][#7]{
	\path(#1);
	\pgfgetlastxy{\XCoord}{\YCoord};
	\begin{scope}[xshift=\XCoord,yshift=\YCoord]
		\draw (1.8,-.5)to[out=90,in=-10](.6,.3);
		\draw (-.6,-.5)--(-.6,.5);
		\draw (.6,-.5)--(.6,.5);
		\ifthenelse{\equal{#4}{}}{}
		{
			\draw (-.6,-.35)--(.6,.15);
		};
		\node[below,inline text] at (-.6,-.5) {\footnotesize#2};
		\node[below,inline text] at (.6,-.5) {\footnotesize#3};
		\node[below,inline text] at (0,-.1) {\footnotesize#4};
		\node[above,inline text] at (-.6,.5) {\footnotesize#5};
		\node[above,inline text] at (.6,.5) {\footnotesize#6};
		\node[below,inline text] at (1.8,-.5) {\footnotesize#7};
	\end{scope}
}

\def\Laction[#1][#2][#3]{
		\draw (0,-.5)--(0,.5);
		\draw (-.6,-.5)to[out=90,in=210](0,0);
		\node[below,inline text] at (-.6,-.5) {\footnotesize#1};
		\node[below,inline text] at (0,-.5) {\footnotesize#2};
		\node[above,inline text]  at (0,.5) {\footnotesize#3};
}
\def\Raction[#1][#2][#3]{
	\draw (0,-.5)--(0,.5);
	\draw (.6,-.5)to[out=90,in=-30](0,0);
	\node[below,inline text] at (.6,-.5) {\footnotesize#1};
	\node[below,inline text] at (0,-.5) {\footnotesize#2};
	\node[above,inline text] at (0,.5) {\footnotesize#3};
}

\def\Lassociator[#1][#2][#3][#4]{
	\draw (0,-.5)--(0,.5);
	\draw (-.6,-.5)to[out=90,in=210](0,0);
	\draw (.6,-.5)to[out=90,in=-30](0,.25);
	\node[below,inline text] at (-.6,-.5) {\footnotesize#1};
	\node[below,inline text] at (.6,-.5) {\footnotesize#2};
	\node[below,inline text] at (0,-.5) {\footnotesize#3};
	\node[above,inline text] at (0,.5) {\footnotesize#4};
}
\def\Rassociator[#1][#2][#3][#4]{
	\draw (0,-.5)--(0,.5);
	\draw (-.6,-.5)to[out=90,in=210](0,.25);
	\draw (.6,-.5)to[out=90,in=-30](0,0);
	\node[below,inline text] at (-.6,-.5) {\footnotesize#1};
	\node[below,inline text] at (.6,-.5) {\footnotesize#2};
	\node[below,inline text] at (0,-.5) {\footnotesize#3};
	\node[above,inline text] at (0,.5) {\footnotesize#4};
}

\usepackage{hhline}
\usepackage{etoolbox}

\let\originalleft\left
\let\originalright\right
\renewcommand{\left}{\mathopen{}\mathclose\bgroup\originalleft}
\renewcommand{\right}{\aftergroup\egroup\originalright}

\newcommand{\vvec}[1]
{
	\operatorname{\bf Vec}
	\ifstrequal{#1}{}
	{}
	{\left(#1\right)}
}

\newcommand{\vvectwist}[2]
{
	\operatorname{\bf Vec}^{#2}
	\ifstrequal{#1}{}
	{}
	{\left(#1\right)}
}

\newcommand{\lad}[1]
{
	\operatorname{\bf Lad}
	\ifstrequal{#1}{}
	{}
	{\left(#1\right)}
}

\newcommand{\ann}[3]
{
	\operatorname{\bf Ann}_{\mathcal{#1},\mathcal{#2}}
	\ifstrequal{#3}{}
	{}
	{\left(\mathcal{#3}\right)}
}

\newcommand{\kar}[1]
{
	\operatorname{\bf Kar}
	\ifstrequal{#1}{}
	{}
	{\left(#1\right)}
}

\newcommand{\bpr}[1]
{
	\operatorname{\bf BPR}
	\ifstrequal{#1}{}
	{}
	{\left(#1\right)}
}

\newcommand{\dih}[1]
{
	\operatorname{Dih}
	\ifstrequal{#1}{}
	{}
	{_{#1}}
}

\newcommand{\defect}[5]{
	\ifthenelse{\equal{#3}{}}
	{
		\begin{smallmatrix} #2\hfill	\\ #1\hfill \end{smallmatrix}\!\bigr\vert^{#5}
	}
	{
		\ifthenelse{\equal{#4}{}}
		{
			\begin{smallmatrix} #2\hfill	\\ #1\hfill \end{smallmatrix}\!\bigr\vert_{#3}^{#5}
		}
		{
			\begin{smallmatrix} #2\hfill	\\ #1\hfill \end{smallmatrix}\!\bigr\vert_{(#3,#4)}^{#5}
		}
	}
}

\newcommand{\vpantsparams}[3]
{
	\def\ta{#1};
	\def\tb{#2};
	\def\tc{#3};
}
\newcommand{\vgeneralpantsstp}[6]
{
	\ifthenelse{\equal{#1}{}}{}
	{
		\coordinate (A1) at (0,-1);
		\coordinate (B1) at (-.3,-.6);
		\coordinate (C1) at (0,-.2);
		\node[left,inline text,blue] at (C1) {\footnotesize#1};
		\draw[blue] (A1) to[out=90+45,in=270] (B1) to[out=90,in=180+45] (C1);
	}
	\ifthenelse{\equal{#2}{}}{}
	{
		\coordinate (A1) at (0,-1.1);
		\coordinate (B1) at (.3,-.6);
		\coordinate (C1) at (0,-.1);
		\node[right,inline text,blue] at (C1) {\footnotesize#2};
		\draw[blue] (A1) to[out=90-45,in=270] (B1) to[out=90,in=0-45] (C1);
	}
	\ifthenelse{\equal{#3}{}}{}
	{
		\coordinate (A1) at (0,.2);
		\coordinate (B1) at (-.3,.6);
		\coordinate (C1) at (0,1);
		\node[left,inline text,blue] at (A1) {\footnotesize#3};
		\draw[blue] (A1) to[out=90+45,in=270] (B1) to[out=90,in=180+45] (C1);
	}
	\ifthenelse{\equal{#4}{}}{}
	{
		\coordinate (A1) at (0,.1);
		\coordinate (B1) at (.3,.6);
		\coordinate (C1) at (0,1.1);
		\node[right,inline text,blue] at (A1) {\footnotesize#4};
		\draw[blue] (A1) to[out=90-45,in=270] (B1) to[out=90,in=0-45] (C1);
	}
	\ifthenelse{\equal{#5}{}}{}
	{
		\coordinate (A1) at (0,-1.3);
		\coordinate (B1) at (-.9,0);
		\coordinate (C1) at (0,1.3);
		\node[left,inline text,blue] at (B1) {\footnotesize#5};
		\draw[blue] (A1) to[out=90+45,in=270] (B1) to[out=90,in=180+45] (C1);
	}
	\ifthenelse{\equal{#6}{}}{}
	{
		\coordinate (A1) at (0,-1.4);
		\coordinate (B1) at (.9,0);
		\coordinate (C1) at (0,1.4);
		\node[right,inline text,blue] at (B1) {\footnotesize#6};
		\draw[blue] (A1) to[out=90-45,in=270] (B1) to[out=90,in=0-45] (C1);
	}
	\begin{scope}[even odd rule]
		\clip (0,0) ellipse [x radius=1,y radius=1.5] (0,-.6) circle (.2) (0,.6) circle (.2);
		\draw [red,ultra thick] (0,-1.5)--(0,-.6);
		\draw [orange,ultra thick] (0,-.6)--(0,.6);
		\draw [darkgreen,ultra thick] (0,.6)--(0,1.5);
	\end{scope}
	\draw (0,-.6) circle (.2);
	\draw (0,.6) circle (.2);
	\draw (0,0) ellipse [x radius=1.5,y radius=1.5];
	\node[red,inline text,below] at (0,-1.5) {\footnotesize\ta};
	\node[orange,inline text,above] at (0,.4) {\footnotesize\tb};
	\node[darkgreen,inline text,above] at (0,1.5) {\footnotesize\tc};
}

\newcommand{\vpantsstp}[4]
{
	\ifthenelse{\equal{#1}{}}{}
	{
		\coordinate (A1) at (0,-1);
		\coordinate (B1) at (-.3,-.6);
		\coordinate (C1) at (0,-.2);
		\node[left,inline text,blue] at (B1) {\footnotesize#1};
		\draw[blue] (A1) to[out=90+45,in=270] (B1) to[out=90,in=180+45] (C1);
	}
	\ifthenelse{\equal{#2}{}}{}
	{
		\coordinate (A1) at (0,-1.1);
		\coordinate (B1) at (.3,-.6);
		\coordinate (C1) at (0,-.1);
		\node[right,inline text,blue] at (B1) {\footnotesize#2};
		\draw[blue] (A1) to[out=90-45,in=270] (B1) to[out=90,in=0-45] (C1);
	}
	\ifthenelse{\equal{#3}{}}{}
	{
		\coordinate (A1) at (0,.2);
		\coordinate (B1) at (-.3,.6);
		\coordinate (C1) at (0,1);
		\node[left,inline text,blue] at (B1) {\footnotesize#3};
		\draw[blue] (A1) to[out=90+45,in=270] (B1) to[out=90,in=180+45] (C1);
	}
	\ifthenelse{\equal{#4}{}}{}
	{
		\coordinate (A1) at (0,.1);
		\coordinate (B1) at (.3,.6);
		\coordinate (C1) at (0,1.1);
		\node[right,inline text,blue] at (B1) {\footnotesize#4};
		\draw[blue] (A1) to[out=90-45,in=270] (B1) to[out=90,in=0-45] (C1);
	}
	\begin{scope}[even odd rule]
		\clip (0,0) ellipse [x radius=1,y radius=1.5] (0,-.6) circle (.2) (0,.6) circle (.2);
		\draw [red,ultra thick] (0,-1.5)--(0,-.6);
		\draw [orange,ultra thick] (0,-.6)--(0,.6);
		\draw [darkgreen,ultra thick] (0,.6)--(0,1.5);
	\end{scope}
	\draw (0,-.6) circle (.2);
	\draw (0,.6) circle (.2);
	\draw (0,0) ellipse [x radius=1.5,y radius=1.5];
	\node[red,inline text,below] at (0,-1.5) {\footnotesize\ta};
	\node[orange,inline text,above] at (0,.4) {\footnotesize\tb};
	\node[darkgreen,inline text,above] at (0,1.5) {\footnotesize\tc};
}

\usepackage{soul}

\allowdisplaybreaks
\begin{document}

\title{Fusing Binary Interface Defects in Topological Phases: \\The $\mathbb{Z}/p\mathbb{Z}$ case}

\author{Jacob C.\ Bridgeman}
\email{jcbridgeman1@gmail.com}
\affiliation{Centre for Quantum Software and Information, Faculty of Engineering and Information Technology, University of Technology Sydney, Australia}
\author{Daniel Barter}
\email{danielbarter@gmail.com}
\affiliation{Mathematical Sciences Institute, Australian National University, Canberra, Australia}
\author{Corey Jones}
\email{cormjones88@gmail.com} 
\affiliation{Department of Mathematics, The Ohio State University, USA}

\date{\today}

\begin{abstract}
A binary interface defect is any interface between two (not necessarily invertible) domain walls. We compute all possible binary interface defects in Kitaev's $\mathbb{Z}/p\mathbb{Z}$ model and all possible fusions between them. Our methods can be applied to any Levin-Wen model. We also give physical interpretations for each of the defects in the $\mathbb{Z}/p\mathbb{Z}$ model. These physical interpretations provide a new graphical calculus which can be used to compute defect fusion. 
\end{abstract}

\maketitle



Topological phases are promising candidate materials for robust encoding, storage and manipulation of quantum information\cite{Dennis2002,MR1951039,Brown2014,Terhal2015}. Formed by locally interacting degrees of freedom, these quantum systems have emergent global properties that are protected against the presence of environmental noise. In addition to the bulk properties, inclusion of defects has been shown to improve the power of these materials from a quantum computational perspective\cite{Dennis2002,0610082,Bombin2007a,Bombin2010,Brown2013a,Pastawski2014,Yoshida2015a,1606.07116,Brown2016,IrisCong1,IrisCong2,PhysRevB.96.195129,Yoshida2017,SETPaper,Brown2018}. It is therefore important to understand the full theory, including defects. In this paper, we study non-chiral, 2-dimensional, long-range-entangled topological phases with defects.

A defect of a topological phase is a region of positive codimension which differs from the ground state. For example, in a 2-dimensional topological phase, domain walls are codimension 1 defects and anyonic excitations are codimension 2 defects. Much work has been done on defects in topological phases, for example \onlinecites{Dennis2002,0610082,Bombin2007a,Bombin2010,MR2942952,FUCHS2002353,MR3370609,Kong2013,Brown2013a,MR3063919,Barkeshli2013,Barkeshli2014,Pastawski2014,Yoshida2015a,1606.07116,Brown2016,IrisCong1,IrisCong2,PhysRevB.96.195129,Yoshida2017,PhysRevB.96.125104,SETPaper,Bridgeman2017,Brown2018,1809.00245}. In previous work, the term defect frequently refers to a 0-dimensional interface between two invertible domain walls. To avoid confusion, we shall use the term {\em binary interface defect} to refer to \emph{any} 0-dimensional interface between two, not necessarily invertible, domain walls.

This work builds on our previous paper \onlinecite{1806.01279}, in which we computed the fusions of all domain walls in Kitaev's $\vvec{\ZZ{p}}$ models\cite{MR1951039} (with $p$ prime). In this paper, we compute all possible fusions between all possible binary interface defects in the Kitaev $\vvec{\ZZ{p}}$ model. The tools from both \onlinecite{1806.01279} and the present work can be adapted to more general Levin-Wen models. In the physics literature, defect fusion is often synonymous with symmetry gauging. We compute the fusions even when no gauging exists. 

In \onlinecite{MR1951039}, Kitaev defined a 2D lattice model associated to any finite group $G$ with particle-like low energy excitations (known as \emph{anyons}) parameterized by the simple representations of the Drinfeld double of $G$. These Kitaev models are some of the most well known examples of topological phases, and are of great experimental interest\cite{chow2014implementing,Gambetta1}.

These models were generalized in \onlinecite{Levin2005} to allow any fusion category $\mathcal{C}$ as input. The low energy excitations of these \emph{Levin-Wen} models are particle-like and parameterized by simple object of the Drinfeld center $Z(\mathcal{C})$. When $\mathcal{C} = \vvec{G}$, the Levin-Wen model can be transformed into the Kitaev model associated to $G$ with a finite depth quantum circuit, so they describe the same phase of matter. Indeed, the category of representations of the Drinfeld double of $G$ is equivalent to the Drinfeld center of $\vvec{G}$.

There are many interesting examples of fusion categories, so the Levin-Wen construction gives us many interesting 2D lattice models. Unfortunately they are too complicated to simulate or study using conventional lattice quantum field theoretic techniques. Moreover, for most interesting fusion categories, the data required to write down the Hamiltonian in the Levin-Wen construction is not known explicitly. For this reason, we need alternative tools to study these models.

The renormalization invariant properties of a topological phase is described by a topological quantum field theory (TQFT). In mathematics, a TQFT is a functor from a bordism category into a linear algebraic category. It would be counter-productive for us to give a precise definition here. The reader interested in details is encouraged to consult the recent survey \onlinecite{MR3674995} which includes complete definitions. In \onlinecite{MR1357878}, Barrett and Westbury described how to construct a (2+1)D TQFT from a fusion category, generalizing the Turaev-Viro construction from \onlinecite{MR1191386}. It is well understood that the TQFT associated to a fusion category captures the renormalization invariant properties of the corresponding Levin-Wen model. In this paper, we use TQFTs as studied in mathematics to compute renormalization invariant properties of Levin-Wen models which are intractable from the lattice theoretic perspective.

Our work builds on the work of Morrison and Walker in \onlinecite{MR2978449}. If $\mathcal{C}$ is a fusion category, morphisms in $\mathcal{C}$ can be described by 2-dimensional string diagrams up to isotopy. If $\Sigma$ is a 2-manifold with boundary, the TQFT associated to $\mathcal{C}$ sends $\Sigma$ to the vector space of string diagrams from $\mathcal{C}$ drawn on $\Sigma$ modulo local relations. These vector spaces are called Skein modules. The string diagrams are allowed to terminate on boundary components. When $\Sigma$ has a boundary, the vector space is graded by the object labels on the boundary components, and the graded pieces often assemble into an algebraic object. The process of drawing string diagrams from $\mathcal{C}$ on $\Sigma$ also extends to modules over $\mathcal{C}$. Morphisms in these module categories can also be described using string diagrams, and drawing string diagrams from module categories on $\Sigma$ allows us to use diagrammatic techniques to study defects of codimension 1 and 2 in the corresponding topological phase. The TQFT assigns a category to the annulus, which can be additionally decorated with bimodule strings
\begin{align}
	\begin{array}{c}
		\includeTikz{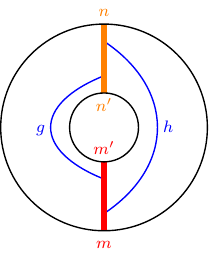}{
			\begin{tikzpicture}[scale=.7,,every node/.style={scale=.6}]
				\annparamst{$m$}{$n$}{$m'$}{$n'$};
				\annst{$g$}{$h$};
			\end{tikzpicture}
		}
	\end{array} : (m,n) \to (m',n').\label{eqn:generic_morphism_in_tube_category}
\end{align}
The representations of this category classify codimension 2 defects (binary interface defects in our language). We refer the reader to the recent survey \onlinecite{1607.05747} for more details. 
In \onlinecite{MR2978449}, representations of this annular category are called sphere modules. Indeed, when $\mathcal{C} = \vvec{G}$ and both of the modules are $\vvec{G}$, then this category is Morita equivalent to the Drinfeld double of $G$. The representations (defects) then correspond to the anyonic excitations of the model. In \onlinecite{MR2942952}, Kitaev and Kong explain that fusion category bimodules correspond to domain walls and bimodule functors correspond to codimension 2 defects in the Levin-Wen models. In \onlinecite{SETPaper}, annular categories, called \emph{dube algebras} by Williamson, Bultinck and Verstraete, are used to study defects interfacing invertible domain walls.

Many of the defect fusions which we compute have multiplicity, which correspond to multiple fusion channels. These multiplicities are somewhat mysterious from the physical perspective.

From a mathematical perspective, we are computing decomposition rules for relative tensor product and composition of bimodule functors. In an upcoming paper \footnote{D.~Barter, J.~C.~Bridgeman, C.~Jones, \emph{in preparation}}, we will provide a rigorous proof of this fact using a robust theory of skeletalization of fusion categories and their bimodules.

\subsection*{What is being computed in this work}

Suppose that $\mathcal{A},\,\mathcal{B}$ are 2-dimensional phases of matter and $\mathcal{M},\,\mathcal{N}$ are domain walls between the two phases. A defect $\alpha$ interfaces between two domain walls
\begin{align}
\begin{array}{c}
\includeTikz{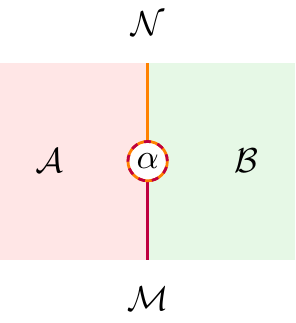}
{
	\begin{tikzpicture}
	\fill[red!10](-1.5,-1) rectangle (0,1);
	\fill[darkgreen!10](1.5,-1) rectangle (0,1);
	\draw[thick,purple] (0,-1)--(0,0);
	\draw[thick,orange](0,0)--(0,1);
	\draw[thick,orange,fill=white] (0,0) ellipse (.2 and .2);
	\draw[thick,purple,dash pattern = on 2 off 2,dash phase=0] (0,0) ellipse (.2 and .2);
	\node[] at (0,0) {$\alpha$};
	\node at (-1,0) {$\mathcal{A}$};
	\node at (1,0) {$\mathcal{B}$};
  \node at (0,1.4) {$\mathcal{N}$};
  \node at (0,-1.4) {$\mathcal{M}$};
	\end{tikzpicture}	
}
\end{array}.
\end{align}
In this paper, we are interested in defects when $\mathcal{A}$ and $\mathcal{B}$ are the phase associated to the fusion category $\vvec{\ZZ{p}}$, and all possible ways these defects can be fused. 
There are two ways in which defects can fuse:
\begin{itemize}
	\item Horizontally, where the supporting domain walls also fuse
	\begin{align} \label{fig:horizontal_defect_fusion}
		\begin{array}{c}
			\includeTikz{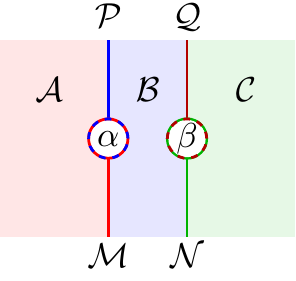}
			{
				\begin{tikzpicture}
				\fill[red!10](-1.5,-1) rectangle (-.4,1);
				\fill[blue!10](-.4,-1) rectangle (.4,1);
				\fill[darkgreen!10](1.5,-1) rectangle (.4,1);
				\draw[thick,red] (-.4,-1)--(-.4,0);
				\draw[thick,blue](-.4,0)--(-.4,1);
				\draw[thick,red,fill=white] (-.4,0) ellipse (.2 and .2);
				\node[fill=white,circle,inner sep=0pt] at (-.4,0) {$\alpha$};
				\draw[thick,blue,dash pattern = on 2 off 2,dash phase=0] (-.4,0) circle (.2);
				\draw[thick,darkgreen] (.4,-1)--(.4,0);
				\draw[thick,darkred](.4,0)--(.4,1);
				\draw[thick,darkgreen,fill=white] (.4,0) circle (.2);
				\node[fill=white,circle,inner sep=0pt] at (.4,0) {$\beta$};
				\draw[thick,darkred,dash pattern = on 2 off 2,dash phase=0] (.4,0) circle (.2);
				\node at (-1,.5) {$\mathcal{A}$};
				\node at (0,.5) {$\mathcal{B}$};
				\node at (1,.5) {$\mathcal{C}$};
				\node[inline text,below] at (-.4,-1) {$\mathcal{M}$};
				\node[inline text,below] at (.4,-1) {$\mathcal{N}$};
				\node[inline text,above] at (-.4,1) {$\mathcal{P}$};
				\node[inline text,above] at (.4,1) {$\mathcal{Q}$};
				\end{tikzpicture}	
			}
		\end{array}
		\to
		\begin{array}{c}
			\includeTikz{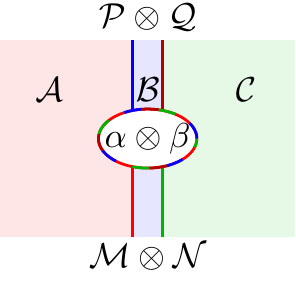}
			{
				\begin{tikzpicture}
				\fill[red!10](-1.5,-1) rectangle (-.15,1);
				\fill[blue!10](-.15,-1) rectangle (.15,1);
				\fill[darkgreen!10](1.5,-1) rectangle (.15,1);
				\draw[thick,red] (-.15,-1)--(-.15,0);
				\draw[thick,blue](-.15,0)--(-.15,1);
				\draw[thick,darkgreen] (.15,-1)--(.15,0);
				\draw[thick,darkred](.15,0)--(.15,1);
				\draw[thick,red,fill=white] (0,0) ellipse (.5 and .3);
				\draw[thick,blue,dash pattern = on 5 off 15] (0,0) ellipse (.5 and .3);
				\draw[thick,darkgreen,dash pattern = on 5 off 15,dash phase=10] (0,0) ellipse (.5 and .3);
				\draw[thick,darkred,dash pattern = on 5 off 15,dash phase=5] (0,0) ellipse (.5 and .3);
				\node[] at (0,0) {$\alpha\otimes\beta$};
				\node at (-1,.5) {$\mathcal{A}$};
				\node at (0,.5) {$\mathcal{B}$};
				\node at (1,.5) {$\mathcal{C}$};
				\node[inline text,below] at (0,-1) {$\mathcal{M}\otimes\mathcal{N}$};
				\node[inline text,above] at (0,1) {$\mathcal{P}\otimes\mathcal{Q}$};
				\end{tikzpicture}	
			}
		\end{array}
		\to
		\sum_{\mathcal{R},\mathcal{S},\gamma}
		\begin{array}{c}
			\includeTikz{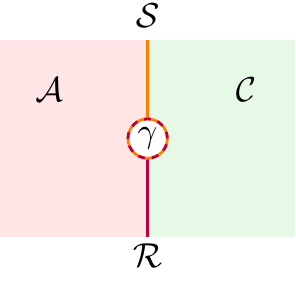}
			{
				\begin{tikzpicture}
				\fill[red!10](-1.5,-1) rectangle (0,1);
				\fill[darkgreen!10](1.5,-1) rectangle (0,1);
				\draw[thick,purple] (0,-1)--(0,0);
				\draw[thick,orange](0,0)--(0,1);
				\draw[thick,orange,fill=white] (0,0) ellipse (.2 and .2);
				\draw[thick,purple,dash pattern = on 2 off 2,dash phase=0] (0,0) ellipse (.2 and .2);
				\node[] at (0,0) {$\gamma$};
				\node at (-1,.5) {$\mathcal{A}$};
				\node at (1,.5) {$\mathcal{C}$};
				\node[inline text,below] at (0,-1) {$\mathcal{R}$};
				\node[inline text,above] at (0,1) {$\mathcal{S}$};
				\end{tikzpicture}	
			}
		\end{array}.
	\end{align}
	\\
\item Vertically, where the defects fuse along a common domain wall
\begin{align}
	\begin{array}{c} \label{fig:vertical_defect_fusion}
		\includeTikz{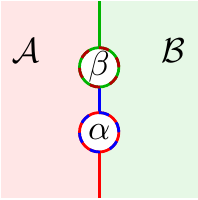}
		{
			\begin{tikzpicture}
			\fill[red!10](-1,-1) rectangle (0,1);
			\fill[darkgreen!10](0,-1) rectangle (1,1);
			\draw[thick,red] (0,-1)--(0,-.33);
			\draw[thick,blue](0,-.33)--(0,.33);
			\draw[thick,darkgreen](0,.33)--(0,1);
			\draw[thick,red,fill=white] (0,-.33) circle (.2);
			\node at (0,-.33) {$\alpha$};
			\draw[thick,blue,dashed] (0,-.33) circle (.2);
			\draw[thick,darkgreen,fill=white] (0,.33) circle (.2);
			\node at (0,.33) {$\beta$};
			\draw[thick,darkred,dashed] (0,.33) circle (.2);
			\node at (-.75,.5) {$\mathcal{A}$};
			\node at (.75,.5) {$\mathcal{B}$};
			\end{tikzpicture}	
		}
	\end{array}
	\to
	\begin{array}{c}
		\includeTikz{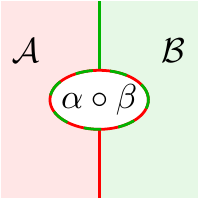}
		{
			\begin{tikzpicture}
			\fill[red!10](-1,-1) rectangle (0,1);
			\fill[darkgreen!10](0,-1) rectangle (1,1);
			\draw[thick,red] (0,-1)--(0,0);
			\draw[thick,darkgreen] (0,0)--(0,1);
			\draw[thick,red,fill=white] (0,0) ellipse (.5 and .3);
			\draw[thick,darkgreen,dash pattern = on 5 off 5] (0,0) ellipse (.5 and .3);
			\node[] at (0,0) {$\alpha\circ\beta$};
			\node at (-.75,.5) {$\mathcal{A}$};
			\node at (.75,.5) {$\mathcal{B}$};
			\end{tikzpicture}	
		}
	\end{array}
	\to
	\sum_{\gamma}
	\begin{array}{c}
		\includeTikz{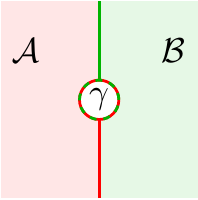}
		{
			\begin{tikzpicture}
			\fill[red!10](-1,-1) rectangle (0,1);
			\fill[darkgreen!10](0,-1) rectangle (1,1);
			\draw[thick,red] (0,-1)--(0,0);
			\draw[thick,darkgreen] (0,0)--(0,1);
			\draw[thick,red,fill=white] (0,0) ellipse (.2 and .2);
			\draw[thick,darkgreen,dash pattern = on 2 off 2,dash phase=0] (0,0) ellipse (.2 and .2);
			\node[] at (0,0) {$\gamma$};
			\node at (-.75,.5) {$\mathcal{A}$};
			\node at (.75,.5) {$\mathcal{B}$};
			\end{tikzpicture}	
		}
	\end{array}.
\end{align}
\end{itemize}
In this work, we present algorithms that can be used to compute both the horizontal and vertical fusions between all binary interface defects.
\subsection*{Structure of the paper}
The paper is structured as follows. In Section~\ref{sec:constructing_idempotents}, we describe how to compute idempotent representations of all the binary interface defects in Kitaev's $\ZZ{p}$ model. For a semi-simple category, the Karoubi envelope agrees with the category of representations, so these idempotents parameterize representations of the annulus categories. In Section~\ref{sec:hotizontal_defect_fusion_algorithm}, we describe the procedure used to compute the horizontal fusion of binary interface defects. Rather than give a formal algorithm, we describe how to proceed in a sufficiently general example. In Section~\ref{sec:horizontal_defect_fusion}, we demonstrate some of the horizontal fusion computations to elucidate some of the complications. In Section~\ref{sec:vertical_fusion_algorithm}, we outline the procedure used to compute vertical fusion. In Section~\ref{sec:vertical_defect_fusion}, we include some example vertical fusion computations. In Section~\ref{sec:physical_interpretations}, we give physical interpretations for all the binary interface defects. The physical interpretations can be used to reproduce all the horizontal and vertical fusion tables, except for the multiplicities, which are still somewhat mysterious from the physical perspective. In Section~\ref{sec:natural_transformations}, we explain how natural transformations between bimodule categories fit into our framework. This will be expanded on in future work.

In Appendix~\ref{sec:idempotents_table}, we tabulate the idempotent representations of all the binary interface defects in Kitaev's $\ZZ{p}$ model. In Appendix~\ref{sec:inflations}, we tabulate all the inflation data required to compute the horizontal fusions. This data was computed as an intermediate step in \onlinecite{1806.01279}. In Appendix~\ref{sec:horizontal_fusion_table}, we tabulate the horizontal fusions and in Appendix~\ref{sec:vertical_fusion_table}, we tabulate the vertical fusions. 


\section{Classifying defects} \label{sec:constructing_idempotents}

For an underlying fusion category $\mathcal{C}$, and given a pair of bimodule categories $\mathcal{M},\mathcal{N}$, the annular category $\ann{M}{N}{C}$ is defined as follows.
The objects are pairs of simples $(m,n)$ and the morphisms are annular diagrams as shown in Eqn.~\ref{eqn:generic_morphism_in_tube_category}.
Representations of $\ann{M}{N}{C}$ classify binary interface defects, which we denote $\defect{\mathcal{M}}{\mathcal{N}}{}{}{}$. If $\mathcal{C}$ is semi-simple, the category of representations is equivalent to the Karoubi envelope $\kar{\ann{M}{N}{C}}$. As in \onlinecite{1806.01279}, we utilize this equivalence to classify the defects.

Objects of $\kar{\ann{M}{N}{C}}$ are pairs $(A,e)$, where $A$ is an object from $\ann{M}{N}{C}$, and $e:A\to A$ is an idempotent annular diagram. The classification of defects is equivalent to construction of inequivalent idempotents. Two idempotent annular morphisms are equivalent if there exists a morphism absorbing the first idempotent on the inside and the second idempotent on the outside.

In the remainder of this section, we show how representative idempotents are constructed for $\mathcal{C}=\vvec{\ZZ{p}}$. The domain wall labels used are those defined in \onlinecite{1806.01279}.

\begin{exmp}[$\defect{T}{T}{}{}{}$]

In many cases, the construction of representative idempotents is a simple task. For example, the basic annulus diagrams for $\defect{T}{T}{}{}{}$ defects are
\begin{align}
	\begin{array}{c}
		\includeTikz{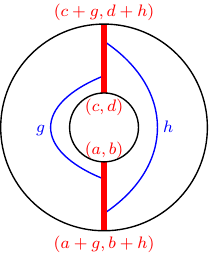}{
			\begin{tikzpicture}[scale=.7,,every node/.style={scale=.6}]
				\annparamss{$(a+g,b+h)$}{$(c+g,d+h)$}{$(a,b)$}{$(c,d)$};
				\annss{$g$}{$h$};
			\end{tikzpicture}
		}
	\end{array}.\label{eqn:TTtube_example_1}
\end{align}
For such a diagram to contribute to an idempotent, the inner and outer bimodule labels must be the same. In this case, the only way for this to happen is $g=h=0$ (the group identity)
\begin{align}
\begin{array}{c}
\includeTikz{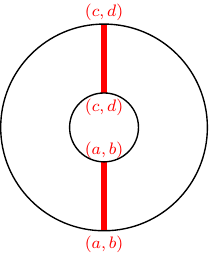}{
	\begin{tikzpicture}[scale=.7,,every node/.style={scale=.6}]
	\annparamss{$(a,b)$}{$(c,d)$}{$(a,b)$}{$(c,d)$};
	\annss{}{};
	\end{tikzpicture}
}
\end{array}.
\end{align}
It remains to find representatives of the isomorphism classes, where any annulus in Eqn.~\ref{eqn:TTtube_example_1} defines an isomorphism between idempotents. Using Eqn.~\ref{eqn:TTtube_example_1}, there is no way to change the value of $\alpha=c-a$ or $\beta=d-b$. These conserved quantities label the isomorphism classes, and we can pick representatives for each of the $p^2$ classes
\begin{align}
\begin{array}{c}
\includeTikz{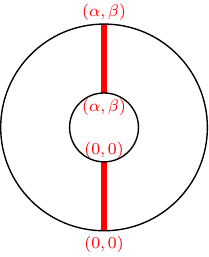}{
	\begin{tikzpicture}[scale=.7,,every node/.style={scale=.6}]
	\annparamss{$(0,0)$}{$(\alpha,\beta)$}{$(0,0)$}{$(\alpha,\beta)$};
	\annss{}{};
	\end{tikzpicture}
}
\end{array},
\end{align}
with $\alpha,\beta\in \ZZ{p}$.

\end{exmp}

\begin{exmp}[$\defect{X_k}{F_r}{}{}{}$: An unusual representation of $\ZZ{p}$]

A more complicated example involves $\defect{X_k}{F_r}{}{}{}$ defects. In this case, the general annulus diagrams are
\begin{align}
\begin{array}{c}
\includeTikz{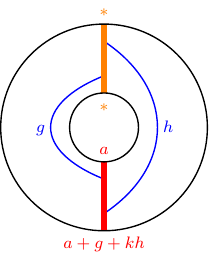}{
	\begin{tikzpicture}[scale=.7,,every node/.style={scale=.6}]
	\annparamst{$a+g+kh$}{$*$}{$a$}{$*$};
	\annst{$g$}{$h$};
	\end{tikzpicture}
}
\end{array},\label{eqn:XkFrtube_example_1}
\end{align}
so the diagrams contributing to idempotents are
\begin{align}
M_g&:=
\begin{array}{c}
\includeTikz{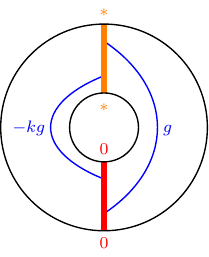}{
	\begin{tikzpicture}[scale=.7,,every node/.style={scale=.6}]
	\annparamst{$0$}{$*$}{$0$}{$*$};
	\annst{$-kg$}{$g$};
	\end{tikzpicture}
}
\end{array},\label{eqn:XkFrtube_example_2}
\end{align}
where $a=0$ has been chosen using the isomorphism Eqn.~\ref{eqn:XkFrtube_example_1}.
Due to the nontrivial associator on $F_r$, the multiplication rule for these annuli is
\begin{align}
M_gM_h&=
\begin{array}{c}
\includeTikz{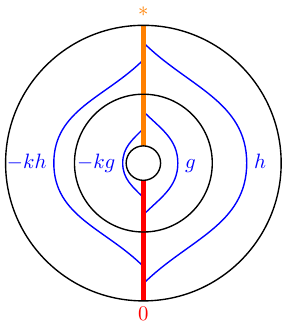}{
	\begin{tikzpicture}[scale=.7,,every node/.style={scale=.6}]
		\draw[blue] (0,-.75) to[out=45,in=-90] (.5,0) to[out=90,in=-45](0,.75);\node[blue,right,inline text] at (.5,0) {$g$};
		\draw[blue] (0,-.5) to[out=180-45,in=-90] (-.3,0) to[out=90,in=180+45](0,.5);\node[blue,left,inline text] at (-.3,0) {$-kg$};
		\draw[blue] (0,-1.75) to[out=45,in=-90] (1.5,0) to[out=90,in=-45](0,1.75);\node[blue,right,inline text] at (1.5,0) {$h$};
		\draw[blue] (0,-1.5) to[out=180-45,in=-90] (-1.3,0) to[out=90,in=180+45](0,1.5);\node[blue,left,inline text] at (-1.3,0) {$-kh$};
		\draw[ultra thick,orange] (0,.25)--(0,2) node[above,inline text] {$*$};
		\draw[ultra thick,red] (0,-.25)--(0,-2) node[below,inline text] {$0$};
		\draw (0,0) circle (.25);
		\draw (0,0) circle (1);
		\draw (0,0) circle (2);
	\end{tikzpicture}
}
\end{array}=
\omega^{kghr}
\begin{array}{c}
\includeTikz{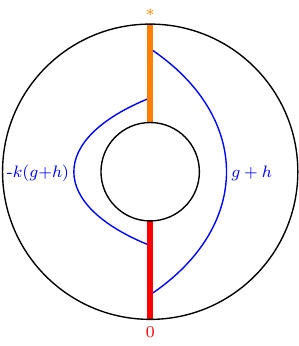}{
	\begin{tikzpicture}[scale=1,,every node/.style={scale=.6}]
	\annparamst{$0$}{$*$}{}{};
	\annst{-$k(g{+}h)$}{$g+h$};
	\end{tikzpicture}
}
\end{array}
=\omega^{kr gh}M_{g+h},
\end{align}
where $\omega=\exp(\frac{2\pi i}{p})$.
These annuli are therefore forming a twisted representation of $\ZZ{p}$ with 2-cocycle $\phi(g,h)=\omega^{krgh}$. Since $H^2(G,U(1))\cong \{1\}$, this is equivalent to a linear representation $U_g$ by some 1-cochain
\begin{align}
	\beta_{kr}(g)M_g&=U_g.
\end{align}
One can obtain an explicit formula for $\beta_{kr}$
\begin{align}
\beta_{kr}(g)&=
\begin{cases}
i^{g}&p=2\\
\omega^{kr g^2 2^{-1}}&p>2
\end{cases},
\end{align}
where $2^{-1}$ is the multiplicative inverse modulo $p$. With this explicit cochain, representatives for the inequivalent idempotents can be found
\begin{align}
	\defect{X_k}{F_r}{x}{}{}&=\frac{1}{p}\sum_g \omega^{gx} \omega U_g\\
	&=\frac{1}{p}\sum_g \omega^{gx}\beta(g) M_g=:\frac{1}{p}\sum_g \Theta_{x,kr}(g) M_g.
\end{align}
The equation $\Theta_{x,a}(g+k) = \Theta_{x,a}(g) \; \Theta_{x,a}(k) \; \omega^{agk}$ implies that the idempotent $\defect{X_k}{F_r}{x}{}{}$ absorbs the diagram
\begin{align}
\begin{array}{c}
				\includeTikz{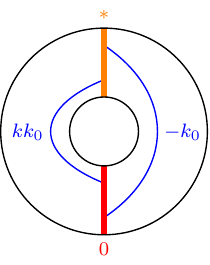}
				{
					\begin{tikzpicture}[scale=.7,every node/.style={scale=.7}]
					\annparamst{$0$}{$*$}{}{};
					\annst{$k k_0$}{$-k_0$};
					\end{tikzpicture}
				}
				\end{array}
\end{align}
up to a global phase.

\end{exmp}

\begin{exmp}[$\defect{F_q}{F_r}{}{}{}$: A genuine projective representation]

A particularly interesting example is given by $\defect{F_q}{F_r}{}{}{}$ defects. Since there is a single object (denoted $*$) in the $F_x$ bimodule, all annulus diagrams potentially contribute to the idempotents
\begin{align}
M_{g,h}&:=
\begin{array}{c}
\includeTikz{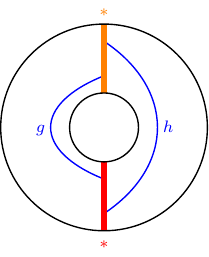}{
	\begin{tikzpicture}[scale=.7,,every node/.style={scale=.6}]
	\annparamst{$*$}{$*$}{}{};
	\annst{$g$}{$h$};
	\end{tikzpicture}
}
\end{array}.
\end{align}
Both modules (potentially) have nontrivial associators, leading to the multiplication rule
\begin{align}
M_{g_0,h_0}M_{g_1,h_1}&=
\begin{array}{c}
\includeTikz{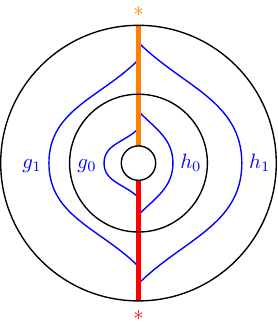}{
	\begin{tikzpicture}[scale=.7,,every node/.style={scale=.6}]
	\draw[blue] (0,-.75) to[out=45,in=-90] (.5,0) to[out=90,in=-45](0,.75);\node[blue,right,inline text] at (.5,0) {$h_0$};
	\draw[blue] (0,-.5) to[out=180-45,in=-90] (-.5,0) to[out=90,in=180+45](0,.5);\node[blue,left,inline text] at (-.5,0) {$g_0$};
	\draw[blue] (0,-1.75) to[out=45,in=-90] (1.5,0) to[out=90,in=-45](0,1.75);\node[blue,right,inline text] at (1.5,0) {$h_1$};
	\draw[blue] (0,-1.5) to[out=180-45,in=-90] (-1.3,0) to[out=90,in=180+45](0,1.5);\node[blue,left,inline text] at (-1.3,0) {$g_1$};
	\draw[ultra thick,orange] (0,.25)--(0,2) node[above,inline text] {$*$};
	\draw[ultra thick,red] (0,-.25)--(0,-2) node[below,inline text] {$*$};
	\draw (0,0) circle (.25);
	\draw (0,0) circle (1);
	\draw (0,0) circle (2);
	\end{tikzpicture}
}
\end{array}=
\omega^{(q-r)h_0g_1}
\begin{array}{c}
\includeTikz{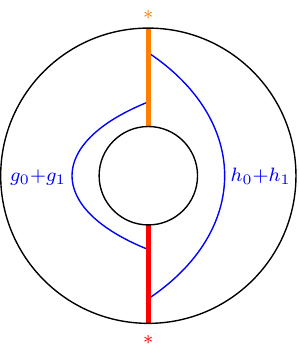}{
	\begin{tikzpicture}[scale=1,,every node/.style={scale=.7}]
	\annparamst{$*$}{$*$}{}{};
	\annst{$g_0{+}g_1$}{$h_0{+}h_1$};
	\end{tikzpicture}
}
\end{array}
=\omega^{(q-r)h_0g_1}M_{g_0+g_1,h_0+h_1}.
\end{align}
For $q\neq r$, this is a nontrivial 2-cocycle for $\ZZ{p}\times\ZZ{p}$\cite{propitius,DEWILDPROPITIUS1997297}, and therefore these annuli form a nontrivial twisted group algebra. There is an algebra isomorphism from this annulus algebra to the $p^2$ dimensional Pauli algebra. This isomorphism is defined by
\begin{align}
M_{g,h}\mapsto X^{(q-r)g}Z^h,\label{eqn:tubetopauli}
\end{align}
where $X$ and $Z$ are Pauli matrices obeying $ZX=\omega XZ$. The Pauli matrices span the full $p\times p$ matrix algebra. Up to isomorphism, there is a unique primitive idempotent of this algebra 
\begin{align}
P&=\begin{pmatrix}
1&0&\cdots&0\\
0&0&\cdots&0\\
\vdots&\vdots&\ddots&\vdots\\
0&0&\cdots&0
\end{pmatrix}\\&=\frac{1}{p}\sum_g Z^g.
\end{align}
Finally, we use the algebra isomorphism  Eqn.~\ref{eqn:tubetopauli} to obtain the idempotent for the single defect
\begin{align}
\defect{F_q}{F_r}{}{}{}&=\frac{1}{p}\sum_g
\begin{array}{c}
\includeTikz{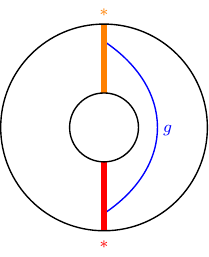}{
	\begin{tikzpicture}[scale=.7,,every node/.style={scale=.6}]
	\annparamst{$*$}{$*$}{}{};
	\annst{}{$g$};
	\end{tikzpicture}
}
\end{array}.
\end{align}

\end{exmp}

Other idempotents can be found using the techniques outlined here. A full set of representative idempotents is provided in Table~\ref{tab:idempotents}.


\section{The Horizontal Fusion Algorithm} \label{sec:hotizontal_defect_fusion_algorithm}

In this section, we explain the algorithm used to compute horizontal defect fusion. Due to the diagrammatic nature of the algorithm and the variety of phenomena that can occur during the computation, we are not going to give a formal specification, suitable for computer implementation. Instead, we shall explain how to proceed by hand in a specific example. The procedure can vary between examples and we demonstrate this in Section~\ref{sec:horizontal_defect_fusion}. This variability only occurs after all the necessary information has been extracted from the correct tables, so for the purpose of explaining how to navigate the tables of data which appear in this paper, a specific example is easier to understand than a general algorithm.

The algorithm proceeds in four key steps:
\begin{enumerate} 
	\item Determine idempotents corresponding to source defects using Table~\ref{tab:idempotents}.
	\item Determine idempotent for target defect using Tables~\ref{tab:zptable} and \ref{tab:idempotents}.
	\item Inflate the target idempotent to a 4-string annulus using Tables~\ref{tab:inflation_1}-\ref{tab:inflation_2}.
	\item Find a nonzero pants diagram that absorbs the source idempotents on the legs and (inflated) target on the waist.
\end{enumerate}

\begin{exmp}[$\defect{X_k}{T}{}{}{} \otimes \defect{F_s}{R}{}{}{}$]
We shall give the step by step procedure for the fusion
\begin{align}
\defect{X_k}{T}{a}{}{} \otimes \defect{F_s}{R}{z}{}{}.
\end{align}
\subsubsection{Writing down the source defect idempotents}
The first step in the procedure is to look up the idempotents representing the defects. These are found in Table~\ref{tab:idempotents}. For our current example, we have
\begin{align}
\defect{X_k}{T}{a}{}{}=
				\begin{array}{c}
				\includeTikz{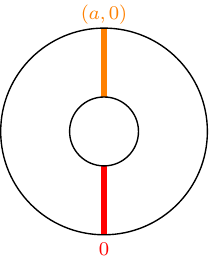}
				{
					\begin{tikzpicture}[scale=.7,every node/.style={scale=.7}]
					\annparamst{$0$}{$(a,0)$}{}{};
					\annst{}{};
					\end{tikzpicture}
				}
        \end{array}, \quad
\defect{F_s}{R}{z}{}{}=\frac{1}{p}\sum_g\omega^{gz}
				\begin{array}{c}
				\includeTikz{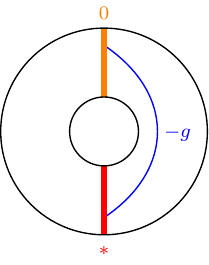}
				{
					\begin{tikzpicture}[scale=.7,every node/.style={scale=.7}]
					\annparamst{$*$}{$0$}{}{};
					\annst{}{$-g$};
					\end{tikzpicture}
				}
				\end{array}.
        \end{align}
\subsubsection{Writing down the target defect idempotent}
To decide the defect species resulting from the fusion, we need to look up the relevant domain wall fusions. These are contained in Table~\ref{tab:zptable}. For our current example, we have
\begin{align}
	T \otimes_{\vvec{\ZZ{p}}} R &= p \cdot R  \label{eq:top_domain_wall_fusion}\\
	X_k \otimes_{\vvec{\ZZ{p}}} F_s &= F_{k^{-1}s}. \label{eq:bottom_domain_wall_fusion}
\end{align}

From Eqns.~\eqref{eq:top_domain_wall_fusion} and \eqref{eq:bottom_domain_wall_fusion}, we can read of the target defect up to the labels. In this case it is $\defect{F_{k^{-1}s}}{R}{\zeta}{}{}$, where $\zeta$ is yet to be determined. From Table~\ref{tab:idempotents}, we find that the defect idempotent is
\begin{align}
	\defect{F_{k^{-1}s}}{R}{\zeta}{}{}=\frac{1}{p}\sum_{\gamma}\omega^{\gamma \zeta}
	\begin{array}{c}
		\includeTikz{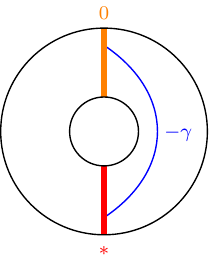}
		{
			\begin{tikzpicture}[scale=.7,every node/.style={scale=.7}]
			\annparamst{$*$}{$0$}{}{};
			\annst{}{$-\gamma$};
			\end{tikzpicture}
		}
	\end{array},
\end{align}
At this stage of the computation, the multiplicities in the domain wall fusion don't play a role. They show up in the inflation stage.
\subsubsection{Inflating the target idempotent}
In order for the target defect idempotent and the initial defect idempotents to interact with each other on a pair-of-pants, we need to \emph{inflate} the target idempotent so that it has four vertical strings. The information required to do this is contained in Table~\ref{tab:inflation_1}-\ref{tab:inflation_2}. These tables contain the information required to explicitly decompose the tensor products $\mathcal{M}\otimes_{\vvec{\ZZ{p}}}\mathcal{N}$ into simple bimodule categories. In our example, we need the entries corresponding to $R \to T \otimes_{\vvec{\ZZ{p}}} R$ and $F_{k^{-1}s} \to X_k \otimes_{\vvec{\ZZ{p}}} F_s$. These entries tell how to replace the trivalent vertices in our target idempotent, giving us
\begin{align}\frac{1}{p} \sum_{\gamma} \omega^{\gamma \zeta}
\begin{array}{c}
\includeTikz{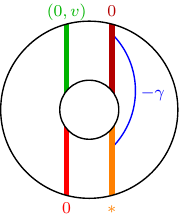}{
	\begin{tikzpicture}[scale=.6,,every node/.style={scale=.6}]
	\annparamstpq{$0$}{$*$}{$(0,v)$}{$0$}{}{}{}{};
	\annstpq{}{}{$-\gamma$}{};
	\end{tikzpicture}
}
\end{array}
\end{align}
The $\nu$ appearing corresponds to the multiplicity in the tensor product $T \otimes_{\vvec{\ZZ{p}}} R = p \cdot R$.

\subsubsection{Decorating the pants}
At this point, we have extracted all the data we need from the tables. Now we need to find all the pants diagrams which absorb our source defect idempotents on the legs and our target inflated defect idempotent at the waist. The most general pair-of-pants which we need to consider looks like:
\begin{align}
\begin{array}{c}
\includeTikz{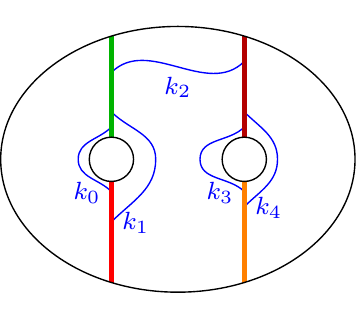}{
	\begin{tikzpicture}[scale=.9,,every node/.style={scale=.9}]
	\pantsparams{}{}{}{}{}{}{}{};
	\pantsstpq{$k_0$}{$k_1$}{$k_3$}{$k_4$}{$k_2$};
	\end{tikzpicture}}
\end{array}.
\end{align}
We use the following equation to bring every horizontal pair-of-pants into the standard form
\begin{align}
\begin{array}{c}
\includeTikz{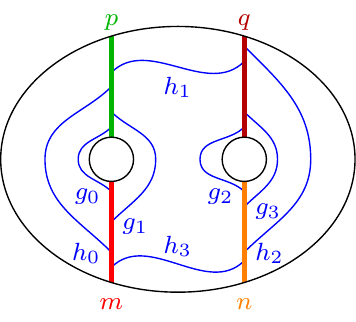}{
	\begin{tikzpicture}[scale=.9,,every node/.style={scale=.9}]
		\pantsparams{$m$}{$n$}{$p$}{$q$}{}{}{}{};
		\generalpantsstpq{$g_0$}{$g_1$}{$g_2$}{$g_3$}{$h_0$}{$h_1$}{$h_2$}{$h_3$};
	\end{tikzpicture}}
\end{array}
&=\frac{\Omega_M(h_0^{-1},g_1^{-1})\Omega_Q(h_1,h_2)\Omega_N(h_3,(g_3h_2)^{-1})\Omega_Q(h_3,g_3h_2)}{\Omega_P(h_0,g_1)}
\begin{array}{c}
\includeTikz{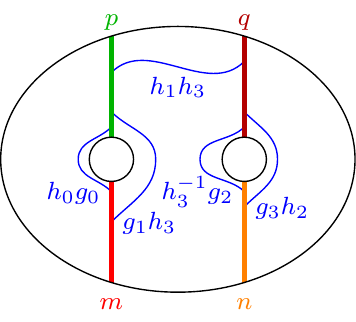}{
	\begin{tikzpicture}[scale=.9,,every node/.style={scale=.9}]
	\pantsparams{$m$}{$n$}{$p$}{$q$}{}{}{}{};
	\pantsstpq{$h_0g_0$}{$g_1h_3$}{$h_3^{-1}g_2$}{$g_3h_2$}{$h_1h_3$};
	\end{tikzpicture}}
\end{array},
\end{align}
where $\Omega_{X}(\bullet,\bullet)$ is the associator for the bimodule $X$\cite{1806.01279}.
Now, when we insert our source defect idempotents on the legs and our inflated target idempotent on the belt, we get
\begin{align}\frac{1}{p^2} \sum_{g,\gamma} \omega^{gz+\gamma\zeta+s\gamma k_3}
\begin{array}{c}
\includeTikz{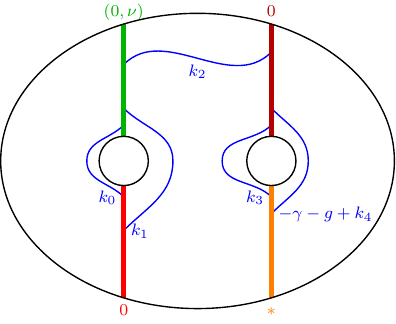}{
	\begin{tikzpicture}[scale=1,,every node/.style={scale=.6}]
		\pantsparams{$0$}{$*$}{$(0,\nu)$}{$0$}{}{}{}{};
		\generalpantsstpq{$k_0$}{$k_1$}{$k_3$}{$-\gamma-g+k_4$}{}{$k_2$}{}{};
	\end{tikzpicture}}
\end{array}
\end{align}
up to a global phase. The term $\omega^{s\gamma k_3}$ appears because we need to use the middle associator on $F_s$ to bring the diagram into the standard form. We have omitted the labels on the leg holes to make the diagram less cluttered. They match the labels on the source defect idempotents. To make all the labels match up, we must have
\begin{align*}
  k_0 &= -a \\
  k_1 &= k^{-1}a \\
  k_2 &= k^{-1}a - \nu \\
  k_3 &= \nu - k^{-1}a
  \end{align*}
The transformation $g \to g + k_4$ only changes the expression by a phase. The whole expression is zero unless $z = \zeta + sk_3$. Rearranging this equation gives
\begin{align}
  \zeta = z + s (k^{-1}a - v).
\end{align}
So we have
\begin{align}
  \defect{X_k}{T}{a}{}{} \otimes \defect{F_s}{R}{z}{}{} = \defect{F_{k^{-1}s}}{R}{z+s(k^{-1}a-\nu)}{}{\nu}
\end{align}
As explained above, the superscript $\nu$ indexes the multiplicity in the top domain wall fusion.
\end{exmp}

\section{Horizontal Defect Fusion} \label{sec:horizontal_defect_fusion}

In this section, we present several horizontal defect fusion computations to elucidate some of the complications that arise.

\begin{exmp}[$\defect{T}{T}{}{}{} \otimes \defect{T}{T}{}{}{}$]

Consider the fusion of two $T$-$T$ defects $\defect{T}{T}{a}{b}{}$ and $\defect{T}{T}{c}{d}{}$. This case is interesting because there is multiplicity in the domain wall fusion $T \otimes_{\ZZ{p}} T = p \cdot T$. The domain wall fusion and the defect fusion are correlated, therefore we represent this defect fusion as a $p \times p$-matrix indexed by the components in the decomposition $T \otimes_{\ZZ{p}} T = p \cdot T$: 
\begin{align}
\left[ \defect{T}{T}{a}{b}{}\otimes \defect{T}{T}{c}{d}{}\right]_{\mu,\nu} \to \defect{T}{T}{\alpha}{\beta}{}.
\end{align}
We want to establish the possible values of $\alpha$ and $\beta$ and the associated multiplicities. We do this by decomposing the tensor product into simple bimodule categories using $T\to T\otimes_{\vvec{\ZZ{p}}}$. By fixing $\mu$ and $\nu$, we can inflate $\defect{T}{T}{\alpha}{\beta}{}$ to a 4-string annulus
\begin{align}
\begin{array}{c}
\includeTikz{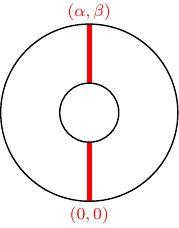}{
	\begin{tikzpicture}[scale=.6,,every node/.style={scale=.6}]
	\annparamss{$(0,0)$}{$(\alpha,\beta)$}{}{};
	\annss{}{};
	\end{tikzpicture}
}
\end{array}
\linf
\begin{array}{c}
\includeTikz{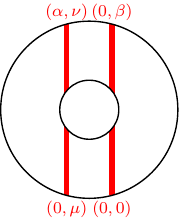}{
	\begin{tikzpicture}[scale=.6,,every node/.style={scale=.6}]
	\annparamstpq{$(0,\mu)$}{$(0,0)$}{$(\alpha,\nu)$}{$(0,\beta)$}{}{}{}{};
	\annssss{}{}{}{};
	\end{tikzpicture}
}
\end{array},
\end{align} 
Next, we can find a pant mapping the pair of 2-string defects to the 4-string defect. The most general form of the pants is
\begin{align}
\begin{array}{c}
\includeTikz{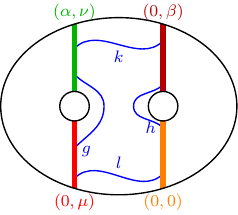}{
	\begin{tikzpicture}[scale=.6,,every node/.style={scale=.6}]
	\pantsparams{$(0,\mu)$}{$(0,0)$}{$(\alpha,\nu)$}{$(0,\beta)$}{}{}{}{};
	\generalpantsstpq{}{$g$}{$h$}{}{}{$k$}{}{$l$};
	\end{tikzpicture}
}
\end{array}
&=
\begin{array}{c}
\includeTikz{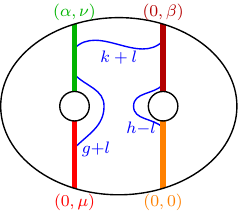}{
	\begin{tikzpicture}[scale=.6,,every node/.style={scale=.6}]
	\pantsparams{$(0,\mu)$}{$(0,0)$}{$(\alpha,\nu)$}{$(0,\beta)$}{}{}{}{};
	\pantsstpq{}{$g{+}l$}{$h{-}l$}{}{$k+l$};
	\end{tikzpicture}
}
\end{array}
=
\begin{array}{c}
\includeTikz{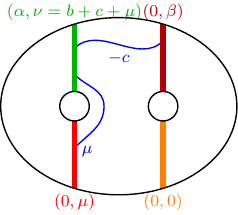}{
	\begin{tikzpicture}[scale=.6,,every node/.style={scale=.6}]
	\pantsparams{$(0,\mu)$}{$(0,0)$}{$(\alpha,\nu=b+c+\mu)$}{$(0,\beta)$}{}{}{}{};
	\pantsstpq{}{$\mu$}{}{}{$-c$};
	\end{tikzpicture}
}
\end{array}.
\end{align}
From this, we conclude that $\alpha = a$, $\beta = d$ and $\nu-\mu \equiv b + c \mod p$, which implies that
\begin{align}
\left[ \defect{T}{T}{a}{b}{}\otimes \defect{T}{T}{c}{d}{} \right]_{\mu,\nu} = \delta_{\nu-\mu}^{b+c} \cdot \defect{T}{T}{a}{d}{}.
\end{align}

\end{exmp}

\begin{exmp}[$\defect{X_k}{X_k}{}{}{} \otimes \defect{X_l}{X_l}{}{}{}$, which includes fusing `anyons']
	
The bimodule $X_1$ corresponds to $\vvec{\ZZ{p}}$ as a self-bimodule. As such, $\defect{X_1}{X_1}{}{}{}$ defects correspond to the excitations of the underlying topological phase. In the physics literature, these excitations are referred to as \emph{anyons}. Consider the fusion
\begin{align}
\defect{X_k}{X_k}{a}{x}{}\otimes \defect{X_l}{X_l}{b}{y}{}\to \defect{X_{kl}}{X_{kl}}{\alpha}{\beta}{}.
\end{align}

We can inflate $ \defect{X_{kl}}{X_{kl}}{\alpha}{\beta}{}$ to a 4-string annulus
\begin{align}
\frac{1}{p}\sum_g\omega^{g\beta}
\begin{array}{c}
\includeTikz{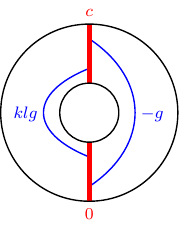}{
	\begin{tikzpicture}[scale=.6,,every node/.style={scale=.6}]
	\annparamss{$0$}{$c$}{}{};
	\annss{$klg$}{$-g$};
	\end{tikzpicture}
}
\end{array}
\linf\frac{1}{p}\sum_g\omega^{g\beta}
\begin{array}{c}
\includeTikz{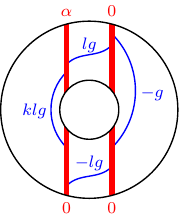}{
	\begin{tikzpicture}[scale=.6,,every node/.style={scale=.6}]
		\annparamstpq{$0$}{$0$}{$\alpha$}{$0$}{}{}{}{};
		\annssss{$klg$}{$lg$}{$-g$}{$-lg$};
	\end{tikzpicture}
}
\end{array}.
\end{align}

Next, we can find a pant mapping the 2-string defects to the 4-string defect. We immediately observe that (for nonzero maps) this forces $\alpha=a+kb$ 
\begin{align}
\frac{1}{p^3}\sum_{g,h,m}\omega^{gx+hy+m\beta}
\begin{array}{c}
\includeTikz{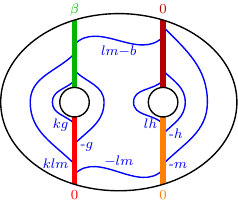}{
	\begin{tikzpicture}[scale=.6,,every node/.style={scale=.5}]
	\pantsparams{$0$}{$0$}{$\beta$}{$0$}{}{}{}{};
	\generalpantsstpq{$kg$}{-$g$}{$lh$}{-$h$}{$klm$}{$lm{-}b$}{-$m$}{$-lm$};
	\end{tikzpicture}
}
\end{array}
&=
\frac{1}{p^3}\sum_{g,h,m}\omega^{gx+hy+m\beta}
\begin{array}{c}
	\includeTikz{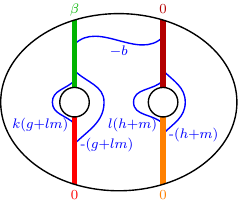}{
		\begin{tikzpicture}[scale=.6,,every node/.style={scale=.5}]
		\pantsparams{$0$}{$0$}{$\beta$}{$0$}{}{}{}{};
		\pantsstpq{$k(g{+}lm)$}{-$(g{+}lm)$}{$l(h{+}m)$}{-$(h{+}m)$}{$-b$};
		\end{tikzpicture}
	}
\end{array}.
\end{align}
Making the replacements
\begin{align}
g^\prime&=g+lm,h^\prime=h+m,
\end{align}
we obtain
\begin{align}
\frac{1}{p^2}\sum_{g,h}\omega^{gx+hy+m\beta}
\begin{array}{c}
\includeTikz{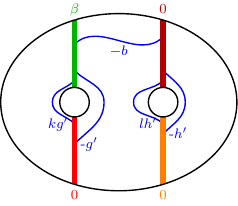}{
	\begin{tikzpicture}[scale=.6,,every node/.style={scale=.5}]
	\pantsparams{$0$}{$0$}{$\beta$}{$0$}{}{}{}{};
	\pantsstpq{$kg^\prime$}{-$g^\prime$}{$lh^\prime$}{-$h^\prime$}{$-b$};
	\end{tikzpicture}
}
\end{array}
\times
\frac{1}{p}\sum_m \omega^{m(\beta-lx-z)}.
\end{align}
The final sum is zero unless the exponent is 0. Therefore, $\beta=lx+z$. The result of the fusion is
\begin{align}
\defect{X_k}{X_k}{a}{x}{}\otimes \defect{X_l}{X_l}{b}{y}{}=\defect{X_{km}}{X_{kl}}{a+kb}{z+lx}{}.
\end{align}
For the special case $k=l=1$, we can identify this with the known anyon fusion rule
\begin{align}
m^ae^x\times m^be^y=m^{a+b}e^{x+y}.
\end{align}
\end{exmp}

\begin{exmp}[$\defect{X_k}{X_l}{}{}{} \otimes \defect{L}{L}{}{}{}$]

Consider the defect fusion
\begin{align}
\defect{X_k}{X_l}{}{}{} \otimes \defect{L}{L}{c}{z}{} \to \defect{L}{L}{\alpha}{\zeta}{}
\end{align}
First, we inflate $\defect{L}{L}{\alpha}{\zeta}{}$ to a 4-string idempotent:
\begin{align}
\frac{1}{p}\sum_h\omega^{h\zeta}
\begin{array}{c}
\includeTikz{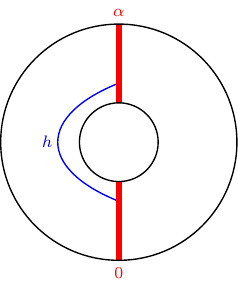}{
	\begin{tikzpicture}[scale=0.8,,every node/.style={scale=.6}]
	\annparamss{$0$}{$\alpha$}{}{};
	\annss{$h$}{};
	\end{tikzpicture}
}
\end{array}
\linf\frac{1}{p}\sum_h\omega^{h\zeta}
\begin{array}{c}
\includeTikz{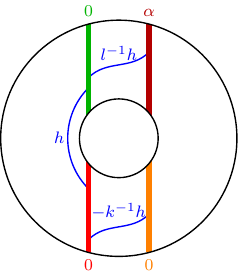}{
	\begin{tikzpicture}[scale=0.8,,every node/.style={scale=0.6}]
		\annparamstpq{$0$}{$0$}{$0$}{$\alpha$}{}{}{}{};
		\annstpq{$h$}{$l^{-1}h$}{}{$-k^{-1}h$};
	\end{tikzpicture}
}
\end{array}.
\end{align}
In order to find a pair-of-pants absorbing the inflated idempotent at the waist, and $\defect{X_k}{X_l}{}{}{}$ and $\defect{L}{L}{c}{z}{}$ on the legs respectively, we must have $\alpha = c$. The general pair-of-pants looks like
\begin{align}
  \frac{1}{p^2} \sum_{g,h} \omega^{gz+h\zeta}
\begin{array}{c}
\includeTikz{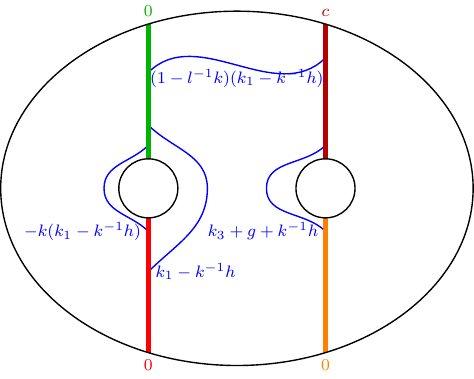}{
	\begin{tikzpicture}[scale=1.2,,every node/.style={scale=.6}]
	\pantsparams{$0$}{$0$}{$0$}{$c$}{}{}{}{};
	\pantsstpq{$-k(k_1-k^{-1}h)$}{$k_1 - k^{-1}h$}{$k_3+g+k^{-1}h$}{}{$(1-l^{-1}k)(k_1-k^{-1}h)$};
	\end{tikzpicture}
}
\end{array} = \frac{\omega^{k k_1 \zeta - k_1 z - k_3 z}}{p^2} \sum_{g,h} \omega^{gz + h k \zeta}
\begin{array}{c}
\includeTikz{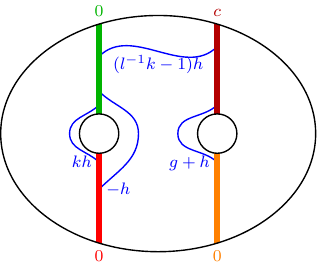}{
	\begin{tikzpicture}[scale=0.8,,every node/.style={scale=.6}]
	\pantsparams{$0$}{$0$}{$0$}{$c$}{}{}{}{};
	\pantsstpq{$kh$}{$-h$}{$g + h$}{}{$(l^{-1}k-1)h$};
	\end{tikzpicture}
}
\end{array}
\end{align}
which is non-zero for all choices of $\zeta$. Therefore
\begin{align}
\defect{X_k}{X_l}{}{}{} \otimes \defect{L}{L}{c}{z}{} = \oplus_{\zeta} \defect{L}{L}{c}{\zeta}{}
\end{align}
\end{exmp}

\begin{exmp}[$\defect{F_0}{X_l}{}{}{} \otimes \defect{L}{L}{}{}{}$]
Consider the defect fusion
\begin{align}
  \defect{F_0}{X_l}{x}{}{} \otimes \defect{L}{L}{c}{z}{} \to \defect{L}{L}{\alpha}{\zeta}{}
\end{align}
Inflating $\defect{L}{L}{\alpha}{\zeta}{}$ to a 4-string idempotent gives
\begin{align}
\frac{1}{p}\sum_{\gamma}\omega^{\gamma \zeta}
\begin{array}{c}
\includeTikz{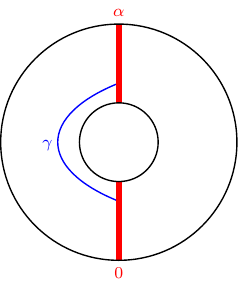}{
	\begin{tikzpicture}[scale=0.8,,every node/.style={scale=.6}]
	\annparamss{$0$}{$\alpha$}{}{};
	\annss{$\gamma$}{};
	\end{tikzpicture}
}
\end{array}
\linf\frac{1}{p^2}\sum_{\gamma,\delta}\omega^{\delta \mu + \gamma \zeta}
\begin{array}{c}
\includeTikz{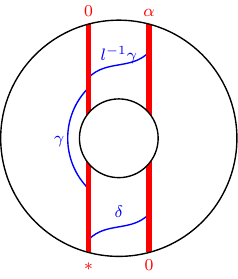}{
	\begin{tikzpicture}[scale=0.8,,every node/.style={scale=0.6}]
		\annparamstpq{$*$}{$0$}{$0$}{$\alpha$}{}{}{}{};
		\annssss{$\gamma$}{$l^{-1}\gamma$}{}{$\delta$};
	\end{tikzpicture}
}
\end{array}.
\end{align}
The general pair-of-pants absorbing $\defect{F_0}{X_l}{x}{}{}$ and $\defect{L}{L}{c}{z}{}$ on the legs and the inflated idempotent on the belt is
\begin{align}
  \frac{1}{p^4} \sum_{g,h,\gamma,\delta} \omega^{gx + hz + \delta \mu + \gamma \zeta}
\begin{array}{c}
\includeTikz{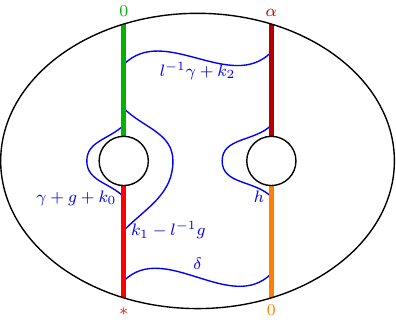}{
	\begin{tikzpicture}[scale=1,,every node/.style={scale=.6}]
	  \pantsparams{$*$}{$0$}{$0$}{$\alpha$}{}{}{}{};
  \generalpantsstpq{$\gamma + g + k_0$}{$k_1-l^{-1}g$}{$h$}{}{}{$l^{-1}\gamma + k_2$}{}{$\delta$};
	\end{tikzpicture}
}
\end{array} = \frac{1}{p^4} \sum_{g,h} \omega^{gx + hz + \delta \mu + \gamma \zeta}
\begin{array}{c}
\includeTikz{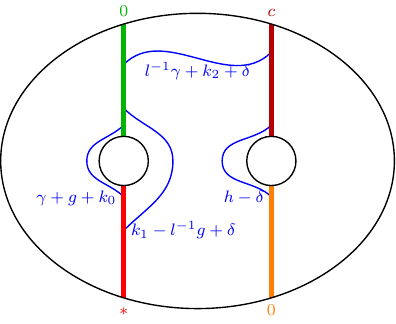}{
	\begin{tikzpicture}[scale=1,,every node/.style={scale=.6}]
	\pantsparams{$*$}{$0$}{$0$}{$c$}{}{}{}{};
	\pantsstpq{$\gamma + g + k_0$}{$k_1 - \l^{-1}g + \delta$}{$h - \delta$}{}{$l^{-1}\gamma + k_2 + \delta$};
	\end{tikzpicture}
}
\end{array}
\end{align}
In order for the objects to match on the green string, we must have $k_0 = l(k_2 - k_1)$. This equation together with the transformations $\gamma \to \gamma - l k_2$ and $g \to g + l k_1$ let us phase away $k_0,k_1$ and $k_2$ giving the following pants:
\begin{align}
\frac{1}{p^4} \sum_{g,h,\gamma,\delta} \omega^{gx + hz + \delta \mu + \gamma \zeta}
\begin{array}{c}
\includeTikz{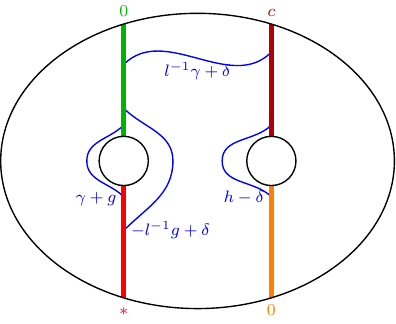}{
	\begin{tikzpicture}[scale=1,,every node/.style={scale=.6}]
	\pantsparams{$*$}{$0$}{$0$}{$c$}{}{}{}{};
	\pantsstpq{$\gamma + g$}{$- \l^{-1}g + \delta$}{$h - \delta$}{}{$l^{-1}\gamma + \delta$};
	\end{tikzpicture}
}
\end{array}
\end{align}
If we define
\begin{align}
  b &= l^{-1} \gamma + \delta \\
  c &= -l^{-1} g + \delta \\
  d &= h - \delta
\end{align}
then $(g,h,\gamma,\delta) \to (b,c,d,\gamma)$ is invertible over $\ZZ{p}$ and we get the following pants:
\begin{align}
\frac{1}{p^3} \sum_{b,c,d} \omega^{(lx + z + \mu)b - lcx + dz}
\begin{array}{c}
\includeTikz{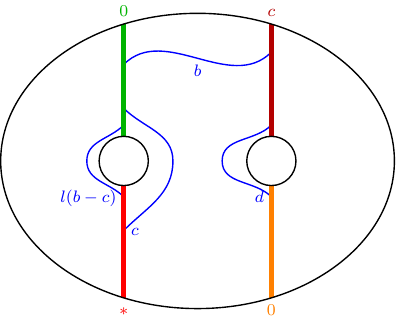}{
	\begin{tikzpicture}[scale=1,,every node/.style={scale=.6}]
	\pantsparams{$*$}{$0$}{$0$}{$c$}{}{}{}{};
	\pantsstpq{$l(b-c)$}{$c$}{$d$}{}{$b$};
	\end{tikzpicture}
}
\end{array}\times\frac{1}{p}\sum_{\gamma}\omega^{(\zeta - x - l^{-1}z - l^{-1}\mu)\gamma}
\end{align}
which is nonzero only if $\zeta = x + l^{-1}(\mu + z)$. Therefore we have
\begin{align}
\defect{X_k}{X_l}{}{}{} \otimes \defect{L}{L}{c}{z}{} \to \defect{L}{L}{c}{x + l^{-1}(\mu + z)}{} 
\end{align}
\end{exmp}

\begin{exmp}[$\defect{X_k}{X_l}{}{}{} \otimes \defect{X_m}{X_n}{}{}{}$]

Consider the defect fusion
\begin{align}
	\defect{X_k}{X_l}{}{}{} \otimes \defect{X_m}{X_n}{}{}{}\to 
	\begin{cases} 
		\defect{X_{km}}{X_{ln}}{}{}{} & km \neq ln\\
		\defect{X_{km}}{X_{km}}{\alpha}{\beta}{} & km = ln
	\end{cases}.
\end{align}
where we assume $k\neq l$ and $m\neq n$. We will present the calculation for these two cases separately.
\subsubsection*{Case I: $km=ln$}
Inflating $\defect{X_{km}}{X_{km}}{\alpha}{\beta}{}$ to a 4-string idempotent gives
\begin{align}
	\frac{1}{p}\sum_{\gamma}\omega^{\gamma\beta}
	\begin{array}{c}
		\includeTikz{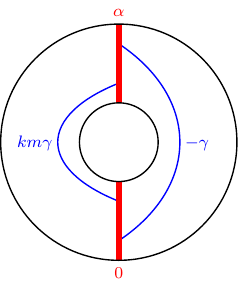}{
			\begin{tikzpicture}[scale=0.8,,every node/.style={scale=.6}]
				\annparamss{$0$}{$\alpha$}{}{};
				\annss{$km\gamma$}{$-\gamma$};
			\end{tikzpicture}
		}
	\end{array}
	\linf\frac{1}{p}\sum_{\gamma}\omega^{\gamma\beta}
	\begin{array}{c}
		\includeTikz{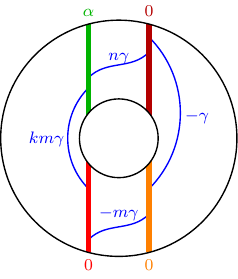}{
			\begin{tikzpicture}[scale=0.8,,every node/.style={scale=0.6}]
			\annparamstpq{$0$}{$0$}{$\alpha$}{$0$}{}{}{}{};
			\annstpq{$km\gamma$}{$n\gamma$}{$-\gamma$}{$-m\gamma$};
			\end{tikzpicture}
		}
	\end{array}.
\end{align}
The general pair-of-pants absorbing $\defect{X_k}{X_l}{}{}{}$ and $\defect{X_m}{X_n}{}{}{}$ on the legs and the inflated idempotent on the belt is
\begin{align}
	\frac{1}{p}\sum_{\gamma}\omega^{\gamma\beta}
	\begin{array}{c}
		\includeTikz{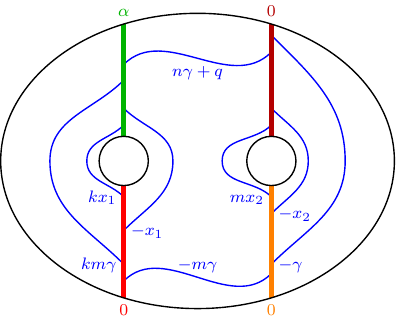}{
			\begin{tikzpicture}[scale=1,,every node/.style={scale=.6}]
			\pantsparams{$0$}{$0$}{$\alpha$}{$0$}{}{}{}{};
			\generalpantsstpq{$kx_1$}{$-x_1$}{$mx_2$}{$-x_2$}{$km\gamma$}{$n\gamma+q$}{$-\gamma$}{$-m\gamma$};
			\end{tikzpicture}
		}
	\end{array} = \frac{1}{p}\sum_{\gamma}\omega^{\gamma\beta}
	\begin{array}{c}
		\includeTikz{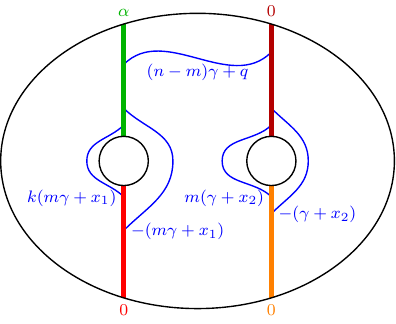}{
			\begin{tikzpicture}[scale=1,,every node/.style={scale=.6}]
			\pantsparams{$0$}{$0$}{$\alpha$}{$0$}{}{}{}{};
			\pantsstpq{$k(m\gamma+x_1)$}{$-(m\gamma+x_1)$}{$m(\gamma+x_2)$}{$-(\gamma+x_2)$}{$(n-m)\gamma+q$};
			\end{tikzpicture}
		}
	\end{array},
\end{align}
where $q=(l^{-1}k-1)x_1-l^{-1}\alpha$, $x_2=-q(m-n)^{-1}$.
Making the change of variables $\gamma^\prime=m\gamma+x_1$ eliminate $x_1$ up to a global phase. Therefore, there is a single, nonzero map
\begin{align}
	\defect{X_k}{X_l}{}{}{} \otimes \defect{X_m}{X_n}{}{}{} \to \defect{X_{km}}{X_{km}}{\alpha}{\beta}{} 
\end{align}
for all $\alpha,\,\beta$.

\subsubsection*{Case II: $km\neq ln$}

Inflating $\defect{X_{km}}{X_{ln}}{}{}{}$ to a 4-string idempotent gives
\begin{align}
	\begin{array}{c}
		\includeTikz{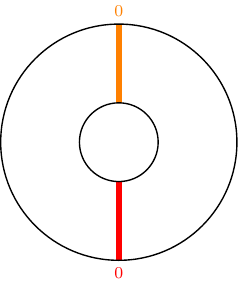}{
			\begin{tikzpicture}[scale=0.8,,every node/.style={scale=.6}]
			\annparamss{$0$}{$0$}{}{};
			\annst{}{};
			\end{tikzpicture}
		}
	\end{array}
	\linf
	\begin{array}{c}
		\includeTikz{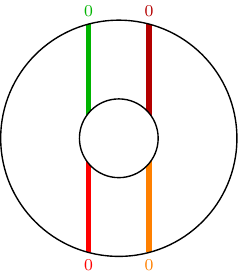}{
			\begin{tikzpicture}[scale=0.8,,every node/.style={scale=0.6}]
			\annparamstpq{$0$}{$0$}{$0$}{$0$}{}{}{}{};
			\annstpq{}{}{}{};
			\end{tikzpicture}
		}
	\end{array}.
\end{align}

The general pair-of-pants absorbing $\defect{X_k}{X_l}{}{}{}$ and $\defect{X_m}{X_n}{}{}{}$ on the legs and the inflated idempotent on the belt is
\begin{align}
	\begin{array}{c}
		\includeTikz{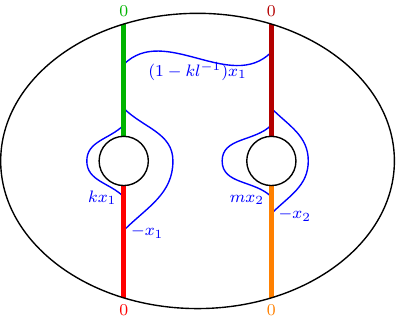}{
			\begin{tikzpicture}[scale=1,,every node/.style={scale=.6}]
				\pantsparams{$0$}{$0$}{$0$}{$0$}{}{}{}{};
				\pantsstpq{$kx_1$}{$-x_1$}{$mx_2$}{$-x_2$}{$(1-kl^{-1})x_1$};
			\end{tikzpicture}
		}
	\end{array},
\end{align}
where $x_2=x_1(kl^{-1}-1)(m-n)^{-1}$.
where $q=(l^{-1}k-1)x_1-l^{-1}\alpha$, $x_2=-q(m-n)^{-1}$. There are no further constraints on the maps. Therefore, there are $p$ distinct maps between these defects
\begin{align}
	\defect{X_k}{X_l}{}{}{} \otimes \defect{X_m}{X_n}{}{}{} = p\cdot \defect{X_{km}}{X_{ln}}{}{}{}. 
\end{align}

\end{exmp}

\begin{exmp}[$\defect{F_q}{X_l}{}{}{} \otimes \defect{F_s}{X_n}{}{}{}$, which includes fusing `twists']

A particular set of defects of the $\ZZ{2}$ Kitaev model identified in \onlinecite{Bombin2010}, where they were referred to as \emph{twists}. In our notation, twists occur at the interface of an $X_1$ and $F_1$ domain wall. In this example, we compute the fusion of generalizations of twists in the $\ZZ{p}$ Kitaev model.

Consider the defect fusion
\begin{align}
	\defect{F_q}{X_l}{x}{}{} \otimes \defect{F_s}{X_n}{z}{}{}\to 
	\begin{cases} 
		\defect{X_{q^{-1}s}}{X_{ln}}{}{}{} & q^{-1}s \neq ln\\
		\defect{X_{ln}}{X_{ln}}{\alpha}{\beta}{} & q^{-1}s = ln
	\end{cases}.
\end{align}

\subsubsection*{Case I: $q^{-1}s\neq ln$}

When $q^{-1}s \not= ln$, the fusion looks like
\begin{align}
\defect{F_q}{X_l}{x}{}{}\otimes \defect{F_s}{X_n}{z}{}{} \to \defect{X_{q^{-1}s}}{X_{ln}}{}{}{},
\end{align}
but we need to check for multiplicity. We can inflate $\defect{X_{q^{-1}s}}{X_{ln}}{}{}{}$ to a 4-string annulus
\begin{align}
\begin{array}{c}
\includeTikz{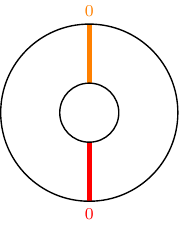}{
	\begin{tikzpicture}[scale=.6,,every node/.style={scale=.6}]
	\annparamst{$0$}{$0$}{}{};
	\annst{}{};
	\end{tikzpicture}
}
\end{array}
\linf\frac{1}{p}\sum_{\gamma}
\begin{array}{c}
\includeTikz{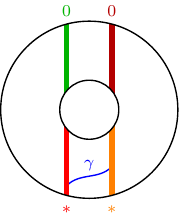}{
	\begin{tikzpicture}[scale=.6,,every node/.style={scale=.6}]
	\annparamstpq{$*$}{$*$}{$0$}{$0$}{}{}{}{};
	\annstpq{}{}{}{$\gamma$};
	\end{tikzpicture}
}
\end{array}.
\end{align}
This implies that the general pair-of-pants absorbing the 2-string defects on the legs and the 4-string defect on the waist is
\begin{align}
\frac{1}{p^3}\sum_{g,h,\gamma} \Theta_{x,-ql}(g) \Theta_{z,-ns}(h) \omega^{-qk_0g - sk_3h}
&\begin{array}{c}
\includeTikz{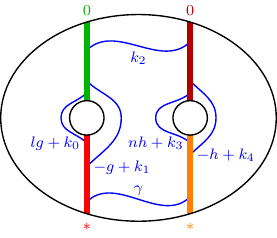}{
	\begin{tikzpicture}[scale=.7,,every node/.style={scale=.55}]
	\pantsparams{$*$}{$*$}{$0$}{$0$}{}{}{}{};
	\generalpantsstpq{$lg+k_0$}{$-g+k_1$}{$nh+k_3$}{$-h+k_4$}{}{$k_2$}{}{$\gamma$};
	\end{tikzpicture}
}
\end{array}
\\=
\frac{1}{p^3}\sum_{g,h,\gamma} \Theta_{x,-ql}(g) \Theta_{z,-ns}(h) \omega^{-qk_0g - sk_3h+s \gamma (h-k_4)}
&\begin{array}{c}
\includeTikz{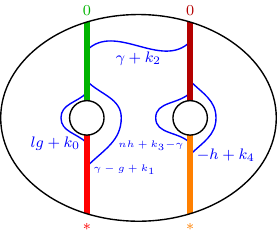}{
	\begin{tikzpicture}[scale=.7,,every node/.style={scale=.55}]
	\pantsparams{$*$}{$*$}{$0$}{$0$}{}{}{}{};
	\pantsstpq{$lg+k_0$}{\tiny$\gamma-g+k_1$}{\tiny$nh+k_3{-}\gamma$}{$-h+k_4$}{$\gamma+k_2$};
	\end{tikzpicture}
}
\end{array}.
\end{align}
For the labels to match, we must have
\begin{align}
  k_0 + l k_1 - l k_2 &= 0 \\
  k_3 + n k_4 + k_2 &= 0.
  \end{align}
If we make the transformation
\begin{align}
  g &= a - l^{-1}k_0 \\
  \gamma &= b - k_2 \\
  h &= c + k_4
\end{align}
then we get
\begin{align}
  \frac{1}{p^3}\sum_{a,b,c} \Theta_{x,-ql}(a)\Theta_{z,-ns}(c) \omega^{sbc}
\begin{array}{c}
\includeTikz{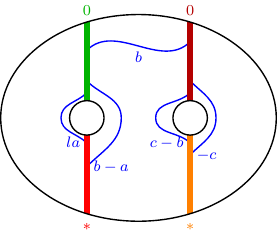}{
	\begin{tikzpicture}[scale=.7,,every node/.style={scale=.55}]
	\pantsparams{$*$}{$*$}{$0$}{$0$}{}{}{}{};
	\pantsstpq{$la$}{$b-a$}{$c-b$}{$-c$}{$b$};
	\end{tikzpicture}
}
\end{array}.
\end{align}
So there is 1 map 
\begin{align}
\defect{F_q}{X_l}{x}{}{}\otimes \defect{F_s}{X_n}{z}{}{}\to \defect{X_{q^{-1}s}}{X_{ln}}{}{}{}.
\end{align}

\subsubsection*{Case II: $q^{-1}s=ln$}
When $q^{-1}s=ln$, the fusion looks like
\begin{align}
  \defect{F_q}{X_l}{x}{}{}\otimes \defect{F_s}{X_n}{z}{}{}\to \defect{X_{ln}}{X_{ln}}{\alpha}{\zeta}{}.
\end{align}
We can inflate $\defect{X_{ln}}{X_{ln}}{\alpha}{\zeta}{}$ to a 4-string annulus
\begin{align}
\frac{1}{p}\sum_{\gamma}\omega^{\gamma \zeta}
\begin{array}{c}
\includeTikz{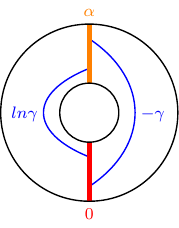}{
	\begin{tikzpicture}[scale=.6,,every node/.style={scale=.6}]
	\annparamst{$0$}{$\alpha$}{}{};
	\annst{$ln\gamma$}{$-\gamma$};
	\end{tikzpicture}
}
\end{array}
\linf\frac{1}{p^2}\sum_{\gamma,\delta}\omega^{\gamma \zeta}
\begin{array}{c}
\includeTikz{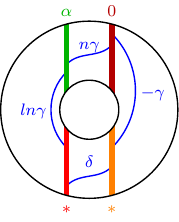}{
	\begin{tikzpicture}[scale=.6,,every node/.style={scale=.6}]
	\annparamstpq{$*$}{$*$}{$\alpha$}{$0$}{}{}{}{};
	\annstpq{$ln\gamma$}{$n\gamma$}{$-\gamma$}{$\delta$};
	\end{tikzpicture}
}
\end{array}.
\end{align}

Next, we can find the pant that absorbs the 2-string idempotents on the legs and 4-string idempotent on the waist. The most general pant that will absorb the appropriate idempotents is 
\begin{align}
\frac{1}{p^4}\sum_{g,h,\gamma,\delta}\Theta_{x,-ql}(g) \Theta_{z,-ns}(h) \omega^{\gamma \zeta-hk_3s - gk_0q}
\begin{array}{c}
\includeTikz{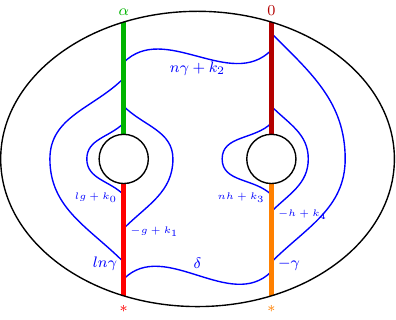}{
	\begin{tikzpicture}[scale=1,,every node/.style={scale=.55}]
	\pantsparams{$*$}{$*$}{$\alpha$}{$0$}{}{}{}{};
	\generalpantsstpq{\tiny$lg+k_0$}{\tiny$-g+k_1$}{\tiny$nh+k_3$}{\tiny$-h+k_4$}{$ln\gamma$}{$n\gamma + k_2$}{$-\gamma$}{$\delta$};
	\end{tikzpicture}
}
\end{array}\\
=
\frac{1}{p^4}\sum_{g,h,\gamma,\delta}
\Theta_{x,-ql}(g) \Theta_{z,-ns}(h) \omega^{\gamma \zeta-hk_3s - gk_0q + ln \gamma(-g+k_1)q +\delta(\gamma + h - k_4)s}
\begin{array}{c}
\includeTikz{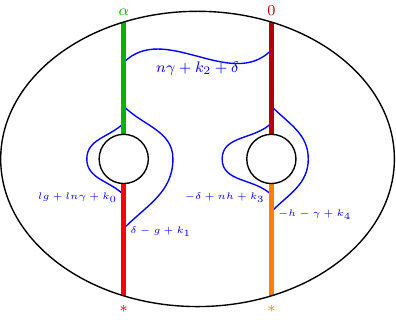}{
	\begin{tikzpicture}[scale=1,,every node/.style={scale=.55}]
	\pantsparams{$*$}{$*$}{$\alpha$}{$0$}{}{}{}{};
	\pantsstpq{\tiny$lg + ln\gamma + k_0$}{\tiny$\delta-g+k_1$}{\tiny$-\delta+nh+k_3$}{\tiny$-h-\gamma+k_4$}{$n\gamma+k_2+\delta$};
	\end{tikzpicture}
}
\end{array}.
\end{align}
For the labels to match up
\begin{align}
  k_0 + l k_1 - l k_2-\alpha &= 0 \\
  k_3 + n k_4 + k_2 &= 0.
  \end{align}
If we make the change of coordinates
\begin{align}
  a &= \delta - g \\
  b &= n \gamma + \delta \\
  c &= -h - \gamma
\end{align}
then we get
\begin{align}
  \frac{1}{p^4}\sum_{a,b,c,\gamma}\omega^{{\rm exponent}}
\begin{array}{c}
\includeTikz{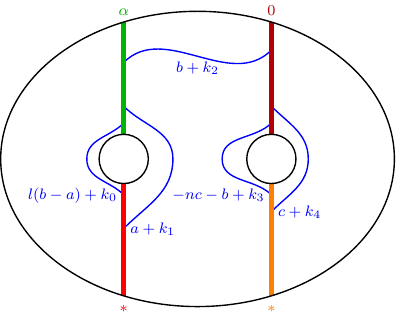}{
	\begin{tikzpicture}[scale=1,,every node/.style={scale=.55}]
	\pantsparams{$*$}{$*$}{$\alpha$}{$0$}{}{}{}{};
	\pantsstpq{$l(b-a)+k_0$}{$a+k_1$}{$-nc-b+k_3$}{$c+k_4$}{$b+k_2$};
	\end{tikzpicture}
}
\end{array}
\end{align}
where
\begin{align*}
  {\rm exponent} = &-2^{-1} a^2 l q+k_0 q (a-b+\gamma  n)+a b l q-a x-2^{-1} b^2 l q+ \\
   &lnq \left(k_4 (\gamma  n-b)+k_3 (c+\gamma )+\gamma  k_1\right)-b c l n q+bx-2^{-1} c^2 l n^2 q+\gamma  \zeta -c z-\gamma  n x-\gamma  z
  \end{align*}
Making the transformation
\begin{align}
  a &\to a - k_1 \\
  b &\to b - k_2 \\
  c &\to c - k_4
\end{align}
gives
\begin{align}
  \frac{1}{p^4}\sum_{a,b,c,\gamma}\omega^{{\rm exponent}}
\begin{array}{c}
\includeTikz{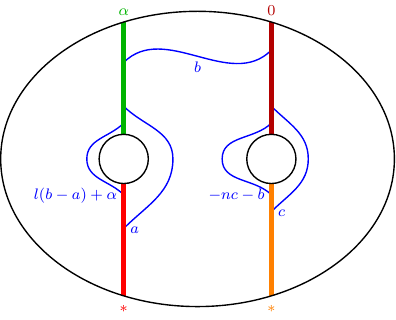}{
	\begin{tikzpicture}[scale=1,,every node/.style={scale=.55}]
	\pantsparams{$*$}{$*$}{$\alpha$}{$0$}{}{}{}{};
	\pantsstpq{$l(b-a)+\alpha$}{$a$}{$-nc-b$}{$c$}{$b$};
	\end{tikzpicture}
}
\end{array}
\end{align}
where
\begin{align*}
  {\rm exponent} =& -2^{-1} a^2 l q+a b l q+\alpha  a q-a x-2^{-1} b^2 l q-b c l n q-\\ 
  &\alpha b q+b x-2^{-1} c^2 l n^2 q-c z+ \gamma  (\zeta +\alpha  n q-n x-z)
\end{align*}
This implies that $\zeta = nx + z - \alpha n q$, so
\begin{align}
\defect{F_q}{X_l}{x}{}{} \otimes \defect{F_s}{X_n}{z}{}{} \to \oplus_{\alpha} \defect{X_{ln}}{X_{ln}}{\alpha}{z+nx-nq\alpha}{}
  \end{align}
for odd primes. The same calculation for $p=2$ yields
\begin{align}
	\defect{F_1}{X_1}{x}{}{} \otimes \defect{F_1}{X_1}{z}{}{} \to \oplus_{\alpha} \defect{X_{1}}{X_{1}}{\alpha}{z+x+\alpha}{}.
\end{align}
This result was presented in \onlinecite{Bombin2010} as $\sigma^x\sigma^z=\sum_j e^{j+x+z} m^{j}$.

\end{exmp}

\begin{exmp}[$\defect{F_q}{X_l}{}{}{} \otimes \defect{X_m}{X_m}{}{}{}$]
\begin{align}
\defect{F_q}{X_l}{x}{}{}\otimes\defect{X_m}{X_m}{c}{z}{}\to \defect{F_{qm}}{X_{lm}}{\alpha}{}{}
\end{align}

We can inflate $\defect{F_{qm}}{X_{lm}}{\alpha}{}{}$ to a 4-string annulus
\begin{align}
\frac{1}{p}\sum_g\Theta_{\alpha,-qlm^2}(g)
\begin{array}{c}
\includeTikz{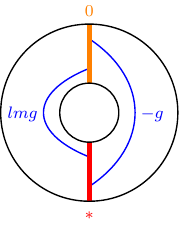}{
	\begin{tikzpicture}[scale=.6,,every node/.style={scale=.6}]
	\annparamst{$*$}{$0$}{}{};
	\annst{$lmg$}{$-g$};
	\end{tikzpicture}
}
\end{array}
\linf
\frac{1}{p}\sum_g\Theta_{\alpha,-qlm^2}(g)
\begin{array}{c}
\includeTikz{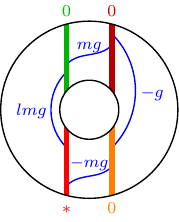}{
	\begin{tikzpicture}[scale=.6,,every node/.style={scale=.6}]
	\annparamstpq{$*$}{$0$}{$0$}{$0$}{}{}{}{};
	\annstpq{$lmg$}{$mg$}{$-g$}{$-mg$};
	\end{tikzpicture}
}
\end{array}.
\end{align}

Next, we can find a pant mapping the 2-string defects to the 4-string defect 
\begin{align}
\frac{1}{p^3}&\sum_{g,h,k}\Theta_{x,-ql}(g)\Theta_{\alpha,-qlm^2}(k)\omega^{hz-ql\mu g}
\begin{array}{c}
\includeTikz{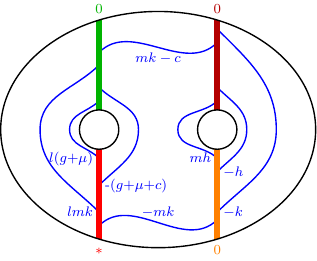}{
	\begin{tikzpicture}[scale=.8,every node/.style={scale=.5}]
	\pantsparams{$*$}{$0$}{$0$}{$0$}{}{}{}{};
	\generalpantsstpq{$l(g{+}\mu)$}{-$(g{+}\mu{+}c)$}{$mh$}{$-h$}{$lmk$}{$mk-c$}{$-k$}{$-mk$};
	\end{tikzpicture}
}
\end{array}\\
&=
\frac{1}{p^3}\sum_{g,h,k}\Theta_{x,-ql}(g)\Theta_{\alpha,-qlm^2}(k)\omega^{hz-ql\mu g-qlmk(g+\mu+c)}
\begin{array}{c}
\includeTikz{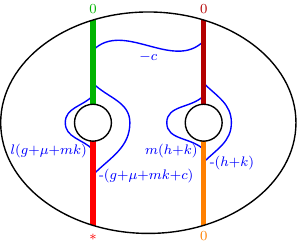}{
	\begin{tikzpicture}[scale=.75,,every node/.style={scale=.5}]
	\pantsparams{$*$}{$0$}{$0$}{$0$}{}{}{}{};
	\pantsstpq{$l(g{+}\mu{+}mk)$}{-$(g$+$\mu{+}mk$+$c)$}{$m(h{+}k)$}{-$(h{+}k)$}{$-c$};
	\end{tikzpicture}
}
\end{array}\\
&\propto  
\frac{1}{p^2}\sum_{g,h}\Theta_{x,-ql}(g)\omega^{hz}
\begin{array}{c}
\includeTikz{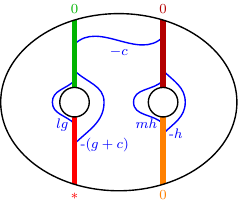}{
	\begin{tikzpicture}[scale=.6,,every node/.style={scale=.5}]
	\pantsparams{$*$}{$0$}{$0$}{$0$}{}{}{}{};
	\pantsstpq{$lg$}{-$(g+c)$}{$mh$}{-$h$}{$-c$};
	\end{tikzpicture}
}
\end{array}
\times \frac{1}{p}\sum_{k}\omega^{k(\alpha-z-m(x+qlc))}.
\end{align}
Therefore
\begin{align}
\defect{F_q}{X_l}{x}{}{} \otimes \defect{X_m}{X_m}{c}{z}{}=\defect{F_{qm}}{X_{lm}}{z+m(x+qlc)}{}{}.
\end{align}
For the special case $q=l=m=1$, this recovers the known fusion rule of `twists' with anyons\cite{Bombin2010,Barkeshli2014,Bridgeman2017}
\begin{align}
\sigma^{x}\times m^c e^z=\sigma^{x+z+c}.
\end{align}

\end{exmp}

\section{The Vertical Fusion Algorithm} \label{sec:vertical_fusion_algorithm}

In this section, we explain the vertical defect fusion algorithm. The structure of this section parallels Section~\ref{sec:hotizontal_defect_fusion_algorithm}. Since there is no domain wall fusion, vertical fusion is simpler than horizontal fusion. No inflation is required, so the algorithm proceeds in three key steps:
\begin{enumerate} 
	\item Determine idempotents corresponding to source defects using Table~\ref{tab:idempotents}.
	\item Determine idempotent for target defect using Tables~\ref{tab:zptable} and \ref{tab:idempotents}.
	\item Find a nonzero pants diagram that absorbs the source idempotents on the legs and (inflated) target on the waist.
\end{enumerate}

We shall explain how the algorithm works, step by step, using an example.

\begin{exmp}[$\defect{R}{R}{a}{x}{} \circ \defect{R}{F_s}{z}{}{}$]
Consider the vertical defect fusion
\begin{align}
  \defect{R}{R}{a}{x}{} \circ \defect{R}{F_s}{z}{}{}.
\end{align}
The convention for the $\circ$ operator is that the left argument is the bottom defect and the right argument is the top defect.
\subsubsection{Writing down the defect idempotents to compose}
The first step in the procedure is to look up the idempotents that we are composing. These are found in Table~\ref{tab:idempotents}. In our current example, we have
\begin{align}
\defect{R}{F_r}{z}{}{}=\frac{1}{p}\sum_h\omega^{hz}
				\begin{array}{c}
				\includeTikz{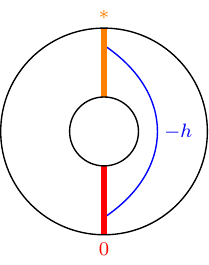}
				{
					\begin{tikzpicture}[scale=.7,every node/.style={scale=.7}]
					\annparamst{$0$}{$*$}{}{};
					\annst{}{$-h$};
					\end{tikzpicture}
				}
				\end{array}, \quad
\defect{R}{R}{a}{x}{}=\frac{1}{p}\sum_g\omega^{gx}
				\begin{array}{c}
				\includeTikz{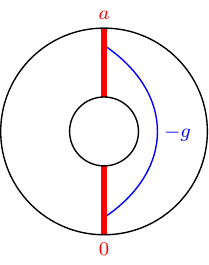}
				{
					\begin{tikzpicture}[scale=.7,every node/.style={scale=.7}]
					\annparamss{$0$}{$a$}{}{};
					\annss{}{$-g$};
					\end{tikzpicture}
				}
				\end{array}.
\end{align}
\subsubsection{Writing down the target idempotent}
Next, we need to write down the target defect idempotent. These can be found in Table~\ref{tab:idempotents}.
\begin{align}
\defect{R}{F_r}{\zeta}{}{}=\frac{1}{p}\sum_{\gamma}\omega^{\gamma\zeta}
				\begin{array}{c}
				\includeTikz{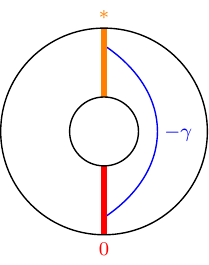}
				{
					\begin{tikzpicture}[scale=.7,every node/.style={scale=.7}]
					\annparamst{$0$}{$*$}{}{};
					\annst{}{$-\gamma$};
					\end{tikzpicture}
				}
				\end{array}.
        \end{align}
\subsubsection{Decorating the pants}
At this point we have extracted all the data we need from the tables. Now we need to find all the pants diagrams which absorb our source defect idempotents on the legs and our target defect idempotent on the waist. In this case, the, the pants are oriented perpendicular to the horizontal case. The most general pair of pants for any 3 domain walls $\mathcal{M},\mathcal{N},\mathcal{P}$ is
\begin{align}
\begin{array}{c}
		\includeTikz{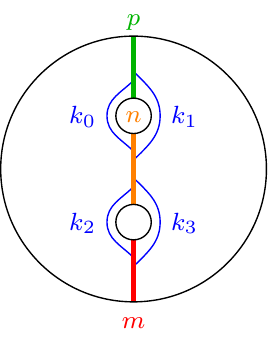}{
			\begin{tikzpicture}[scale=.9,,every node/.style={scale=.9}]
			\vpantsparams{$m$}{$n$}{$p$};
			\vpantsstp{$k_2$}{$k_3$}{$k_0$}{$k_1$};
			\end{tikzpicture}}
	\end{array}.
\end{align}
We use the following equation to bring every vertical pair of pants into a standard form
\begin{align}
	\begin{array}{c}
		\includeTikz{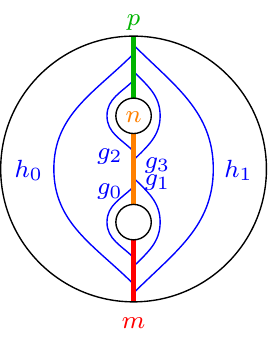}{
			\begin{tikzpicture}[scale=.9,,every node/.style={scale=.9}]
			\vpantsparams{$m$}{$n$}{$p$};
			\vgeneralpantsstp{$g_0$}{$g_1$}{$g_2$}{$g_3$}{$h_0$}{$h_1$};
			\end{tikzpicture}}
	\end{array}
	&=\frac{\Omega_M(h_0^{-1},g_1^{-1})\Omega_N(h_0^{-1},g_1)}{\Omega_P(h_0,g_3)\Omega_N(h_0,g_3^{-1})}
	\begin{array}{c}
		\includeTikz{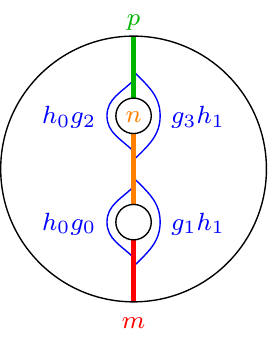}{
			\begin{tikzpicture}[scale=.9,,every node/.style={scale=.9}]
			\vpantsparams{$m$}{$n$}{$p$};
			\vpantsstp{$h_0g_0$}{$g_1h_1$}{$h_0g_2$}{$g_3h_1$};
			\end{tikzpicture}}
	\end{array}.
\end{align}
When we insert our source defect idempotents on the legs and our target defect idempotent on the waist, we get
\begin{align}\frac{1}{p^3}\sum_{g,h,\gamma}\omega^{hz+gx+\gamma\zeta+sk_0h}
\begin{array}{c}
		\includeTikz{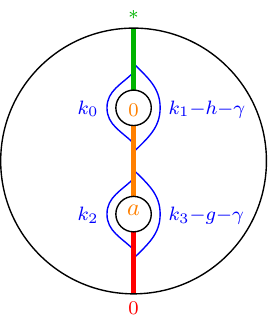}{
			\begin{tikzpicture}[scale=.9,,every node/.style={scale=.7}]
			\vpantsparams{$0$}{$0$}{$*$};
			\vpantsstp{$k_2$}{$k_3{-}g{-}\gamma$}{$k_0$}{$k_1{-}h{-}\gamma$};
			\node[orange,inline text] at (0,-.55) {$a$};
			\end{tikzpicture}}
	\end{array}.
\end{align}
The variables $k_1$ and $k_2$ can be converted into a global phase by the translations $h \to h + k_1, g \to g + k_3$. For the object labels to match up, we need $k_0 = a$ and $k_2 = 0$. Now, if we make the change of variables,
\begin{align}
  c &= g + \gamma \\
  b &= h + \gamma
\end{align}
we get
\begin{align}\frac{1}{p^3}\sum_{c,b,\gamma}\omega^{(b-\gamma)z+(c-\gamma)x+\gamma\zeta+sa(b-\gamma)}
\begin{array}{c}
		\includeTikz{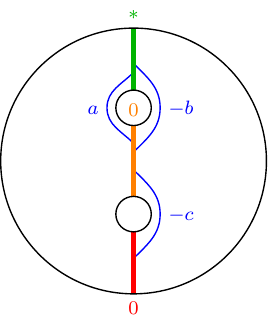}{
			\begin{tikzpicture}[scale=.9,,every node/.style={scale=.7}]
			\vpantsparams{$0$}{$0$}{$*$};
			\vpantsstp{}{$-c$}{$a$}{$-b$};
			\end{tikzpicture}}
	\end{array}.
\end{align}
For this to be non zero, we need to have $\zeta=x+z+sa$,
so
\begin{align}
 \defect{R}{R}{a}{x}{}  \circ  \defect{R}{F_s}{z}{}{}= \defect{R}{F_s}{x+z+sa}{}{}
\end{align}

\end{exmp}

\section{Vertical Defect Fusion} \label{sec:vertical_defect_fusion}

In this section, we present several of the more complicated vertical defect fusion computations.

\begin{exmp}[$\defect{F_q}{X_l}{}{}{} \circ \defect{X_l}{F_s}{}{}{}$]

Consider the defect fusion
\begin{align}
	\defect{F_q}{X_l}{}{}{} \circ \defect{X_l}{F_s}{}{}{}\to 
	\begin{cases} 
		\defect{F_s}{F_s}{\alpha}{\beta}{} & q = s\\
		\defect{F_q}{F_s}{}{}{} & q\neq s
	\end{cases}.
\end{align}

\subsubsection*{Case I: $q=s$}

When $q=s$, the vertical fusion looks like
\begin{align}
 \defect{F_q}{X_l}{x}{}{}  \circ \defect{X_l}{F_s}{z}{}{}  \to \defect{F_s}{F_s}{\zeta}{\eta}{}
  \end{align}
The general pants absorbing $\defect{X_l}{F_s}{z}{}{}$ and $\defect{F_q}{X_l}{x}{}{}$ on the legs and $\defect{F_s}{F_s}{\zeta}{\eta}{}$ on the waist looks like
\begin{align}
  \frac{1}{p^4}\sum_{g,h,\gamma,\delta} \omega^{\gamma\zeta+\delta\eta+s\gamma(h-g)} \Theta_{z,ls}(h)\Theta_{x,-ql}(g)
	\begin{array}{c}
		\includeTikz{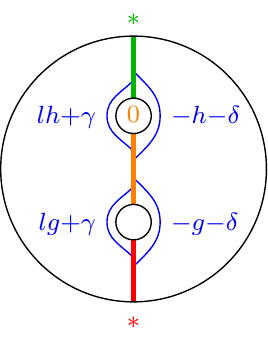}{
			\begin{tikzpicture}[scale=.9,,every node/.style={scale=.9}]
			  \vpantsparams{$*$}{$0$}{$*$};
        \vpantsstp{$lg{+}\gamma$}{$-g{-}\delta$}{$lh{+}\gamma$}{$-h{-}\delta$};
			\end{tikzpicture}}
	\end{array}
\end{align}
If we set
\begin{align}
  a &= -g -\delta \\
  b &= -h - \delta \\
  c &= lh + \gamma
  \end{align}
then this pair of pants becomes
\begin{align}
  \frac{1}{p^4} \sum_{a,b,c,\delta}\omega^{\rm exponent}
	\begin{array}{c}
		\includeTikz{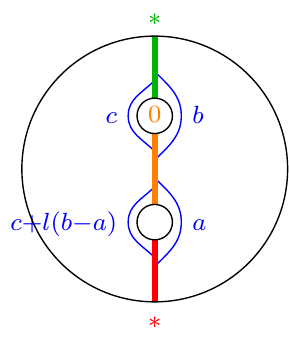}{
			\begin{tikzpicture}[scale=.9,,every node/.style={scale=.9}]
			  \vpantsparams{$*$}{$0$}{$*$};
        \vpantsstp{$c{+}l(b{-}a)$}{$a$}{$c$}{$b$};
			\end{tikzpicture}}
	\end{array}
\end{align}
where
\begin{align}
{\rm exponent}=-2^{-1} a^2 l s+a (b l s+c s-x)-2^{-1} b^2 l s-b (c s-\zeta  l+z)+c \zeta +\delta  (\eta +\zeta  l-x-z).
\end{align}
This implies $\eta = x + z - l\zeta$, so
\begin{align}
 \defect{F_q}{X_l}{x}{}{} \circ   \defect{X_l}{F_s}{z}{}{} \to \oplus_{\zeta} \defect{F_s}{F_s}{\zeta}{x+z-l\zeta}{}
\end{align}
\subsubsection*{Case II: $q\neq s$}
When $q=s$, the vertical fusion looks like
\begin{align}
  \defect{F_q}{X_l}{x}{}{} \circ    \defect{X_l}{F_s}{z}{}{} \to \defect{F_q}{F_s}{}{}{}
  \end{align}
All the variables in general pants transform into global phases, so there is no multiplicity.

\end{exmp}

\begin{exmp}[$\defect{X_k}{X_l}{}{}{} \circ  \defect{X_l}{X_m}{}{}{}$]

	Consider the defect fusion
	\begin{align}
		\defect{X_k}{X_l}{}{}{} \circ \defect{X_l}{X_m}{}{}{}\to 
		\begin{cases} 
			\defect{X_k}{X_k}{\alpha}{\beta}{} & k = m\\
			\defect{X_k}{X_m}{}{}{} & k\neq m
		\end{cases}.
	\end{align}
	
\subsubsection*{Case I: $k=m$}
The general pants absorbing $\defect{X_l}{X_m}{}{}{}$ and $\defect{X_k}{X_l}{}{}{}$ on the legs and $\defect{X_k}{X_k}{\alpha}{\zeta}{}$ on the waist is
\begin{align}
  \frac{1}{p} \sum_{\gamma}\omega^{\gamma\zeta}
	\begin{array}{c}
		\includeTikz{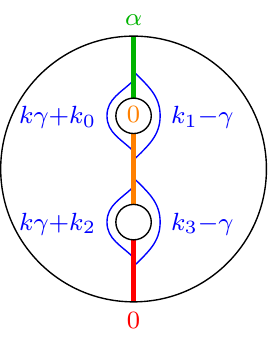}{
			\begin{tikzpicture}[scale=.9,,every node/.style={scale=.9}]
			  \vpantsparams{$0$}{$0$}{$\alpha$};
        \vpantsstp{$k\gamma{+}k_2$}{$k_3{-}\gamma$}{$k\gamma{+}k_0$}{$k_1{-}\gamma$};
			\end{tikzpicture}}
	\end{array}
\end{align}
where
\begin{align}
  k_0 + kk_1 &= \alpha \\
  k_0 + lk_1 - k_2 - lk_3 &= 0 \\
  k_2 + k k_3 &= 0.
\end{align}
The solution to these equations is
\begin{align*}
  k_1 &= k^{-1}(\alpha-k_0) \\
  k_2 &= l\alpha(k-l)^{-1}+k_0 \\
  k_3 &= - l\alpha(k-l)^{-1}k^{-1}-k_0k^{-1}
\end{align*}
so the substitution $\gamma \to \gamma - k^{-1}t$ transforms $k_0$ into a global phase. Therefore
\begin{align}
 \defect{X_k}{X_l}{}{}{}  \circ  \defect{X_l}{X_k}{}{}{}  = \oplus_{\alpha,\zeta}\defect{X_k}{X_k}{\alpha}{\zeta}{}
\end{align}
\subsubsection*{Case II: $k\not=m$}
The general pant absorbing $\defect{X_l}{X_m}{}{}{}$ and $\defect{X_k}{X_l}{}{}{}$ on the legs and $\defect{X_k}{X_m}{}{}{}$ on the waist is
\begin{align}
	\begin{array}{c}
		\includeTikz{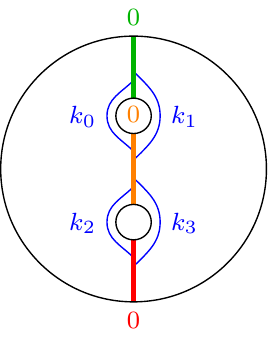}{
			\begin{tikzpicture}[scale=.9,,every node/.style={scale=.9}]
			  \vpantsparams{$0$}{$0$}{$0$};
        \vpantsstp{$k_2$}{$k_3$}{$k_0$}{$k_1$};
			\end{tikzpicture}}
	\end{array}
\end{align}
where
\begin{align}
  k_0 + mk_1 &= 0 \\
  k_0 + lk_1 - k_2 - lk_3 &= 0 \\
  k_2 + k k_3 &= 0.
\end{align}
This system has rank 3, which implies
\begin{align}
 \defect{X_k}{X_l}{}{}{}  \circ  \defect{X_l}{X_m}{}{}{}  = p\cdot\defect{X_k}{X_m}{}{}{}
\end{align}
\end{exmp}

\section{Physical interpretation of defects} \label{sec:physical_interpretations}
\begin{table}
	\resizebox{.9\textwidth}{!}{
		\begin{tabular}{!{\vrule width 1pt}c!{\vrule width 1pt}c|c!{\vrule width 1pt}}
			\toprule[1pt]
				\rowcolor[gray]{.9}[\tabcolsep]Bimodule label & Domain wall & Action on particles\\
			\toprule[1pt]
			$T$&$\begin{array}{c}\includeTikz{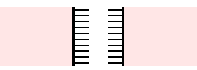}{
				\begin{tikzpicture}[yscale=.3]
				\draw[white](0,-1.1)--(0,1.21);
				\fill[red!10](-1,-1) rectangle (-.25,1);
				\fill[red!10](1,-1) rectangle (.25,1);
				\draw[thick](-.25,-1)--(-.25,1);
				\draw[thick](.25,-1)--(.25,1);
				\foreach \x in {0,...,9}{\draw (.1,-.9+.2*\x)--(.25,-.9+.2*\x);};
				\foreach \x in {0,...,9}{\draw (-.1,-.9+.2*\x)--(-.25,-.9+.2*\x);};
				\end{tikzpicture}}
			\end{array}$&Condenses $e$ on both sides\\
			$L$&$\begin{array}{c}\includeTikz{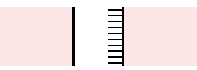}{\begin{tikzpicture}[yscale=.3]
				\draw[white](0,-1.1)--(0,1.21);
				\fill[red!10](-1,-1) rectangle (-.25,1);
				\fill[red!10](1,-1) rectangle (.25,1);
				\draw[thick](-.25,-1)--(-.25,1);
				\draw[thick](.25,-1)--(.25,1);
				\foreach \x in {0,...,9}{\draw (.1,-.9+.2*\x)--(.25,-.9+.2*\x);};
				\end{tikzpicture}}\end{array}$&Condenses $m$ on left and $e$ on right\\
			$R$&$\begin{array}{c}\includeTikz{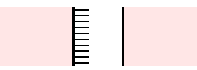}{\begin{tikzpicture}[yscale=.3]
				\draw[white](0,-1.1)--(0,1.21);
				\fill[red!10](-1,-1) rectangle (-.25,1);
				\fill[red!10](1,-1) rectangle (.25,1);
				\draw[thick](-.25,-1)--(-.25,1);
				\draw[thick](.25,-1)--(.25,1);
				\foreach \x in {0,...,9}{\draw (-.1,-.9+.2*\x)--(-.25,-.9+.2*\x);};
				\end{tikzpicture}}\end{array}$&Condenses $e$ on left and $m$ on right\\
			$F_0$&$\begin{array}{c}\includeTikz{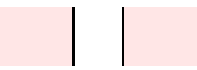}{\begin{tikzpicture}[yscale=.3]
				\draw[white](0,-1.1)--(0,1.21);
				\fill[red!10](-1,-1) rectangle (-.25,1);
				\fill[red!10](1,-1) rectangle (.25,1);
				\draw[thick](-.25,-1)--(-.25,1);
				\draw[thick](.25,-1)--(.25,1);
				\end{tikzpicture}}\end{array}$&Condenses $m$ on both sides\\
			\toprule[1pt]
			$X_k$&$\begin{array}{c}\includeTikz{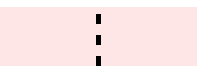}{\begin{tikzpicture}[yscale=.3]
				\draw[white](0,-1.1)--(0,1.21);
				\fill[red!10](-1,-1) rectangle (1,1);
				\draw[ultra thick,dashed] (0,-1)--(0,1);
				\end{tikzpicture}}\end{array}$&$X_k:e^am^b\mapsto e^{ka}m^{k^{-1} b}$ (moving left to right), where $k^{-1}$ is taken multiplicatively modulo $p$\\
			$F_{q}=F_1 X_q$&$\begin{array}{c}\includeTikz{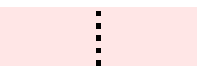}{\begin{tikzpicture}[yscale=.3]
				\draw[white](0,-1.1)--(0,1.21);
				\fill[red!10](-1,-1) rectangle (1,1);
				\draw[ultra thick,dotted] (0,-1)--(0,1);\end{tikzpicture}}\end{array}$&
			$F_1:e^am^b\mapsto e^{b}m^{a}$\\
			\toprule[1pt]
		\end{tabular}
	}
	\caption{Domain walls on the lattice corresponding to bimodules. Reproduced from \onlinecite{1806.01279}.}\label{tab:domainwallinterp}
\end{table}

Associated to a fusion category $\mathcal{C}$ is a topological phase described by its Drinfeld center $Z(\mathcal{C})$. The Levin-Wen procedure\cite{Levin2005} can be used to construct a lattice model which realizes the topological phase in its low energy space. In \onlinecite{1806.01279}, we showed how the bimodule categories for $\vvec{\ZZ{p}}$ can be interpreted in terms of boundaries and domain walls in the physical theory. This data is reproduced in Table~\ref{tab:domainwallinterp} for completeness. 

In this section, we discuss a physical interpretation for all defects in Table~\ref{tab:idempotents}. We will also show how the fusion rules in Tables~\ref{tab:horizontal_table_1}-\ref{tab:horizontal_table_6} and \ref{tab:vertical_fusion_tables} can be obtained from the physical theory up to multiplicity. The multiplicities remain mysterious from the physical perspective, but they can be computed from the Frobenius-Perron dimensions of the defects.

The simplest defect to interpret is $\defect{X_1}{X_1}{a}{x}{}$. By studying Table~\ref{tab:zptable} we observe that $X_1$ is the `identity' domain wall, corresponding to $\vvec{\ZZ{p}}$ as a bimodule over itself. The only defects of such a domain wall are the anyons of the theory. Therefore
\begin{align}
\defect{X_1}{X_1}{a}{x}{}\mapsto e^xm^a,
\end{align}
where the $p^2$ particles of the (domain wall free) theory associated to $\vvec{\ZZ{p}}$ are usually denoted $e^xm^a$, with $a,x\in \ZZ{p}$\cite{MR1951039}. 

For other defects, the physical interpretation is found by studying how anyons can be introduced to obtain distinct states. For example, consider the defects $\defect{X_k}{X_l}{\cdot}{}{}$. If an anyon $m^a$ is introduced on the left of the wall and $e^x$ is introduced on the right, they can be split and passed through the upper and lower walls to find the equivalent states
\begin{align}
	\begin{array}{c}
		\includeTikz{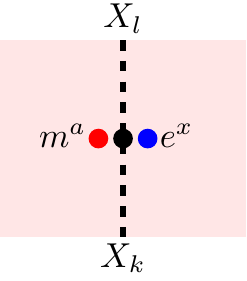}{
			\begin{tikzpicture}
					\fill[red!10] (-1.25,-1) rectangle (1.25,1);
					\draw[ultra thick,dashed](0,1)--(0,-1);
					\fill[black] (0,0) circle (.1) node[right,inline text,text=black] {};
					\node[below,inline text,text=black]  at (0,-1) {$X_k$};
					\node[above,inline text,text=black]  at (0,1) {$X_l$};
					\fill[blue] (.25,0) circle (.1) node[right,inline text,text=black] {$e^{x}$};
					\fill[red] (-.25,0) circle (.1) node[left,inline text,text=black] {$m^{a}$};
			\end{tikzpicture}
			}
	\end{array}
	=
	\begin{array}{c}
		\includeTikz{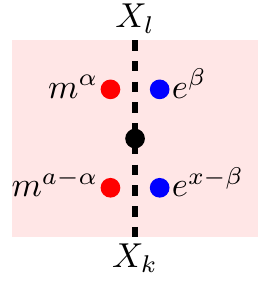}{
			\begin{tikzpicture}
			\fill[red!10] (-1.25,-1) rectangle (1.25,1);
			\draw[ultra thick,dashed](0,1)--(0,-1);
			\fill[black] (0,0) circle (.1) node[right,inline text,text=black] {};
					\node[below,inline text,text=black]  at (0,-1) {$X_k$};
					\node[above,inline text,text=black]  at (0,1) {$X_l$};
			\fill[red] (-.25,.5) circle (.1) node[left,inline text,text=black] {$m^{\alpha}$};
			\fill[red] (-.25,-.5) circle (.1) node[left,inline text,text=black] {$m^{a-\alpha}$};
			\fill[blue] (.25,.5) circle (.1) node[right,inline text,text=black] {$e^{\beta}$};
			\fill[blue] (.25,-.5) circle (.1) node[right,inline text,text=black] {$e^{x-\beta}$};
			\end{tikzpicture}
		}
	\end{array}
	=
	\begin{array}{c}
		\includeTikz{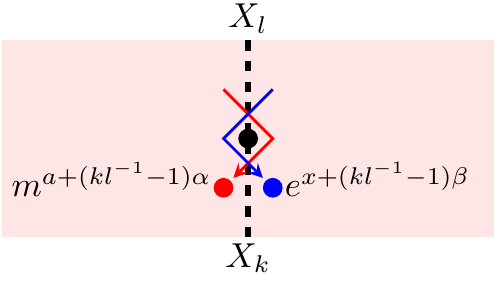}{
			\begin{tikzpicture}
			\fill[red!10] (-2.5,-1) rectangle (2.5,1);
			\draw[ultra thick,dashed](0,1)--(0,-1);
			\fill[black] (0,0) circle (.1) node[right,inline text,text=black] {};
					\node[below,inline text,text=black]  at (0,-1) {$X_k$};
					\node[above,inline text,text=black]  at (0,1) {$X_l$};
			\draw[thick,red,>=stealth,->] (-.25,.5)--(.25,0)->(-.15,-.4);
			\draw[thick,blue,>=stealth,->] (.25,.5)--(-.25,0)->(.15,-.4);
			\fill[red] (-.25,-.5) circle (.1) node[left,inline text,text=black] {$m^{a+(kl^{-1}{-}1)\alpha}$};
			\fill[blue] (.25,-.5) circle (.1) node[right,inline text,text=black] {$e^{x+(kl^{-1}{-}1)\beta}$};
			\end{tikzpicture}
		}
	\end{array}.\label{eqn:XkXldefectstates}
\end{align}
If $k=l$, then all such states are distinct, so there are $p^2$ defects, corresponding to a basic defect with anyons pushed onto it. For $k\neq l$, we are free to choose $\alpha$ and $\beta$ in Eqn.~\ref{eqn:XkXldefectstates}, so all states are equivalent, and there is a unique defect.

A very similar computation can be performed for other domain wall interfaces. In the case of $\defect{F_q}{X_l}{\cdot}{}{}$, the computation is
\begin{align}
	\begin{array}{c}
		\includeTikz{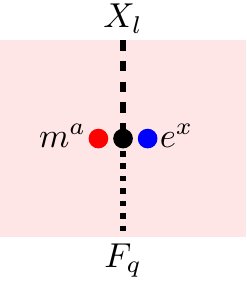}{
			\begin{tikzpicture}
			\fill[red!10] (-1.25,-1) rectangle (1.25,1);
			\draw[ultra thick,dashed](0,1)--(0,0);
			\draw[ultra thick,dotted](0,0)--(0,-1);
			\fill[black] (0,0) circle (.1) node[right,inline text,text=black] {};
			\node[below,inline text,text=black]  at (0,-1) {$F_q$};
			\node[above,inline text,text=black]  at (0,1) {$X_l$};
			\fill[blue] (.25,0) circle (.1) node[right,inline text,text=black] {$e^{x}$};
			\fill[red] (-.25,0) circle (.1) node[left,inline text,text=black] {$m^{a}$};
			\end{tikzpicture}
		}
	\end{array}
	=
	\begin{array}{c}
		\includeTikz{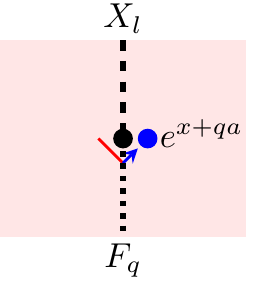}{
			\begin{tikzpicture}
			\fill[red!10] (-1.25,-1) rectangle (1.25,1);
			\fill[black] (0,0) circle (.1) node[right,inline text,text=black] {};
			\node[below,inline text,text=black]  at (0,-1) {$F_q$};
			\node[above,inline text,text=black]  at (0,1) {$X_l$};
			\draw[thick,red](-.25,0)--(0,-.25);\draw[thick,blue,>=stealth,->] (0,-.25)--(.15,-.1);
			\fill[blue] (.25,0) circle (.1) node[right,inline text,text=black] {$e^{x+qa}$};
						\draw[ultra thick,dashed](0,1)--(0,0);
						\draw[ultra thick,dotted](0,0)--(0,-1);
			\end{tikzpicture}
		}
	\end{array},\label{eqn:FqXldefectstates}
\end{align}
with no further equivalence possible. Therefore, there are $p$ distinct states corresponding to a base defect with a number of $e$ particles absorbed from the right. We remark that, just as in the choice of idempotents, there is a choice of labeling. We could have instead chosen to label by some combination of $e$ and $m$ on the left and right. This corresponds to permuting the labels.

When there is a boundary present, particles can be discarded (condensed) into the boundary. In the case of a $\defect{L}{T}{\cdot}{}{}$ defect, we have
\begin{align}
	\begin{array}{c}
		\includeTikz{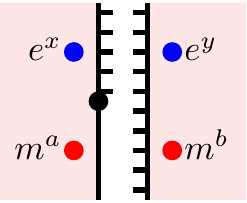}{
			\begin{tikzpicture}
				\begin{scope}[xscale=1]
				\fill[red!10](-1.25,-1) rectangle (-.25,1);
				\draw[ultra thick] (-.25,-1)--(-.25,1);
				\foreach \x in {5,...,9}{\draw[ultra thick] (-.1,-.9+.2*\x)--(-.25,-.9+.2*\x);};
				\end{scope}
				\begin{scope}[xscale=-1]
				\fill[red!10](-1.25,-1) rectangle (-.25,1);
				\draw[ultra thick] (-.25,-1)--(-.25,1);
				\foreach \x in {0,...,9}{\draw[ultra thick] (-.1,-.9+.2*\x)--(-.25,-.9+.2*\x);};
				\end{scope}
				\fill[black] (-.25,0) circle (.1) node[left,inline text,text=black] {};		
				\fill[red] (-.5,-.5) circle (.1) node[left,inline text,text=black] {$m^a$};
				\fill[red] (.5,-.5) circle (.1) node[right,inline text,text=black] {$m^b$};
				\fill[blue] (-.5,.5) circle (.1) node[left,inline text,text=black] {$e^x$};
				\fill[blue] (.5,.5) circle (.1) node[right,inline text,text=black] {$e^y$};
			\end{tikzpicture}
		}
	\end{array}
	=
	\begin{array}{c}
		\includeTikz{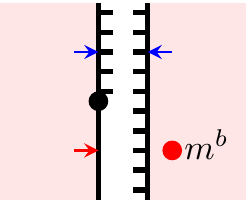}{
			\begin{tikzpicture}
				\begin{scope}[xscale=1]
				\fill[red!10](-1.25,-1) rectangle (-.25,1);
				\draw[ultra thick] (-.25,-1)--(-.25,1);
				\foreach \x in {5,...,9}{\draw[ultra thick] (-.1,-.9+.2*\x)--(-.25,-.9+.2*\x);};
				\end{scope}
				\begin{scope}[xscale=-1]
				\fill[red!10](-1.25,-1) rectangle (-.25,1);
				\draw[ultra thick] (-.25,-1)--(-.25,1);
				\foreach \x in {0,...,9}{\draw[ultra thick] (-.1,-.9+.2*\x)--(-.25,-.9+.2*\x);};
				\end{scope}
				\fill[black] (-.25,0) circle (.1) node[left,inline text,text=black] {};		
				\fill[red] (.5,-.5) circle (.1) node[right,inline text,text=black] {$m^b$};
				\draw[thick,red,>=stealth,->] (-.5,-.5)--(-.25,-.5);
				\draw[thick,blue,>=stealth,->] (-.5,.5)--(-.25,.5);
				\draw[thick,blue,>=stealth,->] (.5,.5)--(.25,.5);
			\end{tikzpicture}
		}
	\end{array},\label{eqn:LTdefectstates}
\end{align}
so there are $p$ distinct defects, corresponding to pinning a number of $m$ defects to the right hand boundary.

The full set of defect interpretations are listed in Table~\ref{tab:defectinterps}. 
We remark that in the case $\ZZ{2}$, many of the defects studied here are the subject of previous work, for example 
\begin{align}
	\defect{X_1}{F_1}{x}{}{},\defect{F_1}{X_1}{x}{}{}\mapsto \sigma_x
\end{align}
was referred to as a `twist' in \onlinecite{Bombin2010}, and defects involving rough/smooth interfaces in Table~\ref{tab:defectinterps} were named `corners' in \onlinecite{Brown2016}.

\begin{table}
	\resizebox{.9\textwidth}{!}{

			\\
			\toprule[1pt]
		\end{tabular}
	}
	\caption{Physical interpretation of all defects.}\label{tab:defectinterps}
\end{table}

\subsection{Fusing defects: Horizontal}
The physical interpretations from Table~\ref{tab:defectinterps} can be used to compute the fusion rules. In this section, we will illustrate how this is done using a few examples.

\begin{exmp}[$\defect{X_k}{X_k}{}{}{} \otimes \defect{X_l}{X_l}{}{}{}$]

Consider the fusion $\defect{X_k}{X_k}{a}{x}{}\otimes \defect{X_l}{X_l}{c}{z}{}$. In the physical theory, we begin by clearing any anyons occurring on the central region
\begin{align}
	\begin{array}{c}
		\includeTikz{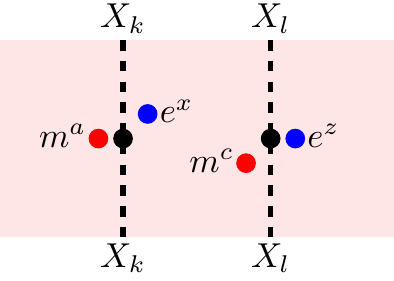}{
			\begin{tikzpicture}
					\fill[red!10] (-2,-1) rectangle (2,1);
					\draw[ultra thick,dashed](-.75,1)--(-.75,-1);
					\draw[ultra thick,dashed](.75,1)--(.75,-1);
					\fill[black] (-.75,0) circle (.1) node[right,inline text,text=black] {};
					\fill[black] (.75,0) circle (.1) node[right,inline text,text=black] {};
					\node[below,inline text,text=black]  at (-.75,-1) {$X_k$};
					\node[above,inline text,text=black]  at (-.75,1) {$X_k$};
					\node[below,inline text,text=black]  at (.75,-1) {$X_l$};
					\node[above,inline text,text=black]  at (.75,1) {$X_l$};
					\fill[blue] (-.5,.25) circle (.1) node[right,inline text,text=black] {$e^{x}$};
					\fill[red] (-1,0) circle (.1) node[left,inline text,text=black] {$m^{a}$};
					\fill[blue] (1,0) circle (.1) node[right,inline text,text=black] {$e^{z}$};
					\fill[red] (.5,-.25) circle (.1) node[left,inline text,text=black] {$m^{c}$};
			\end{tikzpicture}}
	\end{array}
	=
	\begin{array}{c}
		\includeTikz{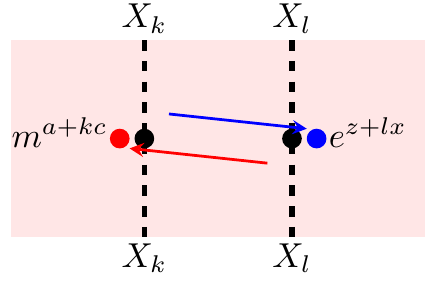}{
			\begin{tikzpicture}
			\fill[red!10] (-2.1,-1) rectangle (2.1,1);
			\draw[ultra thick,dashed](-.75,1)--(-.75,-1);
			\draw[ultra thick,dashed](.75,1)--(.75,-1);
			\fill[black] (-.75,0) circle (.1) node[right,inline text,text=black] {};
			\fill[black] (.75,0) circle (.1) node[right,inline text,text=black] {};
			\node[below,inline text,text=black]  at (-.75,-1) {$X_k$};
			\node[above,inline text,text=black]  at (-.75,1) {$X_k$};
			\node[below,inline text,text=black]  at (.75,-1) {$X_l$};
			\node[above,inline text,text=black]  at (.75,1) {$X_l$};
			\fill[red] (-1,0) circle (.1) node[left,inline text,text=black] {$m^{a+kc}$};
			\fill[blue] (1,0) circle (.1) node[right,inline text,text=black] {$e^{z+lx}$};
\draw[thick,blue,>=stealth,->] (-.5,.25)->(.9,.1);
\draw[thick,red,>=stealth,->] (.5,-.25)->(-.9,-.1);
			\end{tikzpicture}}
	\end{array}
		=
	\begin{array}{c}
		\includeTikz{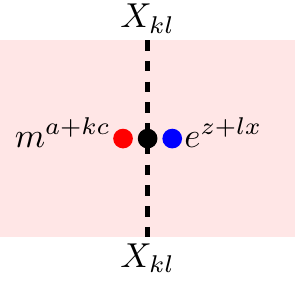}{
			\begin{tikzpicture}
			\fill[red!10] (-1.5,-1) rectangle (1.5,1);
			\draw[ultra thick,dashed](0,1)--(0,-1);
			\fill[black] (0,0) circle (.1) node[right,inline text,text=black] {};
			\node[below,inline text,text=black]  at (0,-1) {$X_{kl}$};
			\node[above,inline text,text=black]  at (0,1) {$X_{kl}$};
			\fill[red] (-.25,0) circle (.1) node[left,inline text,text=black] {$m^{a+kc}$};
			\fill[blue] (.25,0) circle (.1) node[right,inline text,text=black] {$e^{z+lx}$};
			\end{tikzpicture}}
	\end{array}.
\end{align}
This recovers the fusion computed using the annular algebra.
\end{exmp}

\begin{exmp}[$\defect{X_k}{X_l}{}{}{} \otimes \defect{X_m}{X_n}{}{}{}$]
A more interesting computation arises when we allow the upper and lower domain walls to differ
\begin{align}
	\begin{array}{c}
		\includeTikz{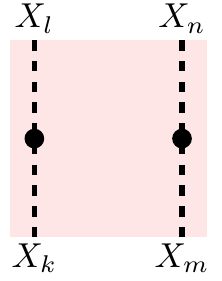}{
			\begin{tikzpicture}
			\fill[red!10] (-1,-1) rectangle (1,1);
			\draw[ultra thick,dashed](-.75,1)--(-.75,-1);
			\draw[ultra thick,dashed](.75,1)--(.75,-1);
			\fill[black] (-.75,0) circle (.1) node[right,inline text,text=black] {};
			\fill[black] (.75,0) circle (.1) node[right,inline text,text=black] {};
			\node[below,inline text,text=black]  at (-.75,-1) {$X_k$};
			\node[above,inline text,text=black]  at (-.75,1) {$X_l$};
			\node[below,inline text,text=black]  at (.75,-1) {$X_m$};
			\node[above,inline text,text=black]  at (.75,1) {$X_n$};
			\end{tikzpicture}}
	\end{array}
	=
	\begin{array}{c}
		\includeTikz{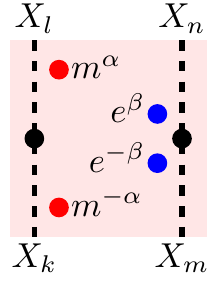}{
			\begin{tikzpicture}
			\fill[red!10] (-1,-1) rectangle (1,1);
			\draw[ultra thick,dashed](-.75,1)--(-.75,-1);
			\draw[ultra thick,dashed](.75,1)--(.75,-1);
			\fill[black] (-.75,0) circle (.1) node[right,inline text,text=black] {};
			\fill[black] (.75,0) circle (.1) node[right,inline text,text=black] {};
			\node[below,inline text,text=black]  at (-.75,-1) {$X_k$};
			\node[above,inline text,text=black]  at (-.75,1) {$X_l$};
			\node[below,inline text,text=black]  at (.75,-1) {$X_m$};
			\node[above,inline text,text=black]  at (.75,1) {$X_n$};
						\fill[blue] (.5,.25) circle (.1) node[left,inline text,text=black] {$e^{\beta}$};
						\fill[blue] (.5,-.25) circle (.1) node[left,inline text,text=black] {$e^{-\beta}$};
						\fill[red] (-.5,.7) circle (.1) node[right,inline text,text=black] {$m^{\alpha}$};
						\fill[red] (-.5,-.7) circle (.1) node[right,inline text,text=black] {$m^{-\alpha}$};
			\end{tikzpicture}}
	\end{array}
	=
	\begin{array}{c}
		\includeTikz{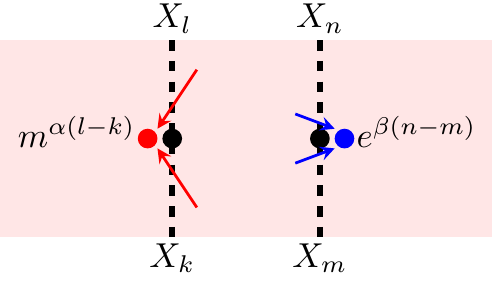}{
			\begin{tikzpicture}
			\fill[red!10] (-2.5,-1) rectangle (2.5,1);
			\draw[ultra thick,dashed](-.75,1)--(-.75,-1);
			\draw[ultra thick,dashed](.75,1)--(.75,-1);
			\fill[black] (-.75,0) circle (.1) node[right,inline text,text=black] {};
			\fill[black] (.75,0) circle (.1) node[right,inline text,text=black] {};
			\node[below,inline text,text=black]  at (-.75,-1) {$X_k$};
			\node[above,inline text,text=black]  at (-.75,1) {$X_l$};
			\node[below,inline text,text=black]  at (.75,-1) {$X_m$};
			\node[above,inline text,text=black]  at (.75,1) {$X_n$};
			\fill[blue] (1,0) circle (.1) node[right,inline text,text=black] {$e^{\beta(n-m)}$};
			\fill[red] (-1,0) circle (.1) node[left,inline text,text=black] {$m^{\alpha(l-k)}$};
			\draw[thick,blue,>=stealth,->] (.5,.25)->(.9,.1);\draw[thick,blue,>=stealth,->] (.5,-.25)->(.9,-.1);
			\draw[thick,red,>=stealth,->] (-.5,.7)->(-.9,.1);\draw[thick,red,>=stealth,->] (-.5,-.7)->(-.9,-.1);
			\end{tikzpicture}}
	\end{array}
	=
	\begin{array}{c}
		\includeTikz{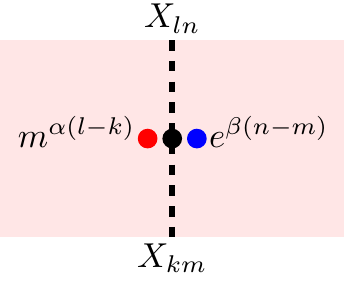}{
			\begin{tikzpicture}
			\fill[red!10] (-1.75,-1) rectangle (1.75,1);
			\draw[ultra thick,dashed](0,1)--(0,-1);
			\fill[black] (0,0) circle (.1) node[right,inline text,text=black] {};
			\node[below,inline text,text=black]  at (0,-1) {$X_{km}$};
			\node[above,inline text,text=black]  at (0,1) {$X_{ln}$};
			\fill[blue] (.25,0) circle (.1) node[right,inline text,text=black] {$e^{\beta(n-m)}$};
			\fill[red] (-.25,0) circle (.1) node[left,inline text,text=black] {$m^{\alpha(l-k)}$};
			\end{tikzpicture}}
	\end{array}.
\end{align}
If $km\neq ln$, then from Eqn.~\ref{eqn:XkXldefectstates}, this is equivalent to the base defect for all choices $\alpha,\beta$. In the case where $km=ln$, these states are distinct for each $\alpha, \beta$, giving the fusion rule
\begin{align}
	\defect{X_k}{X_l}{}{}{}\otimes \defect{X_m}{X_n}{}{}{}&=
	\begin{cases}
		p\cdot \defect{X_{km}}{X_{ln}}{}{}{}&km\neq ln\\
		\oplus_{\alpha,\beta} \defect{X_{km}}{X_{km}}{\alpha}{\beta}{}&km=ln
	\end{cases}.
\end{align}
The coefficient $p$ is not obtained from this computation, but can be calculated using the Frobenius-Perron dimensions computed in Section~\ref{S:FPd}.

\end{exmp}

\begin{exmp}[$\protect \defect{X_k}{F_r}{}{}{} \otimes \protect\defect{X_m}{F_t}{}{}{}$]

Consider the fusion process
\begin{align}
	\defect{X_k}{F_r}{x}{}{}\otimes\defect{X_m}{F_t}{z}{}{}.
\end{align}
In the physical theory, this is computed using
\begin{align}
	\begin{array}{c}
		\includeTikz{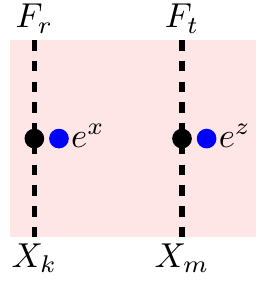}{
			\begin{tikzpicture}
			\fill[red!10] (-1,-1) rectangle (1.5,1);
			\draw[ultra thick,dashed](-.75,1)--(-.75,-1);
			\draw[ultra thick,dashed](.75,1)--(.75,-1);
			\fill[black] (-.75,0) circle (.1) node[right,inline text,text=black] {};
			\fill[black] (.75,0) circle (.1) node[right,inline text,text=black] {};
			\node[below,inline text,text=black]  at (-.75,-1) {$X_k$};
			\node[above,inline text,text=black]  at (-.75,1) {$F_r$};
			\node[below,inline text,text=black]  at (.75,-1) {$X_m$};
			\node[above,inline text,text=black]  at (.75,1) {$F_t$};
			\fill[blue] (-.5,0) circle (.1) node[right,inline text,text=black] {$e^{x}$};
			\fill[blue] (1,0) circle (.1) node[right,inline text,text=black] {$e^{z}$};
			\end{tikzpicture}}
	\end{array}
	=
	\begin{array}{c}
		\includeTikz{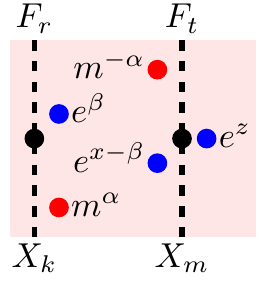}{
			\begin{tikzpicture}
			\fill[red!10] (-1,-1) rectangle (1.5,1);
			\draw[ultra thick,dashed](-.75,1)--(-.75,-1);
			\draw[ultra thick,dashed](.75,1)--(.75,-1);
			\fill[black] (-.75,0) circle (.1) node[right,inline text,text=black] {};
			\fill[black] (.75,0) circle (.1) node[right,inline text,text=black] {};
			\node[below,inline text,text=black]  at (-.75,-1) {$X_k$};
			\node[above,inline text,text=black]  at (-.75,1) {$F_r$};
			\node[below,inline text,text=black]  at (.75,-1) {$X_m$};
			\node[above,inline text,text=black]  at (.75,1) {$F_t$};
			\fill[blue] (-.5,.25) circle (.1) node[right,inline text,text=black] {$e^{\beta}$};
			\fill[blue] (.5,-.25) circle (.1) node[left,inline text,text=black] {$e^{x-\beta}$};
			\fill[red] (.5,.7) circle (.1) node[left,inline text,text=black] {$m^{-\alpha}$};
			\fill[red] (-.5,-.7) circle (.1) node[right,inline text,text=black] {$m^{\alpha}$};
			\fill[blue] (1,0) circle (.1) node[right,inline text,text=black] {$e^{z}$};
			\end{tikzpicture}}
	\end{array}
	=
	\begin{array}{c}
		\includeTikz{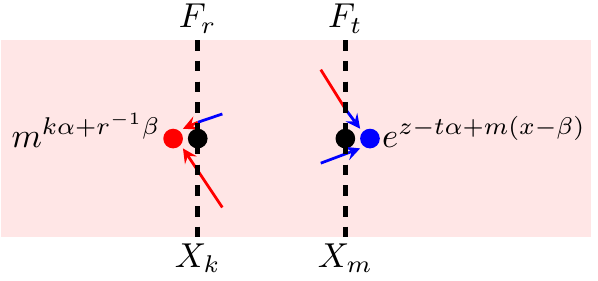}{
			\begin{tikzpicture}
			\fill[red!10] (-2.75,-1) rectangle (3.25,1);
			\draw[ultra thick,dashed](-.75,1)--(-.75,-1);
			\draw[ultra thick,dashed](.75,1)--(.75,-1);
			\fill[black] (-.75,0) circle (.1) node[right,inline text,text=black] {};
			\fill[black] (.75,0) circle (.1) node[right,inline text,text=black] {};
			\node[below,inline text,text=black]  at (-.75,-1) {$X_k$};
			\node[above,inline text,text=black]  at (-.75,1) {$F_r$};
			\node[below,inline text,text=black]  at (.75,-1) {$X_m$};
			\node[above,inline text,text=black]  at (.75,1) {$F_t$};
			\fill[blue] (1,0) circle (.1) node[right,inline text,text=black] {$e^{z-t\alpha+m(x-\beta)}$};
			\draw[thick,blue] (-.5,.25)--(-.75,.165);\draw[thick,red,>=stealth,->] (-.75,.165)->(-.9,.1);
			\draw[thick,red,>=stealth,->] (-.5,-.7)->(-.9,-.1);
			\draw[thick,blue,>=stealth,->] (.5,-.25)->(.9,-.1);
			\draw[thick,red] (.5,.7)--(.75,.3);\draw[thick,blue,>=stealth,->] (.75,.3)->(.9,.1);
			\draw[ultra thick,dashed](-.75,1)--(-.75,-1);
			\draw[ultra thick,dashed](.75,1)--(.75,-1);
			\fill[red] (-1,0) circle (.1) node[left,inline text,text=black] {$m^{k\alpha+r^{-1}\beta}$};
			\end{tikzpicture}}
	\end{array}
	=
	\begin{array}{c}
		\includeTikz{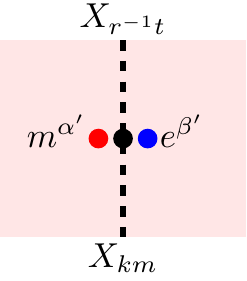}{
			\begin{tikzpicture}
			\fill[red!10] (-1.25,-1) rectangle (1.25,1);
			\draw[ultra thick,dashed](0,1)--(0,-1);
			\fill[black] (0,0) circle (.1) node[right,inline text,text=black] {};
			\node[below,inline text,text=black]  at (0,-1) {$X_{km}$};
			\node[above,inline text,text=black]  at (0,1) {$X_{r^{-1}t}$};
			\fill[blue] (.25,0) circle (.1) node[right,inline text,text=black] {$e^{\beta^\prime}$};
			\fill[red] (-.25,0) circle (.1) node[left,inline text,text=black] {$m^{\alpha^\prime}$};
			\end{tikzpicture}}
	\end{array},
\end{align}
where $\alpha^\prime=k\alpha+r^{-1}\beta$ and $\beta^\prime=z+mx-mr\alpha^\prime+(mrk-t)\alpha$. If $km\neq r^{-1}t$, these states are all equivalent to the base defect following Eqn.~\ref{eqn:XkXldefectstates}. In the case $km=r^{-1}t$, we have
\begin{align}
		\begin{array}{c}
			\includeTikz{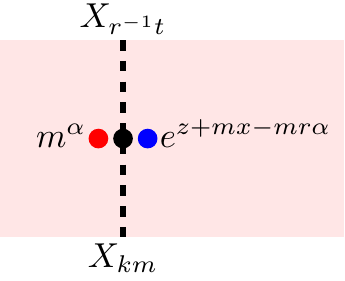}{
				\begin{tikzpicture}
				\fill[red!10] (-1.25,-1) rectangle (2.25,1);
				\draw[ultra thick,dashed](0,1)--(0,-1);
				\fill[black] (0,0) circle (.1) node[right,inline text,text=black] {};
				\node[below,inline text,text=black]  at (0,-1) {$X_{km}$};
				\node[above,inline text,text=black]  at (0,1) {$X_{r^{-1}t}$};
				\fill[blue] (.25,0) circle (.1) node[right,inline text,text=black] {$e^{z+mx-mr\alpha}$};
				\fill[red] (-.25,0) circle (.1) node[left,inline text,text=black] {$m^{\alpha}$};
				\end{tikzpicture}}
		\end{array},
\end{align}
for any $\alpha$. The fusion rule is therefore
\begin{align}
	\defect{X_k}{F_r}{x}{}{}\otimes \defect{X_m}{Ft}{z}{}{}&=
	\begin{cases}
		\defect{X_{km}}{X_{r^{-1}t}}{}{}{}&km\neq r^{-1}t\\
		\oplus_{\alpha} \defect{X_{km}}{X_{km}}{\alpha}{z+m(x-r\alpha)}{}&km=r^{-1}t
	\end{cases},
\end{align}
where the Frobenius-Perron dimensions (Section~\ref{S:FPd}) can be used to check that no multiplicity is required.

\end{exmp}

\begin{exmp}[$\defect{T}{T}{}{}{} \otimes \defect{T}{T}{}{}{}$]

Multiplicity in the domain wall fusion (Table~\ref{tab:zptable}) has an effect on the defect fusion. For the fusion
\begin{align}
	\defect{T}{T}{a}{b}{}\otimes \defect{T}{T}{c}{d}{},
\end{align}
the physical pre-fusion state is
\begin{align}
	\begin{array}{c}
		\includeTikz{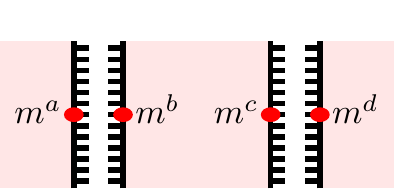}
		{
			\begin{tikzpicture}[yscale=.75]
			\draw[thick,white](1,0)--(1,1) node[above,inline text,pos=1]{\footnotesize$m^{b+c}$};
			\begin{scope}[xscale=1]
			\fill[red!10](-1,-1) rectangle (-.25,1);
			\draw[ultra thick] (-.25,-1)--(-.25,1);
			\foreach \x in {0,...,12}
			{
				\draw [ultra thick] (-.25,-1+.1+\x*.15)--(-.1,-1+.1+\x*.15);
			}
			\end{scope}
			\begin{scope}[xscale=-1]
			\fill[red!10](-1,-1) rectangle (-.25,1);
			\draw[ultra thick] (-.25,-1)--(-.25,1);
			\foreach \x in {0,...,12}
			{
				\draw [ultra thick] (-.25,-1+.1+\x*.15)--(-.1,-1+.1+\x*.15);
			}
			\end{scope}
			\fill[red,text=black] (-.25,0) circle (.1) node[left,inline text] {$m^a$};
			\fill[red,text=black] (.25,0) circle (.1) node[right,inline text] {$m^b$};
			\begin{scope}[shift={(2,0)}]
			\begin{scope}[xscale=1]
			\fill[red!10](-1,-1) rectangle (-.25,1);
			\draw[ultra thick] (-.25,-1)--(-.25,1);
			\foreach \x in {0,...,12}
			{
				\draw [ultra thick] (-.25,-1+.1+\x*.15)--(-.1,-1+.1+\x*.15);
			}
			\end{scope}
			\begin{scope}[xscale=-1]
			\fill[red!10](-1,-1) rectangle (-.25,1);
			\draw[ultra thick] (-.25,-1)--(-.25,1);
			\foreach \x in {0,...,12}
			{
				\draw [ultra thick] (-.25,-1+.1+\x*.15)--(-.1,-1+.1+\x*.15);
			}
			\end{scope}
			\fill[red,text=black] (-.25,0) circle (.1) node[left,inline text] {$m^c$};
			\fill[red,text=black] (.25,0) circle (.1) node[right,inline text] {$m^d$};
			\end{scope}
			\end{tikzpicture}
		}
	\end{array}.
\end{align}
The central strip supports a $p$ dimensional vector space. The qudit state can be read out by exchanging an $e$ particle between the boundaries. The state is changed by inserting an $m$ line vertically. Suppose the strip is in the state $m^\mu$. To perform the fusion, we push the inner $m$ particles away from the fusion region
\begin{align}
	\begin{array}{c}
		\includeTikz{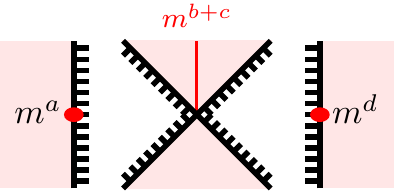}
		{
			\begin{tikzpicture}[yscale=.75]
			\filldraw[red!10](.25,-1)--(1.75,-1)--(1,0)--cycle;
			\filldraw[yscale=-1,red!10](.25,-1)--(1.75,-1)--(1,0)--cycle;
			\draw[thick,red](1,0)--(1,1) node[above,inline text,pos=1]{\footnotesize$m^{b+c}$};
			\begin{scope}[xscale=1]
			\fill[red!10](-1,-1) rectangle (-.25,1);
			\draw[ultra thick] (-.25,-1)--(-.25,1);
			\foreach \x in {0,...,12}
			{
				\draw [ultra thick] (-.25,-1+.1+\x*.15)--(-.1,-1+.1+\x*.15);
			}
			\end{scope}
			\fill[red,text=black] (-.25,0) circle (.1) node[left,inline text] {$m^a$};
			\begin{scope}[shift={(2,0)}]
			\begin{scope}[xscale=-1]
			\fill[red!10](-1,-1) rectangle (-.25,1);
			\draw[ultra thick] (-.25,-1)--(-.25,1);
			\foreach \x in {0,...,12}
			{
				\draw [ultra thick] (-.25,-1+.1+\x*.15)--(-.1,-1+.1+\x*.15);
			}
			\end{scope}
			\fill[red,text=black] (.25,0) circle (.1) node[right,inline text] {$m^d$};
			\end{scope}
			\draw[ultra thick] (.25,-1) --(1.75,1) (.25,1) --(1.75,-1);
			\foreach \x in {0,...,8}
			{
				\draw [ultra thick,shift={(.075+3/40*\x,.1+.1*\x)}] (.25,-1)--(.25-.075,-1+.1);
				\draw [yscale=-1,ultra thick,shift={(.075+3/40*\x,.1+.1*\x)}] (.25,-1)--(.25-.075,-1+.1);
				\draw [xscale=-1,ultra thick,shift={(-2+.075+3/40*\x,.1+.1*\x)}] (.25,-1)--(.25-.075,-1+.1);
				\draw [xscale=-1,yscale=-1,ultra thick,shift={(-2+.075+3/40*\x,.1+.1*\x)}] (.25,-1)--(.25-.075,-1+.1);
			}
			\end{tikzpicture}
		}
	\end{array}.
\end{align}
After fusing the central strip, there are still $p$ states. These are now understood as a subspace of a 2 qudit system, with one supported on each of the upper and lower regions. The subspace is spanned by states of the form $\ket{m^\mu}\otimes\ket{m^{\mu+b+c}}$. The fusion outcome is therefore
\begin{align}
	\left[ \defect{T}{T}{a}{b}{}\otimes \defect{T}{T}{c}{d}{} \right]_{\mu,\nu} = \delta_{\nu-\mu}^{b+c} \cdot \defect{T}{T}{a}{d}{}.
\end{align}
\end{exmp}

\begin{exmp}[$\defect{L}{T}{}{}{} \otimes \defect{L}{T}{}{}{}$]
As a final example, we show how the fusion
\begin{align}
	\defect{L}{T}{a}{}{}\otimes \defect{L}{T}{b}{}{}
\end{align}
can be computed.

Physically, this process is represented by
\begin{align}
	\begin{array}{c}
		\includeTikz{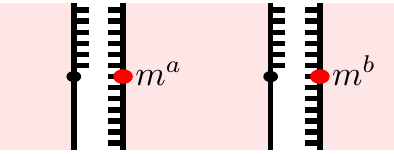}
		{
			\begin{tikzpicture}[yscale=.75]
			\begin{scope}[xscale=1]
			\fill[red!10](-1,-1) rectangle (-.25,1);
			\draw[ultra thick] (-.25,-1)--(-.25,1);
			\foreach \x in {7,...,12}
			{
				\draw [ultra thick] (-.25,-1+.1+\x*.15)--(-.1,-1+.1+\x*.15);
			}
			\end{scope}
			\begin{scope}[xscale=-1]
			\fill[red!10](-1,-1) rectangle (-.25,1);
			\draw[ultra thick] (-.25,-1)--(-.25,1);
			\foreach \x in {0,...,12}
			{
				\draw [ultra thick] (-.25,-1+.1+\x*.15)--(-.1,-1+.1+\x*.15);
			}
			\end{scope}
			\fill[black,text=black] (-.25,0) circle (.075);
			\fill[red,text=black] (.25,0) circle (.1) node[right,inline text] {$m^a$};
			\begin{scope}[shift={(2,0)}]
			\begin{scope}[xscale=1]
			\fill[red!10](-1,-1) rectangle (-.25,1);
			\draw[ultra thick] (-.25,-1)--(-.25,1);
			\foreach \x in {7,...,12}
			{
				\draw [ultra thick] (-.25,-1+.1+\x*.15)--(-.1,-1+.1+\x*.15);
			}
			\end{scope}
			\begin{scope}[xscale=-1]
			\fill[red!10](-1,-1) rectangle (-.25,1);
			\draw[ultra thick] (-.25,-1)--(-.25,1);
			\foreach \x in {0,...,12}
			{
				\draw [ultra thick] (-.25,-1+.1+\x*.15)--(-.1,-1+.1+\x*.15);
			}
			\end{scope}
			\fill[black,text=black] (-.25,0) circle (.075);
			\fill[red,text=black] (.25,0) circle (.1) node[right,inline text] {$m^b$};
			\end{scope}
			\end{tikzpicture}
		}
	\end{array}.
\end{align}
The central strip supports a $p$ dimensional vector space. The qudit state can be read out by exchanging an $e$ particle between the boundaries in the upper region. The state is changed by inserting an $m$ line vertically. To perform the fusion, we must push the $m$ particles away from the fusion region
\begin{align}
	\begin{array}{c}
		\includeTikz{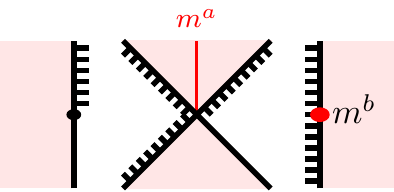}
		{
			\begin{tikzpicture}[yscale=.75]
			\filldraw[red!10](.25,-1)--(1.75,-1)--(1,0)--cycle;
			\filldraw[yscale=-1,red!10](.25,-1)--(1.75,-1)--(1,0)--cycle;
			\draw[thick,red](1,0)--(1,1) node[above,inline text,pos=1]{\footnotesize$m^a$};
			\begin{scope}[xscale=1]
			\fill[red!10](-1,-1) rectangle (-.25,1);
			\draw[ultra thick] (-.25,-1)--(-.25,1);
			\foreach \x in {7,...,12}
			{
				\draw [ultra thick] (-.25,-1+.1+\x*.15)--(-.1,-1+.1+\x*.15);
			}
			\end{scope}
			\fill[black,text=black] (-.25,0) circle (.075);
			\begin{scope}[shift={(2,0)}]
			\begin{scope}[xscale=-1]
			\fill[red!10](-1,-1) rectangle (-.25,1);
			\draw[ultra thick] (-.25,-1)--(-.25,1);
			\foreach \x in {0,...,12}
			{
				\draw [ultra thick] (-.25,-1+.1+\x*.15)--(-.1,-1+.1+\x*.15);
			}
			\end{scope}
			\fill[red,text=black] (.25,0) circle (.1) node[right,inline text] {$m^b$};
			\end{scope}
			\draw[ultra thick] (.25,-1) --(1.75,1) (.25,1) --(1.75,-1);
			\foreach \x in {0,...,8}
			{
				\draw [ultra thick,shift={(.075+3/40*\x,.1+.1*\x)}] (.25,-1)--(.25-.075,-1+.1);
				\draw [yscale=-1,ultra thick,shift={(.075+3/40*\x,.1+.1*\x)}] (.25,-1)--(.25-.075,-1+.1);
				\draw [xscale=-1,yscale=-1,ultra thick,shift={(-2+.075+3/40*\x,.1+.1*\x)}] (.25,-1)--(.25-.075,-1+.1);
			}
			\end{tikzpicture}
		}
	\end{array}.
\end{align}
After fusing the central strip, there are still $p$ states, fully supported on the upper region. The fusion outcome is therefore
\begin{align}
	\left[\defect{L}{T}{a}{}{}\times \defect{L}{T}{b}{}{}\right]_\mu=\defect{L}{T}{b}{}{}
\end{align}

\end{exmp}

\subsection{Fusing defects: Vertical}
Since the vertical fusions are much simpler than the horizontal ones, we will provide a single example to illustrate how the physical interpretation can be used for the fusion calculation. 
\begin{exmp}[$\defect{L}{F_r}{}{}{}\circ \defect{F_r}{L}{}{}{}$]
Consider the fusion
\begin{align}
	\defect{L}{F_r}{x}{}{}\circ \defect{F_r}{L}{z}{}{}. 
\end{align}
Physically, the fusion is
\begin{align}
	\begin{array}{c}
		\includeTikz{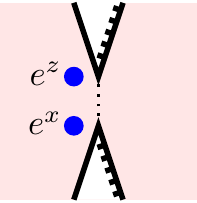}{
			\begin{tikzpicture}
			\fill[red!10] (-1,-1)rectangle(1,1);
			\filldraw[white,draw=black,ultra thick] (-.25,-1)--(0,-.25)--(.25,-1);
			\filldraw[white,draw=black,ultra thick] (-.25,1)--(0,.25)--(.25,1);
			\draw[thick,dotted](0,-.25)--(0,.25);
			\foreach \x in {0,...,4}{\draw[ultra thick,shift={(.25-.04*\x-.025,-1+.12*\x+.075)}] (0,0)--(-.075,-.025);};
			\foreach \x in {0,...,4}{\draw[ultra thick,shift={(.25-.04*\x-.025,1-.12*\x-.075)}] (0,0)--(-.075,.025);};
			\fill[blue] (-.25,-.25) circle (.1) node[left,inline text,text=black] {$e^{x}$};
			\fill[blue] (-.25,.25) circle (.1) node[left,inline text,text=black] {$e^{z}$};
			\end{tikzpicture}}
	\end{array}
=
	\begin{array}{c}
		\includeTikz{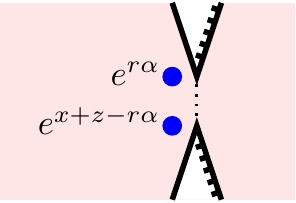}{
			\begin{tikzpicture}
			\fill[red!10] (-2,-1)rectangle(1,1);
			\filldraw[white,draw=black,ultra thick] (-.25,-1)--(0,-.25)--(.25,-1);
			\filldraw[white,draw=black,ultra thick] (-.25,1)--(0,.25)--(.25,1);
			\draw[thick,dotted](0,-.25)--(0,.25);
			\foreach \x in {0,...,4}{\draw[ultra thick,shift={(.25-.04*\x-.025,-1+.12*\x+.075)}] (0,0)--(-.075,-.025);};
			\foreach \x in {0,...,4}{\draw[ultra thick,shift={(.25-.04*\x-.025,1-.12*\x-.075)}] (0,0)--(-.075,.025);};
			\fill[blue] (-.25,-.25) circle (.1) node[left,inline text,text=black] {$e^{x+z-r\alpha}$};
			\fill[blue] (-.25,.25) circle (.1) node[left,inline text,text=black] {$e^{r\alpha}$};
			\end{tikzpicture}}
	\end{array}
=
	\begin{array}{c}
		\includeTikz{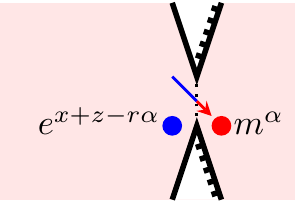}{
			\begin{tikzpicture}
			\fill[red!10] (-2,-1)rectangle(1,1);
			\filldraw[white,draw=black,ultra thick] (-.25,-1)--(0,-.25)--(.25,-1);
			\filldraw[white,draw=black,ultra thick] (-.25,1)--(0,.25)--(.25,1);
			\draw[thick,dotted](0,-.25)--(0,.25);
			\foreach \x in {0,...,4}{\draw[ultra thick,shift={(.25-.04*\x-.025,-1+.12*\x+.075)}] (0,0)--(-.075,-.025);};
			\foreach \x in {0,...,4}{\draw[ultra thick,shift={(.25-.04*\x-.025,1-.12*\x-.075)}] (0,0)--(-.075,.025);};
			\fill[blue] (-.25,-.25) circle (.1) node[left,inline text,text=black] {$e^{x+z-r\alpha}$};
			\fill[red] (.25,-.25) circle (.1) node[right,inline text,text=black] {$m^{\alpha}$};
			\draw[thick,blue](-.25,.25)--(0,0);
			\draw[thick,red,>=stealth,->] (0,0)--(.15,-.15);
			\end{tikzpicture}}
	\end{array}
=
	\begin{array}{c}
		\includeTikz{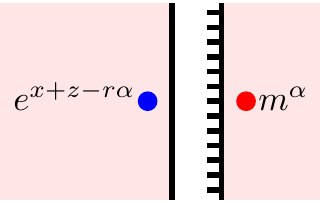}{
			\begin{tikzpicture}
			\fill[red!10] (-2,-1)rectangle(1.25,1);
			\filldraw[white] (-.25,-1) rectangle (.25,1);
			\draw[ultra thick](-.25,-1)--(-.25,1);\draw[ultra thick](.25,-1)--(.25,1);
			\foreach \x in {0,...,12}
			{
				\draw [ultra thick] (.25,-1+.1+\x*.15)--(.1,-1+.1+\x*.15);
			}
			\fill[blue] (-.5,0) circle (.1) node[left,inline text,text=black] {$e^{x+z-r\alpha}$};
			\fill[red] (.5,0) circle (.1) node[right,inline text,text=black] {$m^{\alpha}$};
			\end{tikzpicture}}
	\end{array}.
\end{align}
Since any $\alpha$ could have been chosen for this calculation, the fusion rule is
\begin{align}
	\defect{L}{F_r}{x}{}{}\circ\defect{F_r}{L}{z}{}{}=\oplus_\alpha \defect{L}{L}{\alpha}{x+z-r\alpha}{}.
\end{align}
\end{exmp}

\subsection{Frobenius-Perron Dimension}\label{S:FPd}
The Frobenius-Perron dimension (FPd) of the defects can be computed using the fusion table. For defects $a$, $b$ and $c$, the FPd obeys 
\begin{align}
d_a d_b=\sum_c N_{a,b}^c d_c,
\end{align}
where $N_{a,b}^c$ is the multiplicity of the fusion. The FPd of the defects of $\vvec{\ZZ{p}}$ do not depend on the defect label, only on the domain walls involved. The FPds are 
\begin{align}
\dim\left(\defect{X_k}{X_l}{\bullet}{}{}\right)&=\dim\left(\defect{F_k}{F_l}{\bullet}{}{}\right)=\begin{cases}1&k=l\\p&k\neq l\end{cases},\\
\dim\left(\defect{X_k}{F_r}{\bullet}{}{}\right)&=\dim\left(\defect{F_r}{X_k}{\bullet}{}{}\right)=\sqrt{p},\\
\dim\left(\defect{W_1}{W_2}{\bullet}{}{}\right)&=0&\text{for all other defects}.
\end{align}


\section{Natural Transformations} \label{sec:natural_transformations}

As explained in \onlinecite{MR2942952}, fusion category bimodules correspond to domain walls, and bimodule functors correspond to defects. On the mathematics side, we also have natural transformations. Using the diagrammatic framework from this paper, these natural transformations are easy to compute: they are just morphisms in the Karoubi envelope of the annular category $\ann{M}{N}{}$. Since the annular category is semi-simple, there are no morphisms between distinct objects in the Karoubi envelope and the endomorphism algebra of any simple object is just $\mathbb{C}$. Interesting natural transformations show up when we consider fusion. For example, in the vertical fusion
\begin{align}
  \defect{R}{L}{}{}{} \circ \defect{L}{X_m}{}{}{} = p \cdot \defect{R}{X_m}{}{}{}
\end{align}
There are $p$ distinct pants diagrams
\begin{align}
	\begin{array}{c}
		\includeTikz{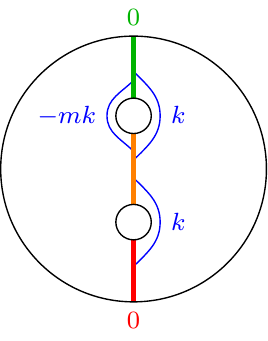}{
			\begin{tikzpicture}[scale=.9,,every node/.style={scale=.9}]
			  \vpantsparams{$0$}{}{$0$};
        \vpantsstp{}{$k$}{$-mk$}{$k$};
			\end{tikzpicture}}
	\end{array}
\end{align}
parameterized by $k$, which absorb $\defect{R}{L}{}{}{}, \defect{L}{X_m}{}{}{}$ on the legs and $\defect{R}{X_m}{}{}{}$ on the waist. These pants diagrams represent different natural transformations. For every defect fusion in this paper, the pants diagram which is computed represents a natural transformation.

A reasonable physical interpretation is that natural transformations capture certain aspects of the renormalization process. When computing defect fusion, we are witnessing an isomorphism between the horizontal or vertical concatenation of the defects and another defect. Physically, this corresponds to bringing the defects close together and then locally renormalizing, or \emph{zooming out}.


\section{Conclusions}

In this work, we have studied binary interface defects and their fusion. Using string diagrams and the annular category, we have classified the full set of defects occurring interfacing a pair of (not necessarily invertible) domain walls for the tensor category $\vvec{\ZZ{p}}$. Further, we have provided algorithms for computing both horizontal (tensor product) and vertical (composition) fusion of arbitrary pairs of defects. For the theory $\vvec{\ZZ{p}}$, we have provided complete fusion tables. 

We have specialized to $\vvec{\ZZ{p}}$ for simplicity. The framework outlined here is not restricted to this class of fusion categories. Of particular interest is the Color code\cite{Bombin2007a,Yoshida2015a} ($\vvec{\ZZ{2}\times\ZZ{2}}$) due to its importance in quantum computation. Defects between invertible domain walls and the trivial wall ($X_1$ in this paper) were studied in \onlinecite{Brown2018}, but the full theory is currently open. Additionally, the tools presented here are expected to be useful for studying nonabelian theories\cite{PhysRevB.96.195129}. Additionally, one could study domain walls and defects between distinct phases, such as the Color code and $\vvec{\ZZ{4}}$, which may prove useful for quantum computing tasks.
Although we restrict to binary interface defects, generalizations of the techniques developed here can be applied to higher defects such as those occurring at the interface of three or more domain walls. Such defects allow the meeting of many distinct topological phases.

In the physics literature, defect fusion is often synonymous with symmetry gauging\cite{MR2677836,MR2609644,Barkeshli2014,SETPaper,1804.01657}. In this work, we have computed the fusions without consideration of gauging. It would be extremely useful if the techniques developed in this work can say something about (obstructions to) gauging invertible defects.

\acknowledgments
This work is supported by the Australian Research Council (ARC) via Discovery Project ``Subfactors and symmetries'' DP140100732 and Discovery Project ``Low dimensional categories'' DP16010347. Support is also provided by the US Army Research Office for Basic Scientific Research via Grant W911NF-17-1-0401. We thank Christopher Chubb for feedback on the draft manuscript.

\vspace{5mm}

\bibliographystyle{apsrev_jacob}
\bibliography{refs}

\begin{thebibliography}{43}%
\makeatletter
\providecommand \@ifxundefined [1]{%
 \@ifx{#1\undefined}
}%
\providecommand \@ifnum [1]{%
 \ifnum #1\expandafter \@firstoftwo
 \else \expandafter \@secondoftwo
 \fi
}%
\providecommand \@ifx [1]{%
 \ifx #1\expandafter \@firstoftwo
 \else \expandafter \@secondoftwo
 \fi
}%
\providecommand \natexlab [1]{#1}%
\providecommand \enquote  [1]{``#1''}%
\providecommand \bibnamefont  [1]{#1}%
\providecommand \bibfnamefont [1]{#1}%
\providecommand \citenamefont [1]{#1}%
\providecommand \href@noop [0]{\@secondoftwo}%
\providecommand \href [0]{\begingroup \@sanitize@url \@href}%
\providecommand \@href[1]{\@@startlink{#1}\@@href}%
\providecommand \@@href[1]{\endgroup#1\@@endlink}%
\providecommand \@sanitize@url [0]{\catcode `\\12\catcode `\$12\catcode
  `\&12\catcode `\#12\catcode `\^12\catcode `\_12\catcode `\%12\relax}%
\providecommand \@@startlink[1]{}%
\providecommand \@@endlink[0]{}%
\providecommand \url  [0]{\begingroup\@sanitize@url \@url }%
\providecommand \@url [1]{\endgroup\@href {#1}{\urlprefix }}%
\providecommand \urlprefix  [0]{URL }%
\providecommand \Eprint [0]{\href }%
\providecommand \doibase [0]{http://dx.doi.org/}%
\providecommand \selectlanguage [0]{\@gobble}%
\providecommand \bibinfo  [0]{\@secondoftwo}%
\providecommand \bibfield  [0]{\@secondoftwo}%
\providecommand \translation [1]{[#1]}%
\providecommand \BibitemOpen [0]{}%
\providecommand \bibitemStop [0]{}%
\providecommand \bibitemNoStop [0]{.\EOS\space}%
\providecommand \EOS [0]{\spacefactor3000\relax}%
\providecommand \BibitemShut  [1]{\csname bibitem#1\endcsname}%
\let\auto@bib@innerbib\@empty
\bibitem [{\citenamefont {Dennis}\ \emph {et~al.}(2002)\citenamefont {Dennis},
  \citenamefont {Kitaev}, \citenamefont {Landahl},\ and\ \citenamefont
  {Preskill}}]{Dennis2002}%
  \BibitemOpen
  \bibfield  {author} {\bibinfo {author} {\bibfnamefont {E.}~\bibnamefont
  {Dennis}}, \bibinfo {author} {\bibfnamefont {A.}~\bibnamefont {Kitaev}},
  \bibinfo {author} {\bibfnamefont {A.}~\bibnamefont {Landahl}}, \ and\
  \bibinfo {author} {\bibfnamefont {J.}~\bibnamefont {Preskill}},\
  }{Topological quantum memory},\ \href {\doibase 10.1063/1.1499754} {\bibfield
   {journal} {\bibinfo  {journal} {Journal of Mathematical Physics}\ }\textbf
  {\bibinfo {volume} {43}},\ \bibinfo {pages} {4452}},\ \Eprint
  {http://arxiv.org/abs/quant-ph/0110143} {arXiv:quant-ph/0110143}  (\bibinfo
  {year} {2002})\BibitemShut {NoStop}%
\bibitem [{\citenamefont {Kitaev}(2003)}]{MR1951039}%
  \BibitemOpen
  \bibfield  {author} {\bibinfo {author} {\bibfnamefont {A.~Y.}\ \bibnamefont
  {Kitaev}},\ }Fault-tolerant quantum computation by anyons,\ \href {\doibase
  10.1016/S0003-4916(02)00018-0} {\bibfield  {journal} {\bibinfo  {journal}
  {Annals of Physics}\ }\textbf {\bibinfo {volume} {303}},\ \bibinfo {pages}
  {2}},\ \Eprint {http://arxiv.org/abs/quant-ph/9707021}
  {arXiv:quant-ph/9707021}  (\bibinfo {year} {2003})\BibitemShut {NoStop}%
\bibitem [{\citenamefont {Brown}\ \emph {et~al.}(2016)\citenamefont {Brown},
  \citenamefont {Loss}, \citenamefont {Pachos}, \citenamefont {Self},\ and\
  \citenamefont {Wootton}}]{Brown2014}%
  \BibitemOpen
  \bibfield  {author} {\bibinfo {author} {\bibfnamefont {B.~J.}\ \bibnamefont
  {Brown}}, \bibinfo {author} {\bibfnamefont {D.}~\bibnamefont {Loss}},
  \bibinfo {author} {\bibfnamefont {J.~K.}\ \bibnamefont {Pachos}}, \bibinfo
  {author} {\bibfnamefont {C.~N.}\ \bibnamefont {Self}}, \ and\ \bibinfo
  {author} {\bibfnamefont {J.~R.}\ \bibnamefont {Wootton}},\ }Quantum memories
  at finite temperature,\ \href {\doibase 10.1103/RevModPhys.88.045005}
  {\bibfield  {journal} {\bibinfo  {journal} {Reviews of Modern Physics}\
  }\textbf {\bibinfo {volume} {88}},\ \bibinfo {pages} {045005}},\ \Eprint
  {http://arxiv.org/abs/1411.6643} {arXiv:1411.6643}  (\bibinfo {year}
  {2016})\BibitemShut {NoStop}%
\bibitem [{\citenamefont {Terhal}(2015)}]{Terhal2015}%
  \BibitemOpen
  \bibfield  {author} {\bibinfo {author} {\bibfnamefont {B.~M.}\ \bibnamefont
  {Terhal}},\ }Quantum error correction for quantum memories,\ \href {\doibase
  10.1103/RevModPhys.87.307} {\bibfield  {journal} {\bibinfo  {journal}
  {Reviews of Modern Physics}\ }\textbf {\bibinfo {volume} {87}},\ \bibinfo
  {pages} {307}},\ \Eprint {http://arxiv.org/abs/1302.3428} {arXiv:1302.3428}
  (\bibinfo {year} {2015})\BibitemShut {NoStop}%
\bibitem [{\citenamefont {Raussendorf}\ and\ \citenamefont
  {Harrington}(2007)}]{0610082}%
  \BibitemOpen
  \bibfield  {author} {\bibinfo {author} {\bibfnamefont {R.}~\bibnamefont
  {Raussendorf}}\ and\ \bibinfo {author} {\bibfnamefont {J.}~\bibnamefont
  {Harrington}},\ }{Fault-tolerant quantum computation with high threshold in
  two dimensions},\ \href {\doibase 10.1103/PhysRevLett.98.190504} {\bibfield
  {journal} {\bibinfo  {journal} {Physical Review Letters}\ }\textbf {\bibinfo
  {volume} {98}},\ \bibinfo {pages} {190504}},\ \Eprint
  {http://arxiv.org/abs/quant-ph/0610082} {arXiv:quant-ph/0610082}  (\bibinfo
  {year} {2007})\BibitemShut {NoStop}%
\bibitem [{\citenamefont {Bombin}\ and\ \citenamefont
  {Martin-Delgado}(2009)}]{Bombin2007a}%
  \BibitemOpen
  \bibfield  {author} {\bibinfo {author} {\bibfnamefont {H.}~\bibnamefont
  {Bombin}}\ and\ \bibinfo {author} {\bibfnamefont {M.}~\bibnamefont
  {Martin-Delgado}},\ }{Quantum Measurements and Gates by Code Deformation},\
  \href {\doibase 10.1088/1751-8113/42/9/095302} {\bibfield  {journal}
  {\bibinfo  {journal} {Journal of Physics A: Mathematical and Theoretical}\
  }\textbf {\bibinfo {volume} {42}},\ \bibinfo {pages} {095302}},\ \Eprint
  {http://arxiv.org/abs/0704.2540} {arXiv:0704.2540}  (\bibinfo {year}
  {2009})\BibitemShut {NoStop}%
\bibitem [{\citenamefont {Bombin}(2010)}]{Bombin2010}%
  \BibitemOpen
  \bibfield  {author} {\bibinfo {author} {\bibfnamefont {H.}~\bibnamefont
  {Bombin}},\ }{Topological order with a twist: Ising anyons from an Abelian
  model},\ \href {\doibase 10.1103/PhysRevLett.105.030403} {\bibfield
  {journal} {\bibinfo  {journal} {Physical Review Letters}\ }\textbf {\bibinfo
  {volume} {105}},\ \bibinfo {pages} {030403}},\ \Eprint
  {http://arxiv.org/abs/1004.1838} {arXiv:1004.1838}  (\bibinfo {year}
  {2010})\BibitemShut {NoStop}%
\bibitem [{\citenamefont {Brown}\ \emph {et~al.}(2014)\citenamefont {Brown},
  \citenamefont {Al-Shimary},\ and\ \citenamefont {Pachos}}]{Brown2013a}%
  \BibitemOpen
  \bibfield  {author} {\bibinfo {author} {\bibfnamefont {B.~J.}\ \bibnamefont
  {Brown}}, \bibinfo {author} {\bibfnamefont {A.}~\bibnamefont {Al-Shimary}}, \
  and\ \bibinfo {author} {\bibfnamefont {J.~K.}\ \bibnamefont {Pachos}},\
  }{Entropic Barriers for Two-Dimensional Quantum Memories},\ \href {\doibase
  10.1103/PhysRevLett.112.120503} {\bibfield  {journal} {\bibinfo  {journal}
  {Physical Review Letters}\ }\textbf {\bibinfo {volume} {112}},\ \bibinfo
  {pages} {120503}},\ \Eprint {http://arxiv.org/abs/1307.6222}
  {arXiv:1307.6222}  (\bibinfo {year} {2014})\BibitemShut {NoStop}%
\bibitem [{\citenamefont {Pastawski}\ and\ \citenamefont
  {Yoshida}(2015)}]{Pastawski2014}%
  \BibitemOpen
  \bibfield  {author} {\bibinfo {author} {\bibfnamefont {F.}~\bibnamefont
  {Pastawski}}\ and\ \bibinfo {author} {\bibfnamefont {B.}~\bibnamefont
  {Yoshida}},\ }{Fault-tolerant logical gates in quantum error-correcting
  codes},\ \href {\doibase 10.1103/PhysRevA.91.012305} {\bibfield  {journal}
  {\bibinfo  {journal} {Physical Review A}\ }\textbf {\bibinfo {volume} {91}},\
  \bibinfo {pages} {012305}},\ \Eprint {http://arxiv.org/abs/1408.1720}
  {arXiv:1408.1720}  (\bibinfo {year} {2015})\BibitemShut {NoStop}%
\bibitem [{\citenamefont {Yoshida}(2015)}]{Yoshida2015a}%
  \BibitemOpen
  \bibfield  {author} {\bibinfo {author} {\bibfnamefont {B.}~\bibnamefont
  {Yoshida}},\ }{Topological color code and symmetry-protected topological
  phases},\ \href {\doibase 10.1103/PhysRevB.91.245131} {\bibfield  {journal}
  {\bibinfo  {journal} {Physical Review B}\ }\textbf {\bibinfo {volume} {91}},\
  \bibinfo {pages} {245131}},\ \Eprint {http://arxiv.org/abs/1503.07208}
  {arXiv:1503.07208}  (\bibinfo {year} {2015})\BibitemShut {NoStop}%
\bibitem [{\citenamefont {Delfosse}\ \emph {et~al.}(2016)\citenamefont
  {Delfosse}, \citenamefont {Iyer},\ and\ \citenamefont {Poulin}}]{1606.07116}%
  \BibitemOpen
  \bibfield  {author} {\bibinfo {author} {\bibfnamefont {N.}~\bibnamefont
  {Delfosse}}, \bibinfo {author} {\bibfnamefont {P.}~\bibnamefont {Iyer}}, \
  and\ \bibinfo {author} {\bibfnamefont {D.}~\bibnamefont {Poulin}},\
  }{Generalized surface codes and packing of logical qubits},\ \href@noop {} {\
  }\Eprint {http://arxiv.org/abs/1606.07116} {arXiv:1606.07116}  (\bibinfo
  {year} {2016})\BibitemShut {NoStop}%
\bibitem [{\citenamefont {Brown}\ \emph {et~al.}(2017)\citenamefont {Brown},
  \citenamefont {Laubscher}, \citenamefont {Kesselring},\ and\ \citenamefont
  {Wootton}}]{Brown2016}%
  \BibitemOpen
  \bibfield  {author} {\bibinfo {author} {\bibfnamefont {B.~J.}\ \bibnamefont
  {Brown}}, \bibinfo {author} {\bibfnamefont {K.}~\bibnamefont {Laubscher}},
  \bibinfo {author} {\bibfnamefont {M.~S.}\ \bibnamefont {Kesselring}}, \ and\
  \bibinfo {author} {\bibfnamefont {J.~R.}\ \bibnamefont {Wootton}},\ }{Poking
  Holes and Cutting Corners to Achieve Clifford Gates with the Surface Code},\
  \href {\doibase 10.1103/PhysRevX.7.021029} {\bibfield  {journal} {\bibinfo
  {journal} {Physical Review X}\ }\textbf {\bibinfo {volume} {7}},\ \bibinfo
  {pages} {021029}},\ \Eprint {http://arxiv.org/abs/1609.04673}
  {arXiv:1609.04673}  (\bibinfo {year} {2017})\BibitemShut {NoStop}%
\bibitem [{\citenamefont {Cong}\ \emph {et~al.}(2016)\citenamefont {Cong},
  \citenamefont {Cheng},\ and\ \citenamefont {Wang}}]{IrisCong1}%
  \BibitemOpen
  \bibfield  {author} {\bibinfo {author} {\bibfnamefont {I.}~\bibnamefont
  {Cong}}, \bibinfo {author} {\bibfnamefont {M.}~\bibnamefont {Cheng}}, \ and\
  \bibinfo {author} {\bibfnamefont {Z.}~\bibnamefont {Wang}},\ }{Topological
  Quantum Computation with Gapped Boundaries},\ \href@noop {} {\ }\Eprint
  {http://arxiv.org/abs/1609.02037} {arXiv:1609.02037}  (\bibinfo {year}
  {2016})\BibitemShut {NoStop}%
\bibitem [{\citenamefont {Cong}\ \emph
  {et~al.}(2017{\natexlab{a}})\citenamefont {Cong}, \citenamefont {Cheng},\
  and\ \citenamefont {Wang}}]{IrisCong2}%
  \BibitemOpen
  \bibfield  {author} {\bibinfo {author} {\bibfnamefont {I.}~\bibnamefont
  {Cong}}, \bibinfo {author} {\bibfnamefont {M.}~\bibnamefont {Cheng}}, \ and\
  \bibinfo {author} {\bibfnamefont {Z.}~\bibnamefont {Wang}},\ }{Universal
  Quantum Computation with Gapped Boundaries},\ \href {\doibase
  10.1103/PhysRevLett.119.170504} {\bibfield  {journal} {\bibinfo  {journal}
  {Physical Review Letters}\ }\textbf {\bibinfo {volume} {119}},\ \bibinfo
  {pages} {170504}},\ \Eprint {http://arxiv.org/abs/1707.05490}
  {arXiv:1707.05490}  (\bibinfo {year} {2017}{\natexlab{a}})\BibitemShut
  {NoStop}%
\bibitem [{\citenamefont {Cong}\ \emph
  {et~al.}(2017{\natexlab{b}})\citenamefont {Cong}, \citenamefont {Cheng},\
  and\ \citenamefont {Wang}}]{PhysRevB.96.195129}%
  \BibitemOpen
  \bibfield  {author} {\bibinfo {author} {\bibfnamefont {I.}~\bibnamefont
  {Cong}}, \bibinfo {author} {\bibfnamefont {M.}~\bibnamefont {Cheng}}, \ and\
  \bibinfo {author} {\bibfnamefont {Z.}~\bibnamefont {Wang}},\ }Defects between
  gapped boundaries in two-dimensional topological phases of matter,\ \href
  {\doibase 10.1103/PhysRevB.96.195129} {\bibfield  {journal} {\bibinfo
  {journal} {Physical Review B}\ }\textbf {\bibinfo {volume} {96}},\ \bibinfo
  {pages} {195129}},\ \Eprint {http://arxiv.org/abs/1703.03564}
  {arXiv:1703.03564}  (\bibinfo {year} {2017}{\natexlab{b}})\BibitemShut
  {NoStop}%
\bibitem [{\citenamefont {Yoshida}(2017)}]{Yoshida2017}%
  \BibitemOpen
  \bibfield  {author} {\bibinfo {author} {\bibfnamefont {B.}~\bibnamefont
  {Yoshida}},\ }{Gapped boundaries, group cohomology and fault-tolerant logical
  gates},\ \href {\doibase 10.1016/j.aop.2016.12.014} {\bibfield  {journal}
  {\bibinfo  {journal} {Annals of Physics}\ }\textbf {\bibinfo {volume}
  {377}},\ \bibinfo {pages} {387}},\ \Eprint {http://arxiv.org/abs/1509.03626}
  {arXiv:1509.03626}  (\bibinfo {year} {2017})\BibitemShut {NoStop}%
\bibitem [{\citenamefont {Williamson}\ \emph {et~al.}(2017)\citenamefont
  {Williamson}, \citenamefont {Bultinck},\ and\ \citenamefont
  {Verstraete}}]{SETPaper}%
  \BibitemOpen
  \bibfield  {author} {\bibinfo {author} {\bibfnamefont {D.~J.}\ \bibnamefont
  {Williamson}}, \bibinfo {author} {\bibfnamefont {N.}~\bibnamefont
  {Bultinck}}, \ and\ \bibinfo {author} {\bibfnamefont {F.}~\bibnamefont
  {Verstraete}},\ }{Symmetry-enriched topological order in tensor networks:
  Defects, gauging and anyon condensation},\ \href@noop {} {\ }\Eprint
  {http://arxiv.org/abs/1711.07982} {arXiv:1711.07982}  (\bibinfo {year}
  {2017})\BibitemShut {NoStop}%
\bibitem [{\citenamefont {Kesselring}\ \emph {et~al.}(2018)\citenamefont
  {Kesselring}, \citenamefont {Pastawski}, \citenamefont {Eisert},\ and\
  \citenamefont {Brown}}]{Brown2018}%
  \BibitemOpen
  \bibfield  {author} {\bibinfo {author} {\bibfnamefont {M.~S.}\ \bibnamefont
  {Kesselring}}, \bibinfo {author} {\bibfnamefont {F.}~\bibnamefont
  {Pastawski}}, \bibinfo {author} {\bibfnamefont {J.}~\bibnamefont {Eisert}}, \
  and\ \bibinfo {author} {\bibfnamefont {B.~J.}\ \bibnamefont {Brown}},\ }The
  boundaries and twist defects of the color code and their applications to
  topological quantum computation,\ \href {\doibase 10.22331/q-2018-10-19-101}
  {\bibfield  {journal} {\bibinfo  {journal} {{Quantum}}\ }\textbf {\bibinfo
  {volume} {2}},\ \bibinfo {pages} {101}},\ \Eprint
  {http://arxiv.org/abs/1806.02820} {arXiv:1806.02820}  (\bibinfo {year}
  {2018})\BibitemShut {NoStop}%
\bibitem [{\citenamefont {Kitaev}\ and\ \citenamefont
  {Kong}(2012)}]{MR2942952}%
  \BibitemOpen
  \bibfield  {author} {\bibinfo {author} {\bibfnamefont {A.}~\bibnamefont
  {Kitaev}}\ and\ \bibinfo {author} {\bibfnamefont {L.}~\bibnamefont {Kong}},\
  }Models for gapped boundaries and domain walls,\ \href {\doibase
  10.1007/s00220-012-1500-5} {\bibfield  {journal} {\bibinfo  {journal}
  {Communications in Mathematical Physics}\ }\textbf {\bibinfo {volume}
  {313}},\ \bibinfo {pages} {351}},\ \Eprint {http://arxiv.org/abs/1104.5047}
  {arXiv:1104.5047}  (\bibinfo {year} {2012})\BibitemShut {NoStop}%
\bibitem [{\citenamefont {Fuchs}\ \emph {et~al.}(2002)\citenamefont {Fuchs},
  \citenamefont {Runkel},\ and\ \citenamefont {Schweigert}}]{FUCHS2002353}%
  \BibitemOpen
  \bibfield  {author} {\bibinfo {author} {\bibfnamefont {J.}~\bibnamefont
  {Fuchs}}, \bibinfo {author} {\bibfnamefont {I.}~\bibnamefont {Runkel}}, \
  and\ \bibinfo {author} {\bibfnamefont {C.}~\bibnamefont {Schweigert}},\ }{TFT
  construction of RCFT correlators I: partition functions},\ \href {\doibase
  10.1016/S0550-3213(02)00744-7} {\bibfield  {journal} {\bibinfo  {journal}
  {Nuclear Physics B}\ }\textbf {\bibinfo {volume} {646}},\ \bibinfo {pages}
  {353 }},\ \Eprint {http://arxiv.org/abs/hep-th/0204148}
  {arXiv:hep-th/0204148}  (\bibinfo {year} {2002})\BibitemShut {NoStop}%
\bibitem [{\citenamefont {Fuchs}\ \emph {et~al.}(2015)\citenamefont {Fuchs},
  \citenamefont {Priel}, \citenamefont {Schweigert},\ and\ \citenamefont
  {Valentino}}]{MR3370609}%
  \BibitemOpen
  \bibfield  {author} {\bibinfo {author} {\bibfnamefont {J.}~\bibnamefont
  {Fuchs}}, \bibinfo {author} {\bibfnamefont {J.}~\bibnamefont {Priel}},
  \bibinfo {author} {\bibfnamefont {C.}~\bibnamefont {Schweigert}}, \ and\
  \bibinfo {author} {\bibfnamefont {A.}~\bibnamefont {Valentino}},\ }On the
  {B}rauer groups of symmetries of abelian {D}ijkgraaf-{W}itten theories,\
  \href {\doibase 10.1007/s00220-015-2420-y} {\bibfield  {journal} {\bibinfo
  {journal} {Communications in Mathematical Physics}\ }\textbf {\bibinfo
  {volume} {339}},\ \bibinfo {pages} {385}},\ \Eprint
  {http://arxiv.org/abs/1404.6646} {arXiv:1404.6646}  (\bibinfo {year}
  {2015})\BibitemShut {NoStop}%
\bibitem [{\citenamefont {Kong}(2014)}]{Kong2013}%
  \BibitemOpen
  \bibfield  {author} {\bibinfo {author} {\bibfnamefont {L.}~\bibnamefont
  {Kong}},\ }{Anyon condensation and tensor categories},\ \href {\doibase
  10.1016/j.nuclphysb.2014.07.003} {\bibfield  {journal} {\bibinfo  {journal}
  {Nuclear Physics B}\ }\textbf {\bibinfo {volume} {886}},\ \bibinfo {pages}
  {436}},\ \Eprint {http://arxiv.org/abs/1307.8244} {arXiv:1307.8244}
  (\bibinfo {year} {2014})\BibitemShut {NoStop}%
\bibitem [{\citenamefont {Fuchs}\ \emph {et~al.}(2013)\citenamefont {Fuchs},
  \citenamefont {Schweigert},\ and\ \citenamefont {Valentino}}]{MR3063919}%
  \BibitemOpen
  \bibfield  {author} {\bibinfo {author} {\bibfnamefont {J.}~\bibnamefont
  {Fuchs}}, \bibinfo {author} {\bibfnamefont {C.}~\bibnamefont {Schweigert}}, \
  and\ \bibinfo {author} {\bibfnamefont {A.}~\bibnamefont {Valentino}},\
  }{Bicategories for boundary conditions and for surface defects in 3-D TFT},\
  \href {\doibase 10.1007/s00220-013-1723-0} {\bibfield  {journal} {\bibinfo
  {journal} {Communications in Mathematical Physics}\ }\textbf {\bibinfo
  {volume} {321}},\ \bibinfo {pages} {543}} (\bibinfo {year}
  {2013})\BibitemShut {NoStop}%
\bibitem [{\citenamefont {Barkeshli}\ \emph {et~al.}(2013)\citenamefont
  {Barkeshli}, \citenamefont {Jian},\ and\ \citenamefont {Qi}}]{Barkeshli2013}%
  \BibitemOpen
  \bibfield  {author} {\bibinfo {author} {\bibfnamefont {M.}~\bibnamefont
  {Barkeshli}}, \bibinfo {author} {\bibfnamefont {C.-M.}\ \bibnamefont {Jian}},
  \ and\ \bibinfo {author} {\bibfnamefont {X.-L.}\ \bibnamefont {Qi}},\
  }{Theory of defects in Abelian topological states},\ \href {\doibase
  10.1103/PhysRevB.88.235103} {\bibfield  {journal} {\bibinfo  {journal}
  {Physical Review B}\ }\textbf {\bibinfo {volume} {88}},\ \bibinfo {pages}
  {235103}},\ \Eprint {http://arxiv.org/abs/1305.7203} {arXiv:1305.7203}
  (\bibinfo {year} {2013})\BibitemShut {NoStop}%
\bibitem [{\citenamefont {Barkeshli}\ \emph {et~al.}(2014)\citenamefont
  {Barkeshli}, \citenamefont {Bonderson}, \citenamefont {Cheng},\ and\
  \citenamefont {Wang}}]{Barkeshli2014}%
  \BibitemOpen
  \bibfield  {author} {\bibinfo {author} {\bibfnamefont {M.}~\bibnamefont
  {Barkeshli}}, \bibinfo {author} {\bibfnamefont {P.}~\bibnamefont
  {Bonderson}}, \bibinfo {author} {\bibfnamefont {M.}~\bibnamefont {Cheng}}, \
  and\ \bibinfo {author} {\bibfnamefont {Z.}~\bibnamefont {Wang}},\ }{Symmetry,
  Defects, and Gauging of Topological Phases},\ \href
  {http://arxiv.org/abs/1410.4540} {\ }\Eprint {http://arxiv.org/abs/1410.4540}
  {arXiv:1410.4540}  (\bibinfo {year} {2014})\BibitemShut {NoStop}%
\bibitem [{\citenamefont {Bridgeman}\ and\ \citenamefont
  {Williamson}(2017)}]{PhysRevB.96.125104}%
  \BibitemOpen
  \bibfield  {author} {\bibinfo {author} {\bibfnamefont {J.~C.}\ \bibnamefont
  {Bridgeman}}\ and\ \bibinfo {author} {\bibfnamefont {D.~J.}\ \bibnamefont
  {Williamson}},\ }Anomalies and entanglement renormalization,\ \href {\doibase
  10.1103/PhysRevB.96.125104} {\bibfield  {journal} {\bibinfo  {journal}
  {Physical Review B}\ }\textbf {\bibinfo {volume} {96}},\ \bibinfo {pages}
  {125104}},\ \Eprint {http://arxiv.org/abs/1703.07782} {arXiv:1703.07782}
  (\bibinfo {year} {2017})\BibitemShut {NoStop}%
\bibitem [{\citenamefont {Bridgeman}\ \emph {et~al.}(2017)\citenamefont
  {Bridgeman}, \citenamefont {Doherty},\ and\ \citenamefont
  {Bartlett}}]{Bridgeman2017}%
  \BibitemOpen
  \bibfield  {author} {\bibinfo {author} {\bibfnamefont {J.~C.}\ \bibnamefont
  {Bridgeman}}, \bibinfo {author} {\bibfnamefont {A.~C.}\ \bibnamefont
  {Doherty}}, \ and\ \bibinfo {author} {\bibfnamefont {S.~D.}\ \bibnamefont
  {Bartlett}},\ }Tensor networks with a twist: Anyon-permuting domain walls and
  defects in projected entangled pair states,\ \href {\doibase
  10.1103/PhysRevB.96.245122} {\bibfield  {journal} {\bibinfo  {journal}
  {Physical Review B}\ }\textbf {\bibinfo {volume} {96}},\ \bibinfo {pages}
  {245122}},\ \Eprint {http://arxiv.org/abs/1708.08930} {arXiv:1708.08930}
  (\bibinfo {year} {2017})\BibitemShut {NoStop}%
\bibitem [{\citenamefont {Cui}\ \emph {et~al.}(2019)\citenamefont {Cui},
  \citenamefont {Zini},\ and\ \citenamefont {Wang}}]{1809.00245}%
  \BibitemOpen
  \bibfield  {author} {\bibinfo {author} {\bibfnamefont {S.~X.}\ \bibnamefont
  {Cui}}, \bibinfo {author} {\bibfnamefont {M.~S.}\ \bibnamefont {Zini}}, \
  and\ \bibinfo {author} {\bibfnamefont {Z.}~\bibnamefont {Wang}},\ }On
  generalized symmetries and structure of modular categories,\ \href {\doibase
  10.1007/s11425-018-9455-5} {\bibfield  {journal} {\bibinfo  {journal}
  {Science China Mathematics}\ }10.1007/s11425-018-9455-5},\ \Eprint
  {http://arxiv.org/abs/1809.00245} {arXiv:1809.00245}  (\bibinfo {year}
  {2019})\BibitemShut {NoStop}%
\bibitem [{\citenamefont {Barter}\ \emph {et~al.}(2018)\citenamefont {Barter},
  \citenamefont {Bridgeman},\ and\ \citenamefont {Jones}}]{1806.01279}%
  \BibitemOpen
  \bibfield  {author} {\bibinfo {author} {\bibfnamefont {D.}~\bibnamefont
  {Barter}}, \bibinfo {author} {\bibfnamefont {J.~C.}\ \bibnamefont
  {Bridgeman}}, \ and\ \bibinfo {author} {\bibfnamefont {C.}~\bibnamefont
  {Jones}},\ }{Domain walls in topological phases and the Brauer-Picard ring
  for $\operatorname{Vec}(\mathbb{Z}/p\mathbb{Z})$},\ \href {\doibase
  10.1007/s00220-019-03338-2} {\bibfield  {journal} {\bibinfo  {journal}
  {Communications in Mathematical Physics}\ }10.1007/s00220-019-03338-2},\
  \bibinfo {note} {in press},\ \Eprint {http://arxiv.org/abs/1806.01279}
  {arXiv:1806.01279}  (\bibinfo {year} {2018})\BibitemShut {NoStop}%
\bibitem [{\citenamefont {Chow}\ \emph {et~al.}(2014)\citenamefont {Chow},
  \citenamefont {Gambetta}, \citenamefont {Magesan}, \citenamefont {Abraham},
  \citenamefont {Cross}, \citenamefont {Johnson}, \citenamefont {Masluk},
  \citenamefont {Ryan}, \citenamefont {Smolin}, \citenamefont {Srinivasan}
  \emph {et~al.}}]{chow2014implementing}%
  \BibitemOpen
  \bibfield  {author} {\bibinfo {author} {\bibfnamefont {J.~M.}\ \bibnamefont
  {Chow}}, \bibinfo {author} {\bibfnamefont {J.~M.}\ \bibnamefont {Gambetta}},
  \bibinfo {author} {\bibfnamefont {E.}~\bibnamefont {Magesan}}, \bibinfo
  {author} {\bibfnamefont {D.~W.}\ \bibnamefont {Abraham}}, \bibinfo {author}
  {\bibfnamefont {A.~W.}\ \bibnamefont {Cross}}, \bibinfo {author}
  {\bibfnamefont {B.}~\bibnamefont {Johnson}}, \bibinfo {author} {\bibfnamefont
  {N.~A.}\ \bibnamefont {Masluk}}, \bibinfo {author} {\bibfnamefont {C.~A.}\
  \bibnamefont {Ryan}}, \bibinfo {author} {\bibfnamefont {J.~A.}\ \bibnamefont
  {Smolin}}, \bibinfo {author} {\bibfnamefont {S.~J.}\ \bibnamefont
  {Srinivasan}},  \emph {et~al.},\ }Implementing a strand of a scalable
  fault-tolerant quantum computing fabric,\ \href {\doibase 10.1038/ncomms5015}
  {\bibfield  {journal} {\bibinfo  {journal} {Nature Communications}\ }\textbf
  {\bibinfo {volume} {5}},\ \bibinfo {pages} {4015}},\ \Eprint
  {http://arxiv.org/abs/1311.6330} {arXiv:1311.6330}  (\bibinfo {year}
  {2014})\BibitemShut {NoStop}%
\bibitem [{\citenamefont {Gambetta}\ \emph {et~al.}(2017)\citenamefont
  {Gambetta}, \citenamefont {Chow},\ and\ \citenamefont {Steffen}}]{Gambetta1}%
  \BibitemOpen
  \bibfield  {author} {\bibinfo {author} {\bibfnamefont {J.~M.}\ \bibnamefont
  {Gambetta}}, \bibinfo {author} {\bibfnamefont {J.~M.}\ \bibnamefont {Chow}},
  \ and\ \bibinfo {author} {\bibfnamefont {M.}~\bibnamefont {Steffen}},\
  }{Building logical qubits in a superconducting quantum computing system},\
  \href {\doibase 10.1038/s41534-016-0004-0} {\bibfield  {journal} {\bibinfo
  {journal} {npj Quantum Information}\ }\textbf {\bibinfo {volume} {3}},\
  \bibinfo {pages} {2}},\ \Eprint {http://arxiv.org/abs/1510.04375}
  {arXiv:1510.04375}  (\bibinfo {year} {2017})\BibitemShut {NoStop}%
\bibitem [{\citenamefont {Levin}\ and\ \citenamefont {Wen}(2005)}]{Levin2005}%
  \BibitemOpen
  \bibfield  {author} {\bibinfo {author} {\bibfnamefont {M.}~\bibnamefont
  {Levin}}\ and\ \bibinfo {author} {\bibfnamefont {X.-G.}\ \bibnamefont
  {Wen}},\ }{String-net condensation: A physical mechanism for topological
  phases},\ \href {\doibase 10.1103/PhysRevB.71.045110} {\bibfield  {journal}
  {\bibinfo  {journal} {Physical Review B}\ }\textbf {\bibinfo {volume} {71}},\
  \bibinfo {pages} {045110}},\ \Eprint {http://arxiv.org/abs/cond-mat/0404617}
  {arXiv:cond-mat/0404617}  (\bibinfo {year} {2005})\BibitemShut {NoStop}%
\bibitem [{\citenamefont {Turaev}\ and\ \citenamefont
  {Virelizier}(2017)}]{MR3674995}%
  \BibitemOpen
  \bibfield  {author} {\bibinfo {author} {\bibfnamefont {V.}~\bibnamefont
  {Turaev}}\ and\ \bibinfo {author} {\bibfnamefont {A.}~\bibnamefont
  {Virelizier}},\ }\href {\doibase 10.1007/978-3-319-49834-8} {\emph {\bibinfo
  {title} {Monoidal categories and topological field theory}}},\ \bibinfo
  {series} {Progress in Mathematics}, Vol.\ \bibinfo {volume} {322}\ (\bibinfo
  {publisher} {Birkh\"auser/Springer, Cham},\ \bibinfo {year} {2017})\ pp.\
  \bibinfo {pages} {xii+523}\BibitemShut {NoStop}%
\bibitem [{\citenamefont {Barrett}\ and\ \citenamefont
  {Westbury}(1996)}]{MR1357878}%
  \BibitemOpen
  \bibfield  {author} {\bibinfo {author} {\bibfnamefont {J.~W.}\ \bibnamefont
  {Barrett}}\ and\ \bibinfo {author} {\bibfnamefont {B.~W.}\ \bibnamefont
  {Westbury}},\ }Invariants of piecewise-linear {$3$}-manifolds,\ \href
  {\doibase 10.1090/S0002-9947-96-01660-1} {\bibfield  {journal} {\bibinfo
  {journal} {Transactions of the American Mathematical Society}\ }\textbf
  {\bibinfo {volume} {348}},\ \bibinfo {pages} {3997}},\ \Eprint
  {http://arxiv.org/abs/hep-th/9311155} {arXiv:hep-th/9311155}  (\bibinfo
  {year} {1996})\BibitemShut {NoStop}%
\bibitem [{\citenamefont {Turaev}\ and\ \citenamefont
  {Viro}(1992)}]{MR1191386}%
  \BibitemOpen
  \bibfield  {author} {\bibinfo {author} {\bibfnamefont {V.~G.}\ \bibnamefont
  {Turaev}}\ and\ \bibinfo {author} {\bibfnamefont {O.~Y.}\ \bibnamefont
  {Viro}},\ }State sum invariants of {$3$}-manifolds and quantum
  {$6j$}-symbols,\ \href {\doibase 10.1016/0040-9383(92)90015-A} {\bibfield
  {journal} {\bibinfo  {journal} {Topology}\ }\textbf {\bibinfo {volume}
  {31}},\ \bibinfo {pages} {865}} (\bibinfo {year} {1992})\BibitemShut
  {NoStop}%
\bibitem [{\citenamefont {Morrison}\ and\ \citenamefont
  {Walker}(2012)}]{MR2978449}%
  \BibitemOpen
  \bibfield  {author} {\bibinfo {author} {\bibfnamefont {S.}~\bibnamefont
  {Morrison}}\ and\ \bibinfo {author} {\bibfnamefont {K.}~\bibnamefont
  {Walker}},\ }Blob homology,\ \href {\doibase 10.2140/gt.2012.16.1481}
  {\bibfield  {journal} {\bibinfo  {journal} {Geometry \& Topology}\ }\textbf
  {\bibinfo {volume} {16}},\ \bibinfo {pages} {1481}} (\bibinfo {year}
  {2012})\BibitemShut {NoStop}%
\bibitem [{\citenamefont {Carqueville}(2016)}]{1607.05747}%
  \BibitemOpen
  \bibfield  {author} {\bibinfo {author} {\bibfnamefont {N.}~\bibnamefont
  {Carqueville}},\ }\href@noop {} {{Lecture notes on 2-dimensional defect
  TQFT}} (\bibinfo {year} {2016}),\ \Eprint {http://arxiv.org/abs/1607.05747}
  {arXiv:1607.05747} \BibitemShut {NoStop}%
\bibitem [{Note1()}]{Note1}%
  \BibitemOpen
  \bibinfo {note} {D.~Barter, J.~C.~Bridgeman, C.~Jones, \protect \emph {in
  preparation}}\BibitemShut {NoStop}%
\bibitem [{\citenamefont {{de Wild Propitius}}(1995)}]{propitius}%
  \BibitemOpen
  \bibfield  {author} {\bibinfo {author} {\bibfnamefont {M.}~\bibnamefont {{de
  Wild Propitius}}},\ }\emph {\bibinfo {title} {{Topological interactions in
  broken gauge theories}}},\ \href@noop {} {Ph.D. thesis},\ \bibinfo  {school}
  {University of Amsterdam},\ \Eprint {http://arxiv.org/abs/hep-th/9511195}
  {arXiv:hep-th/9511195}  (\bibinfo {year} {1995})\BibitemShut {NoStop}%
\bibitem [{\citenamefont {de~Wild~Propitius}(1997)}]{DEWILDPROPITIUS1997297}%
  \BibitemOpen
  \bibfield  {author} {\bibinfo {author} {\bibfnamefont {M.}~\bibnamefont
  {de~Wild~Propitius}},\ }{(Spontaneously broken) Abelian Chern-Simons
  theories},\ \href {\doibase 10.1016/S0550-3213(97)00005-9} {\bibfield
  {journal} {\bibinfo  {journal} {Nuclear Physics B}\ }\textbf {\bibinfo
  {volume} {489}},\ \bibinfo {pages} {297 }},\ \Eprint
  {http://arxiv.org/abs/hep-th/9606029} {arXiv:hep-th/9606029}  (\bibinfo
  {year} {1997})\BibitemShut {NoStop}%
\bibitem [{\citenamefont {Etingof}\ \emph {et~al.}(2010)\citenamefont
  {Etingof}, \citenamefont {Nikshych},\ and\ \citenamefont
  {Ostrik}}]{MR2677836}%
  \BibitemOpen
  \bibfield  {author} {\bibinfo {author} {\bibfnamefont {P.}~\bibnamefont
  {Etingof}}, \bibinfo {author} {\bibfnamefont {D.}~\bibnamefont {Nikshych}}, \
  and\ \bibinfo {author} {\bibfnamefont {V.}~\bibnamefont {Ostrik}},\ }Fusion
  categories and homotopy theory,\ \href {https://doi.org/10.4171/QT/6}
  {\bibfield  {journal} {\bibinfo  {journal} {Quantum Topology}\ }\textbf
  {\bibinfo {volume} {1}},\ \bibinfo {pages} {209}},\ \bibinfo {note} {with an
  appendix by Ehud Meir},\ \Eprint {http://arxiv.org/abs/0909.3140}
  {arXiv:0909.3140}  (\bibinfo {year} {2010})\BibitemShut {NoStop}%
\bibitem [{\citenamefont {Drinfeld}\ \emph {et~al.}(2010)\citenamefont
  {Drinfeld}, \citenamefont {Gelaki}, \citenamefont {Nikshych},\ and\
  \citenamefont {Ostrik}}]{MR2609644}%
  \BibitemOpen
  \bibfield  {author} {\bibinfo {author} {\bibfnamefont {V.}~\bibnamefont
  {Drinfeld}}, \bibinfo {author} {\bibfnamefont {S.}~\bibnamefont {Gelaki}},
  \bibinfo {author} {\bibfnamefont {D.}~\bibnamefont {Nikshych}}, \ and\
  \bibinfo {author} {\bibfnamefont {V.}~\bibnamefont {Ostrik}},\ }On braided
  fusion categories. {I},\ \href {\doibase 10.1007/s00029-010-0017-z}
  {\bibfield  {journal} {\bibinfo  {journal} {Selecta Mathematica. New Series}\
  }\textbf {\bibinfo {volume} {16}},\ \bibinfo {pages} {1}},\ \Eprint
  {http://arxiv.org/abs/0906.0620} {arXiv:0906.0620}  (\bibinfo {year}
  {2010})\BibitemShut {NoStop}%
\bibitem [{\citenamefont {Edie-Michell}\ \emph {et~al.}(2018)\citenamefont
  {Edie-Michell}, \citenamefont {Jones},\ and\ \citenamefont
  {Plavnik}}]{1804.01657}%
  \BibitemOpen
  \bibfield  {author} {\bibinfo {author} {\bibfnamefont {C.}~\bibnamefont
  {Edie-Michell}}, \bibinfo {author} {\bibfnamefont {C.}~\bibnamefont {Jones}},
  \ and\ \bibinfo {author} {\bibfnamefont {J.}~\bibnamefont {Plavnik}},\
  }{Fusion Rules for $\operatorname{Vec}(\mathbb{Z}/2\mathbb{Z})$ Permutation
  Gauging},\ \href@noop {} {\ }\Eprint {http://arxiv.org/abs/1804.01657}
  {arXiv:1804.01657}  (\bibinfo {year} {2018})\BibitemShut {NoStop}%
\end{thebibliography}%

\onecolumngrid
\appendix

\clearpage
\section{Defect idempotents} \label{sec:idempotents_table}
\vspace*{20mm}
	\begin{table}[h]
\begin{minipage}[t][.8\textheight][c]{\textwidth}
	\rotatebox[]{90}{
\resizebox{.98\textheight}{!}{
			\begin{tabular}{!{\vrule width 1pt}>{\columncolor[gray]{.9}[\tabcolsep]}c!{\vrule width 1pt}c!{\color[gray]{.8}\vrule}c!{\color[gray]{.8}\vrule}c!{\color[gray]{.8}\vrule}c!{\color[gray]{.8}\vrule}c!{\color[gray]{.8}\vrule}c!{\vrule width 1pt}}
				\toprule[1pt]
				\rowcolor[gray]{.9}[\tabcolsep]&$T$&$L$&$R$&$F_0$&$X_l$&$F_r$\\
				\toprule[1pt]
				$T$&$\defect{T}{T}{a}{b}{}=
				\begin{array}{c}
				\includeTikz{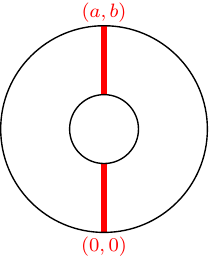}
				{
					\begin{tikzpicture}[scale=.7,every node/.style={scale=.7}]
					\annparamss{$(0,0)$}{$(a,b)$}{}{};
					\annss{}{};
					\end{tikzpicture}
				}
				\end{array}$&$\defect{T}{L}{a}{}{}=
				\begin{array}{c}
				\includeTikz{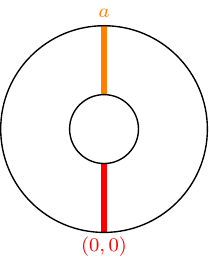}
				{
					\begin{tikzpicture}[scale=.7,every node/.style={scale=.7}]
					\annparamst{$(0,0)$}{$a$}{}{};
					\annst{}{};
					\end{tikzpicture}
				}
				\end{array}$&$\defect{T}{R}{a}{}{}=
				\begin{array}{c}
				\includeTikz{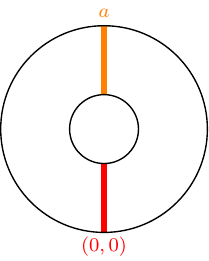}
				{
					\begin{tikzpicture}[scale=.7,every node/.style={scale=.7}]
					\annparamst{$(0,0)$}{$a$}{}{};
					\annst{}{};
					\end{tikzpicture}
				}
				\end{array}$&$\defect{T}{F_0}{}{}{}=
				\begin{array}{c}
				\includeTikz{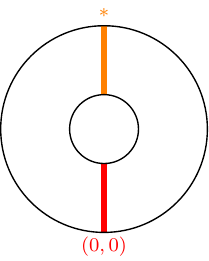}
				{
					\begin{tikzpicture}[scale=.7,every node/.style={scale=.7}]
					\annparamst{$(0,0)$}{$*$}{}{};
					\annst{}{};
					\end{tikzpicture}
				}
				\end{array}$&$\defect{T}{X_l}{a}{}{}=
				\begin{array}{c}
				\includeTikz{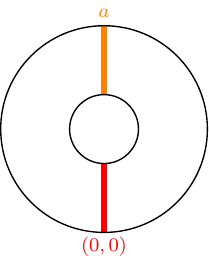}
				{
					\begin{tikzpicture}[scale=.7,every node/.style={scale=.7}]
					\annparamst{$(0,0)$}{$a$}{}{};
					\annst{}{};
					\end{tikzpicture}
				}
				\end{array}$&$\defect{T}{F_r}{}{}{}=
				\begin{array}{c}
				\includeTikz{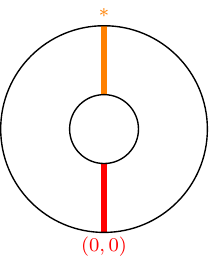}
				{
					\begin{tikzpicture}[scale=.7,every node/.style={scale=.7}]
					\annparamst{$(0,0)$}{$*$}{}{};
					\annst{}{};
					\end{tikzpicture}
				}
				\end{array}$\\
\greycline{2-7}
				$L$&
				$\defect{L}{T}{a}{}{}=
				\begin{array}{c}
				\includeTikz{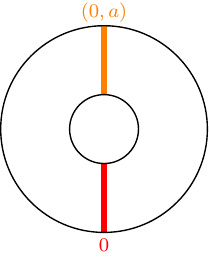}
				{
					\begin{tikzpicture}[scale=.7,every node/.style={scale=.7}]
					\annparamst{$0$}{$(0,a)$}{}{};
					\annst{}{};
					\end{tikzpicture}
				}
				\end{array}$
				&$\defect{L}{L}{a}{x}{}=\frac{1}{p}\sum_{g}\omega^{gx}
				\begin{array}{c}
				\includeTikz{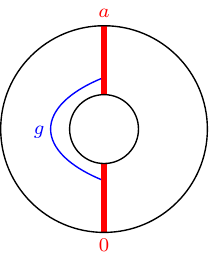}
				{
					\begin{tikzpicture}[scale=.7,every node/.style={scale=.7}]
					\annparamss{$0$}{$a$}{}{};
					\annss{$g$}{};
					\end{tikzpicture}
				}
				\end{array}$&$\defect{L}{R}{}{}{}=
				\begin{array}{c}
				\includeTikz{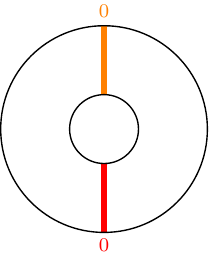}
				{
					\begin{tikzpicture}[scale=.7,every node/.style={scale=.7}]
					\annparamst{$0$}{$0$}{}{};
					\annst{}{};
					\end{tikzpicture}
				}
				\end{array}$&$\defect{L}{F_0}{x}{}{}=\frac{1}{p}\sum_g \omega^{gx}
				\begin{array}{c}
				\includeTikz{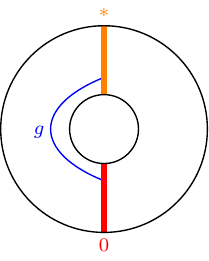}
				{
					\begin{tikzpicture}[scale=.7,every node/.style={scale=.7}]
					\annparamst{$0$}{$*$}{}{};
					\annst{$g$}{};
					\end{tikzpicture}
				}
				\end{array}$&$\defect{L}{X_l}{}{}{}=
				\begin{array}{c}
				\includeTikz{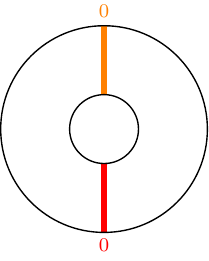}
				{
					\begin{tikzpicture}[scale=.7,every node/.style={scale=.7}]
					\annparamst{$0$}{$0$}{}{};
					\annst{}{};
					\end{tikzpicture}
				}
				\end{array}$&$\defect{L}{F_r}{x}{}{}=\frac{1}{p}\sum_g \omega^{gx}
				\begin{array}{c}
				\includeTikz{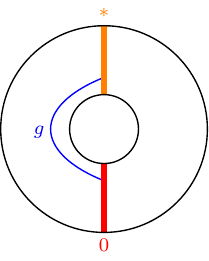}
				{
					\begin{tikzpicture}[scale=.7,every node/.style={scale=.7}]
					\annparamst{$0$}{$*$}{}{};
					\annst{$g$}{};
					\end{tikzpicture}
				}
				\end{array}$\\
\greyhline
				$R$&$\defect{R}{T}{a}{}{}=
				\begin{array}{c}
				\includeTikz{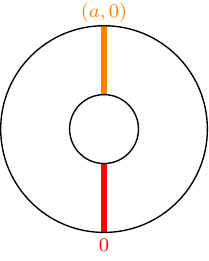}
				{
					\begin{tikzpicture}[scale=.7,every node/.style={scale=.7}]
					\annparamst{$0$}{$(a,0)$}{}{};
					\annst{}{};
					\end{tikzpicture}
				}
				\end{array}$&$\defect{R}{L}{}{}{}=
				\begin{array}{c}
				\includeTikz{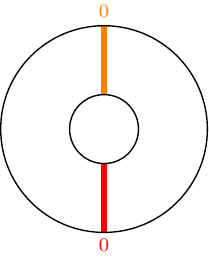}
				{
					\begin{tikzpicture}[scale=.7,every node/.style={scale=.7}]
					\annparamst{$0$}{$0$}{}{};
					\annst{}{};
					\end{tikzpicture}
				}
				\end{array}$&$\defect{R}{R}{a}{x}{}=\frac{1}{p}\sum_g\omega^{gx}
				\begin{array}{c}
				\includeTikz{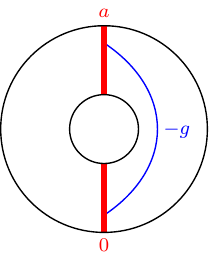}
				{
					\begin{tikzpicture}[scale=.7,every node/.style={scale=.7}]
					\annparamss{$0$}{$a$}{}{};
					\annss{}{$-g$};
					\end{tikzpicture}
				}
				\end{array}$&$\defect{R}{F_0}{x}{}{}=\frac{1}{p}\sum_g\omega^{gx}
				\begin{array}{c}
				\includeTikz{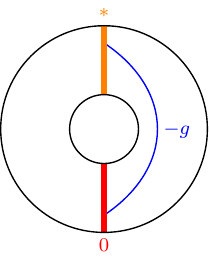}
				{
					\begin{tikzpicture}[scale=.7,every node/.style={scale=.7}]
					\annparamst{$0$}{$*$}{}{};
					\annst{}{$-g$};
					\end{tikzpicture}
				}
				\end{array}$&$\defect{R}{X_l}{}{}{}=
				\begin{array}{c}
				\includeTikz{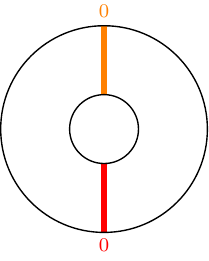}
				{
					\begin{tikzpicture}[scale=.7,every node/.style={scale=.7}]
					\annparamst{$0$}{$0$}{}{};
					\annst{}{};
					\end{tikzpicture}
				}
				\end{array}$&$\defect{R}{F_r}{x}{}{}=\frac{1}{p}\sum_g\omega^{gx}
				\begin{array}{c}
				\includeTikz{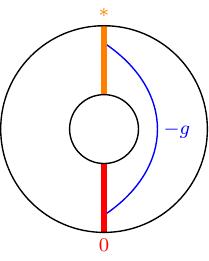}
				{
					\begin{tikzpicture}[scale=.7,every node/.style={scale=.7}]
					\annparamst{$0$}{$*$}{}{};
					\annst{}{$-g$};
					\end{tikzpicture}
				}
				\end{array}$\\
\greycline{2-7}
				$F_0$&$\defect{F_0}{T}{}{}{}=
				\begin{array}{c}
				\includeTikz{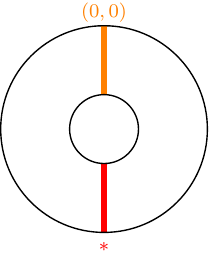}
				{
					\begin{tikzpicture}[scale=.7,every node/.style={scale=.7}]
					\annparamst{$*$}{$(0,0)$}{}{};
					\annst{}{};
					\end{tikzpicture}
				}
				\end{array}$&$\defect{F_0}{L}{x}{}{}=\frac{1}{p}\sum_g\omega^{gx}
				\begin{array}{c}
				\includeTikz{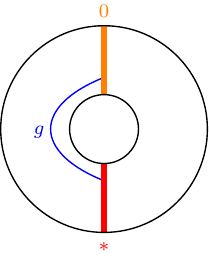}
				{
					\begin{tikzpicture}[scale=.7,every node/.style={scale=.7}]
					\annparamst{$*$}{$0$}{}{};
					\annst{$g$}{};
					\end{tikzpicture}
				}
				\end{array}$&$\defect{F_0}{R}{x}{}{}=\frac{1}{p}\sum_g\omega^{gx}
				\begin{array}{c}
				\includeTikz{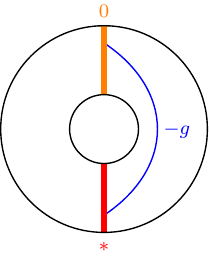}
				{
					\begin{tikzpicture}[scale=.7,every node/.style={scale=.7}]
					\annparamst{$*$}{$0$}{}{};
					\annst{}{$-g$};
					\end{tikzpicture}
				}
				\end{array}$&$\defect{F_0}{F_0}{x}{y}{}=\frac{1}{p^2}\sum_{g,h}\omega^{gx+hy}
				\begin{array}{c}
				\includeTikz{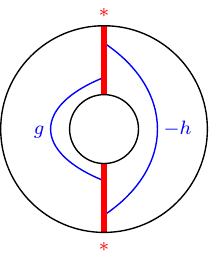}
				{
					\begin{tikzpicture}[scale=.7,every node/.style={scale=.7}]
					\annparamss{$*$}{$*$}{}{};
					\annss{$g$}{$-h$};
					\end{tikzpicture}
				}
				\end{array}$&$\defect{F_0}{X_l}{x}{}{}=\frac{1}{p}\sum_{g}\omega^{gx}
				\begin{array}{c}
				\includeTikz{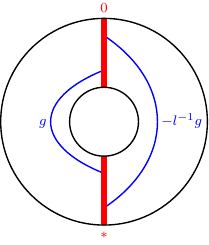}
				{
					\begin{tikzpicture}[scale=.7,every node/.style={scale=.5}]
					\annparamss{$*$}{$0$}{}{};
					\annss{$g$}{$-l^{-1}g$};
					\end{tikzpicture}
				}
				\end{array}$&$\defect{F_0}{F_r}{}{}{}=\frac{1}{p}\sum_{g}
				\begin{array}{c}
				\includeTikz{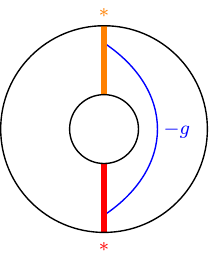}
				{
					\begin{tikzpicture}[scale=.7,every node/.style={scale=.7}]
					\annparamst{$*$}{$*$}{}{};
					\annst{}{$-g$};
					\end{tikzpicture}
				}
				\end{array}$\\
\greycline{2-7}
				$X_k$&$\defect{X_k}{T}{a}{}{}=
				\begin{array}{c}
				\includeTikz{XkT_idempotent}
				{
					\begin{tikzpicture}[scale=.7,every node/.style={scale=.7}]
					\annparamst{$0$}{$(a,0)$}{}{};
					\annst{}{};
					\end{tikzpicture}
				}
				\end{array}$&$\defect{X_k}{L}{}{}{}=
				\begin{array}{c}
				\includeTikz{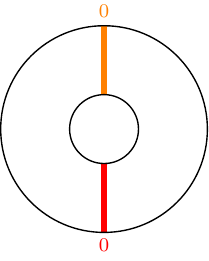}
				{
					\begin{tikzpicture}[scale=.7,every node/.style={scale=.7}]
					\annparamst{$0$}{$0$}{}{};
					\annst{}{};
					\end{tikzpicture}
				}
				\end{array}$&$\defect{X_k}{R}{}{}{}=
				\begin{array}{c}
				\includeTikz{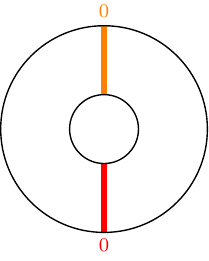}
				{
					\begin{tikzpicture}[scale=.7,every node/.style={scale=.7}]
					\annparamst{$0$}{$0$}{}{};
					\annst{}{};
					\end{tikzpicture}
				}
				\end{array}$&$\defect{X_k}{F_0}{x}{}{}=\frac{1}{p}\sum_g\omega^{gx}
				\begin{array}{c}
				\includeTikz{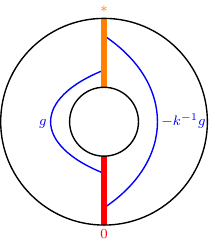}
				{
					\begin{tikzpicture}[scale=.7,every node/.style={scale=.5}]
					\annparamst{$0$}{$*$}{}{};
					\annst{$g$}{$-k^{-1}g$};
					\end{tikzpicture}
				}
				\end{array}$&
				\begin{tabular}{c}
					$\defect{X_k}{X_k}{a}{x}{}=\frac{1}{p}\sum_g\omega^{gx}
					\begin{array}{c}
					\includeTikz{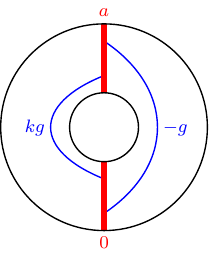}
					{
						\begin{tikzpicture}[scale=.7,every node/.style={scale=.65}]
						\annparamss{$0$}{$a$}{}{};
						\annss{$kg$}{$-g$};
						\end{tikzpicture}
					}
					\end{array}$
					\\
					\greyhline
					$\defect{X_k}{X_l}{}{}{}=
					\begin{array}{c}
					\includeTikz{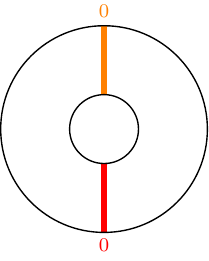}
					{
						\begin{tikzpicture}[scale=.7,every node/.style={scale=.7}]
						\annparamst{$0$}{$0$}{}{};
						\annst{}{};
						\end{tikzpicture}
					}
					\end{array}$
				\end{tabular}
				&	
				$\defect{X_k}{F_r}{x}{}{}=\frac{1}{p}\sum_g\Theta_{x,kr}(g)
				\begin{array}{c}
				\includeTikz{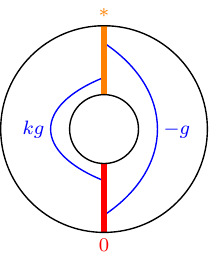}
				{
					\begin{tikzpicture}[scale=.7,every node/.style={scale=.7}]
					\annparamst{$0$}{$*$}{}{};
					\annst{$kg$}{$-g$};
					\end{tikzpicture}
				}
				\end{array}$
				\\
\greycline{2-7}
				$F_q$&$\defect{F_q}{T}{}{}{}=
				\begin{array}{c}
				\includeTikz{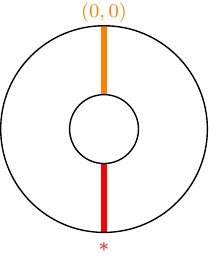}
				{
					\begin{tikzpicture}[scale=.7,every node/.style={scale=.7}]
					\annparamst{$*$}{$(0,0)$}{}{};
					\annst{}{};
					\end{tikzpicture}
				}
				\end{array}$&$\defect{F_q}{L}{x}{}{}=\frac{1}{p}\sum_g\omega^{gx}
				\begin{array}{c}
				\includeTikz{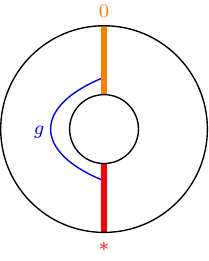}
				{
					\begin{tikzpicture}[scale=.7,every node/.style={scale=.7}]
					\annparamst{$*$}{$0$}{}{};
					\annst{$g$}{};
					\end{tikzpicture}
				}
				\end{array}$&$\defect{F_q}{R}{x}{}{}=\frac{1}{p}\sum_g\omega^{gx}
				\begin{array}{c}
				\includeTikz{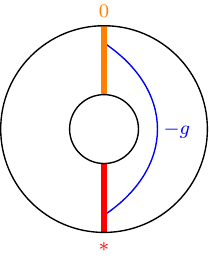}
				{
					\begin{tikzpicture}[scale=.7,every node/.style={scale=.7}]
					\annparamst{$*$}{$0$}{}{};
					\annst{}{$-g$};
					\end{tikzpicture}
				}
				\end{array}$&$\defect{F_q}{F_0}{}{}{}=\frac{1}{p}\sum_{g}
				\begin{array}{c}
				\includeTikz{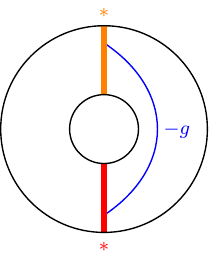}
				{
					\begin{tikzpicture}[scale=.7,every node/.style={scale=.7}]
					\annparamst{$*$}{$*$}{}{};
					\annst{}{$-g$};
					\end{tikzpicture}
				}
				\end{array}$&
				$\defect{F_q}{X_l}{x}{}{}=\frac{1}{p}\sum_g\Theta_{x,-ql}(g)
				\begin{array}{c}
				\includeTikz{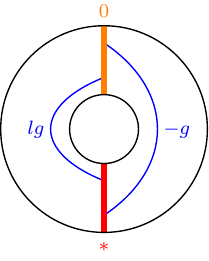}
				{
					\begin{tikzpicture}[scale=.7,every node/.style={scale=.7}]
					\annparamst{$*$}{$0$}{}{};
					\annst{$lg$}{$-g$};
					\end{tikzpicture}
				}
				\end{array}$
				&
				\begin{tabular}{c}
					$\defect{F_q}{F_q}{x}{y}{}=\frac{1}{p^2}\sum_{g,h}\omega^{gx+hy}
					\begin{array}{c}
					\includeTikz{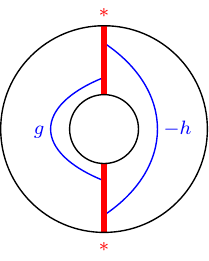}
					{
						\begin{tikzpicture}[scale=.7,every node/.style={scale=.7}]
						\annparamss{$*$}{$*$}{}{};
						\annss{$g$}{$-h$};
						\end{tikzpicture}
					}
					\end{array}$\\
					\greyhline
					$\defect{F_q}{F_r}{}{}{}=\frac{1}{p}\sum_{g}
					\begin{array}{c}
					\includeTikz{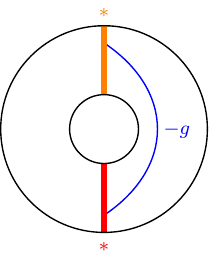}
					{
						\begin{tikzpicture}[scale=.7,every node/.style={scale=.7}]
						\annparamst{$*$}{$*$}{}{};
						\annst{}{$-g$};
						\end{tikzpicture}
					}
					\end{array}$
				\end{tabular}
				\\
				\toprule[1pt]
			\end{tabular}
		}
	}
	\caption{Indecomposable idempotents for 2-string annuli of all domain walls, corresponding to defects. For $p=2$, $\Theta_{x,a}(g)=(-1)^{gx} i^{ag}$, whilst for odd $p$ $\Theta_{x,a}(g)=\omega^{gx+ag^2 2^{-1}}$, where $2^{-1}$ is the modular inverse of 2.}\label{tab:idempotents}
\end{minipage}
	\end{table}


\clearpage

\section{Inflations} \label{sec:inflations}

When performing the horizontal fusion algorithm, the target idempotent needs to be inflated from a 2-string annulus onto a 4-string annulus. The tables in this section contain the data needed for this process. The procedure used to compute this data is explained in \onlinecite{1806.01279}. For completeness, Table~\ref{tab:zptable} contains the domain wall fusion data from \onlinecite{1806.01279}.

\begin{table}[h]
	\begin{tabular}{!{\vrule width 1pt}>{\columncolor[gray]{.9}[\tabcolsep]}c!{\vrule width 1pt}c c c c!{\vrule width 1pt}c c!{\vrule width 1pt}}
		\toprule[1pt]
		\rowcolor[gray]{.9}[\tabcolsep]$\otimes_{\vvec{\ZZ{p}}}$&$T$&$L$&$R$&$F_0$&$X_l$&$F_r$\\
		\toprule[1pt]
		$T$&$p\cdot T$&$T$&$p\cdot R$&$R$&$T$&$R$\\
		$L$&$p\cdot L$&$L$&$p\cdot F_0$&$F_0$&$L$&$F_0$\\
		$R$&$T$&$p\cdot T$&$R$&$p\cdot R$&$R$&$T$\\
		$F_0$&$L$&$p\cdot L$&$F_0$&$p\cdot F_0$&$F_0$&$L$\\
		\toprule[1pt]
		$X_k$&$T$&$L$&$R$&$F_0$&$X_{kl}$&$F_{k^{-1}r}$\\
		$F_q$&$L$&$T$&$F_0$&$R$&$F_{ql}$&$X_{q^{-1}r}$\\
		\toprule[1pt]
	\end{tabular}
	\caption{Multiplication table for $\protect\bpr{\protect\vvec{\ZZ{p}}}$, reproduced from \onlinecite{1806.01279}.}\label{tab:zptable}
\end{table}

\addtocounter{table}{1}

\newcounter{curtable}
\setcounter{curtable}{0}
\addtocounter{curtable}{\value{table}}

\setcounter{table}{0}
\renewcommand{\thetable}{\Roman{curtable}(\alph{table})}

\begin{table}[h!]
			\begin{tabular}{!{\vrule width 1pt}>{\columncolor[gray]{.9}[\tabcolsep]}c!{\vrule width 1pt}c!{\color[gray]{.8}\vrule}c!{\color[gray]{.8}\vrule}c!{\vrule width 1pt}c!{\color[gray]{.8}\vrule}c!{\vrule width 1pt}}
				\toprule[1pt]
&$T\otimes_{\vvec{\ZZ{p}}}T$&\multirow{20}{*}{
	$\begin{array}{c}
		\includeTikz{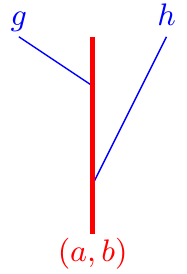}{\begin{tikzpicture}
			\inflationalhs{$g$}{$(a,b)$}{$h$};
			\end{tikzpicture}}
	\end{array}$}&$\begin{array}{c}\includeTikz{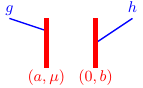}{\begin{tikzpicture}[xscale=.5,yscale=.25,every node/.style={scale=.5}]\inflationarhss{$g$}{$(a,\mu)$}{}{$(0,b)$}{$h$};\end{tikzpicture}}\end{array}$
&\multirow{20}{*}{
	$\begin{array}{c}
	\includeTikz{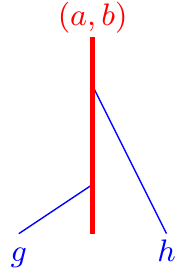}{\begin{tikzpicture}
		\inflationblhs{$g$}{$(a,b)$}{$h$};
		\end{tikzpicture}}
	\end{array}$}&$\begin{array}{c}\includeTikz{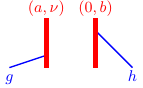}{\begin{tikzpicture}[xscale=.5,yscale=.25,every node/.style={scale=.5}]\inflationbrhss{$g$}{$(a,\nu)$}{}{$(0,b)$}{$h$};\end{tikzpicture}}\end{array}$\\
				\greycline{2-2}\greycline{4-4}\greycline{6-6}
				&$T\otimes_{\vvec{\ZZ{p}}}L$&&$\begin{array}{c}\includeTikz{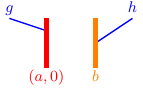}{\begin{tikzpicture}[xscale=.5,yscale=.25,every node/.style={scale=.5}]\inflationarhs{$g$}{$(a,0)$}{}{$b$}{$h$};\end{tikzpicture}}\end{array}$&&$\begin{array}{c}\includeTikz{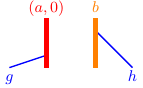}{\begin{tikzpicture}[xscale=.5,yscale=.25,every node/.style={scale=.5}]\inflationbrhs{$g$}{$(a,0)$}{}{$b$}{$h$};\end{tikzpicture}}\end{array}$\\
				\greycline{2-2}\greycline{4-4}\greycline{6-6}
				&$T\otimes_{\vvec{\ZZ{p}}}X_l$&&$\begin{array}{c}\includeTikz{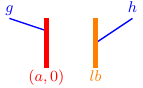}{\begin{tikzpicture}[xscale=.5,yscale=.25,every node/.style={scale=.5}]\inflationarhs{$g$}{$(a,0)$}{}{$lb$}{$h$};\end{tikzpicture}}\end{array}$&&$\begin{array}{c}\includeTikz{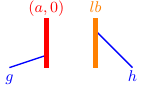}{\begin{tikzpicture}[xscale=.5,yscale=.25,every node/.style={scale=.5}]\inflationbrhs{$g$}{$(a,0)$}{}{$lb$}{$h$};\end{tikzpicture}}\end{array}$\\
				\greycline{2-2}\greycline{4-4}\greycline{6-6}
				&$R\otimes_{\vvec{\ZZ{p}}}T$&&$\begin{array}{c}\includeTikz{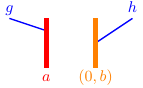}{\begin{tikzpicture}[xscale=.5,yscale=.25,every node/.style={scale=.5}]\inflationarhs{$g$}{$a$}{}{$(0,b)$}{$h$};\end{tikzpicture}}\end{array}$&&$\begin{array}{c}\includeTikz{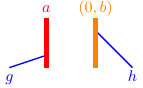}{\begin{tikzpicture}[xscale=.5,yscale=.25,every node/.style={scale=.5}]\inflationbrhs{$g$}{$a$}{}{$(0,b)$}{$h$};\end{tikzpicture}}\end{array}$\\
				\greycline{2-2}\greycline{4-4}\greycline{6-6}
				&$R\otimes_{\vvec{\ZZ{p}}}L$&&$\frac{1}{p}\sum_k\omega^{\mu k}\begin{array}{c}\includeTikz{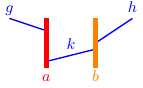}{\begin{tikzpicture}[xscale=.5,yscale=.25,every node/.style={scale=.5}]\inflationarhs{$g$}{$a$}{$k$}{$b$}{$h$};\end{tikzpicture}}\end{array}$&&$\frac{1}{p}\sum_k\omega^{\nu k}\begin{array}{c}\includeTikz{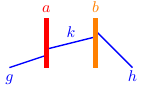}{\begin{tikzpicture}[xscale=.5,yscale=.25,every node/.style={scale=.5}]\inflationbrhs{$g$}{$a$}{$k$}{$b$}{$h$};\end{tikzpicture}}\end{array}$\\
				\greycline{2-2}\greycline{4-4}\greycline{6-6}
				&$R\otimes_{\vvec{\ZZ{p}}}F_r$&&$\frac{1}{p}\sum_k\omega^{kbr}\begin{array}{c}\includeTikz{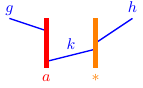}{\begin{tikzpicture}[xscale=.5,yscale=.25,every node/.style={scale=.5}]\inflationarhs{$g$}{$a$}{$k$}{$*$}{$h$};\end{tikzpicture}}\end{array}$&&$\frac{1}{p}\sum_k\omega^{kr(b-h)}\begin{array}{c}\includeTikz{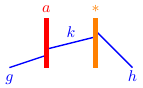}{\begin{tikzpicture}[xscale=.5,yscale=.25,every node/.style={scale=.5}]\inflationbrhs{$g$}{$a$}{$k$}{$*$}{$h$};\end{tikzpicture}}\end{array}$\\
				\greycline{2-2}\greycline{4-4}\greycline{6-6}
				&$X_k\otimes_{\vvec{\ZZ{p}}}T$&&$\begin{array}{c}\includeTikz{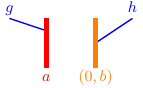}{\begin{tikzpicture}[xscale=.5,yscale=.25,every node/.style={scale=.5}]\inflationarhs{$g$}{$a$}{}{$(0,b)$}{$h$};\end{tikzpicture}}\end{array}$&&$\begin{array}{c}\includeTikz{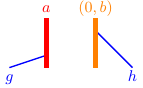}{\begin{tikzpicture}[xscale=.5,yscale=.25,every node/.style={scale=.5}]\inflationbrhs{$g$}{$a$}{}{$(0,b)$}{$h$};\end{tikzpicture}}\end{array}$\\
				\greycline{2-2}\greycline{4-4}\greycline{6-6}
				\multirow{-20}{*}{$T$}&$F_q\otimes_{\vvec{\ZZ{p}}}L$&&$\frac{1}{p}\sum_k\omega^{qak}\begin{array}{c}\includeTikz{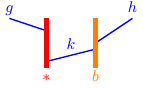}{\begin{tikzpicture}[xscale=.5,yscale=.25,every node/.style={scale=.5}]\inflationarhs{$g$}{$*$}{$k$}{$b$}{$h$};\end{tikzpicture}}\end{array}$&&$\frac{1}{p}\sum_k\omega^{kqa}\begin{array}{c}\includeTikz{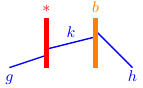}{\begin{tikzpicture}[xscale=.5,yscale=.25,every node/.style={scale=.5}]\inflationbrhs{$g$}{$*$}{$k$}{$b$}{$h$};\end{tikzpicture}}\end{array}$\\
\toprule[1pt]
&$L\otimes_{\vvec{\ZZ{p}}}T$&\multirow{20}{*}{
	$\begin{array}{c}
	\includeTikz{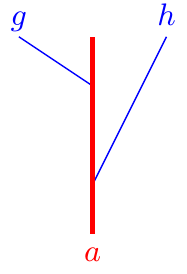}{\begin{tikzpicture}
		\inflationalhs{$g$}{$a$}{$h$};
		\end{tikzpicture}}
	\end{array}$}&$\begin{array}{c}\includeTikz{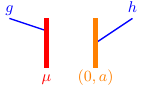}{\begin{tikzpicture}[xscale=.5,yscale=.25,every node/.style={scale=.5}]\inflationarhs{$g$}{$\mu$}{}{$(0,a)$}{$h$};\end{tikzpicture}}\end{array}$&\multirow{20}{*}{
	$\begin{array}{c}
	\includeTikz{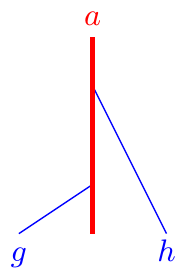}{\begin{tikzpicture}
		\inflationblhs{$g$}{$a$}{$h$};
		\end{tikzpicture}}
	\end{array}$}&$\begin{array}{c}\includeTikz{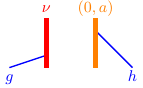}{\begin{tikzpicture}[xscale=.5,yscale=.25,every node/.style={scale=.5}]\inflationbrhs{$g$}{$\nu$}{}{$(0,a)$}{$h$};\end{tikzpicture}}\end{array}$\\
\greycline{2-2}\greycline{4-4}\greycline{6-6}
				&$L\otimes_{\vvec{\ZZ{p}}}L$&&$\begin{array}{c}\includeTikz{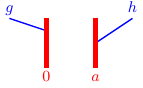}{\begin{tikzpicture}[xscale=.5,yscale=.25,every node/.style={scale=.5}]\inflationarhss{$g$}{$0$}{}{$a$}{$h$};\end{tikzpicture}}\end{array}$&&$\begin{array}{c}\includeTikz{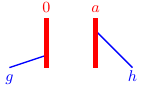}{\begin{tikzpicture}[xscale=.5,yscale=.25,every node/.style={scale=.5}]\inflationbrhss{$g$}{$0$}{}{$a$}{$h$};\end{tikzpicture}}\end{array}$\\
				\greycline{2-2}\greycline{4-4}\greycline{6-6}
				&$L\otimes_{\vvec{\ZZ{p}}}X_l$&&$\begin{array}{c}\includeTikz{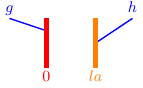}{\begin{tikzpicture}[xscale=.5,yscale=.25,every node/.style={scale=.5}]\inflationarhs{$g$}{$0$}{}{$la$}{$h$};\end{tikzpicture}}\end{array}$&&$\begin{array}{c}\includeTikz{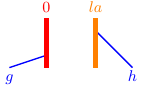}{\begin{tikzpicture}[xscale=.5,yscale=.25,every node/.style={scale=.5}]\inflationbrhs{$g$}{$0$}{}{$la$}{$h$};\end{tikzpicture}}\end{array}$\\
				\greycline{2-2}\greycline{4-4}\greycline{6-6}
				&$F_0\otimes_{\vvec{\ZZ{p}}}T$&&$\begin{array}{c}\includeTikz{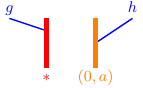}{\begin{tikzpicture}[xscale=.5,yscale=.25,every node/.style={scale=.5}]\inflationarhs{$g$}{$*$}{}{$(0,a)$}{$h$};\end{tikzpicture}}\end{array}$&&$\begin{array}{c}\includeTikz{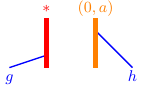}{\begin{tikzpicture}[xscale=.5,yscale=.25,every node/.style={scale=.5}]\inflationbrhs{$g$}{$*$}{}{$(0,a)$}{$h$};\end{tikzpicture}}\end{array}$\\
				\greycline{2-2}\greycline{4-4}\greycline{6-6}
				&$F_0\otimes_{\vvec{\ZZ{p}}}L$&&$\frac{1}{p}\sum_k\omega^{k \mu}\begin{array}{c}\includeTikz{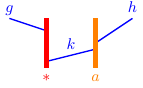}{\begin{tikzpicture}[xscale=.5,yscale=.25,every node/.style={scale=.5}]\inflationarhs{$g$}{$*$}{$k$}{$a$}{$h$};\end{tikzpicture}}\end{array}$&&$\frac{1}{p}\sum_k\omega^{k \nu}\begin{array}{c}\includeTikz{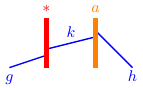}{\begin{tikzpicture}[xscale=.5,yscale=.25,every node/.style={scale=.5}]\inflationbrhs{$g$}{$*$}{$k$}{$a$}{$h$};\end{tikzpicture}}\end{array}$\\
				\greycline{2-2}\greycline{4-4}\greycline{6-6}
				&$F_0\otimes_{\vvec{\ZZ{p}}}F_r$&&$\frac{1}{p}\sum_k\omega^{kar}\begin{array}{c}\includeTikz{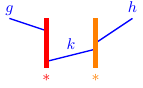}{\begin{tikzpicture}[xscale=.5,yscale=.25,every node/.style={scale=.5}]\inflationarhs{$g$}{$*$}{$k$}{$*$}{$h$};\end{tikzpicture}}\end{array}$&&$\frac{1}{p}\sum_k\omega^{k(a-h)r}\begin{array}{c}\includeTikz{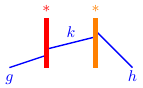}{\begin{tikzpicture}[xscale=.5,yscale=.25,every node/.style={scale=.5}]\inflationbrhs{$g$}{$*$}{$k$}{$*$}{$h$};\end{tikzpicture}}\end{array}$\\
				\greycline{2-2}\greycline{4-4}\greycline{6-6}
				&$X_k\otimes_{\vvec{\ZZ{p}}}L$&&$\begin{array}{c}\includeTikz{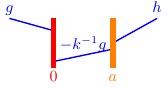}{\begin{tikzpicture}[xscale=.6,yscale=.25,every node/.style={scale=.5}]\inflationarhs{$g$}{$0$}{$-k^{-1}g$}{$a$}{$h$};\end{tikzpicture}}\end{array}$&&$\begin{array}{c}\includeTikz{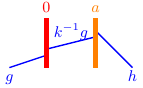}{\begin{tikzpicture}[xscale=.5,yscale=.25,every node/.style={scale=.5}]\inflationbrhs{$g$}{$0$}{$k^{-1}g$}{$a$}{$h$};\end{tikzpicture}}\end{array}$\\
				\greycline{2-2}\greycline{4-4}\greycline{6-6}
				\multirow{-20}{*}{$L$}&$F_q\otimes_{\vvec{\ZZ{p}}}T$&&$\begin{array}{c}\includeTikz{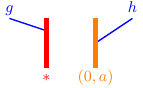}{\begin{tikzpicture}[xscale=.5,yscale=.25,every node/.style={scale=.5}]\inflationarhs{$g$}{$*$}{}{$(0,a)$}{$h$};\end{tikzpicture}}\end{array}$&&$\begin{array}{c}\includeTikz{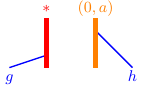}{\begin{tikzpicture}[xscale=.5,yscale=.25,every node/.style={scale=.5}]\inflationbrhs{$g$}{$*$}{}{$(0,a)$}{$h$};\end{tikzpicture}}\end{array}$\\
\toprule[1pt]
			\end{tabular}
	\caption{Inflations (part a). All $\mu,\,\nu$ occurring label components of the tensor decomposition} \label{tab:inflation_1}
	\vspace*{-10mm}
\end{table}

\begin{table}
	\begin{tabular}{!{\vrule width 1pt}>{\columncolor[gray]{.9}[\tabcolsep]}c!{\vrule width 1pt}c!{\color[gray]{.8}\vrule}c!{\color[gray]{.8}\vrule}c!{\vrule width 1pt}c!{\color[gray]{.8}\vrule}c!{\vrule width 1pt}}
		\toprule[1pt]
&$T\otimes_{\vvec{\ZZ{p}}}R$&\multirow{20}{*}{
	$\begin{array}{c}
	\includeTikz{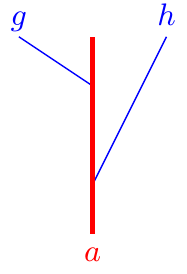}{\begin{tikzpicture}
		\inflationalhs{$g$}{$a$}{$h$};
		\end{tikzpicture}}
	\end{array}$}&$\begin{array}{c}\includeTikz{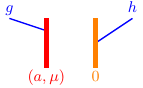}{\begin{tikzpicture}[xscale=.5,yscale=.25,every node/.style={scale=.5}]\inflationarhs{$g$}{$(a,\mu)$}{}{$0$}{$h$};\end{tikzpicture}}\end{array}$&\multirow{20}{*}{
	$\begin{array}{c}
	\includeTikz{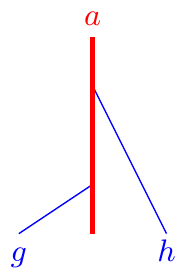}{\begin{tikzpicture}
		\inflationblhs{$g$}{$a$}{$h$};
		\end{tikzpicture}}
	\end{array}$}&$\begin{array}{c}\includeTikz{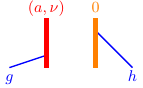}{\begin{tikzpicture}[xscale=.5,yscale=.25,every node/.style={scale=.5}]\inflationbrhs{$g$}{$(a,\nu)$}{}{$0$}{$h$};\end{tikzpicture}}\end{array}$\\
\greycline{2-2}\greycline{4-4}\greycline{6-6}
&$T\otimes_{\vvec{\ZZ{p}}}F_0$&&$\begin{array}{c}\includeTikz{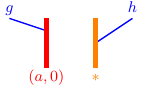}{\begin{tikzpicture}[xscale=.5,yscale=.25,every node/.style={scale=.5}]\inflationarhs{$g$}{$(a,0)$}{}{$*$}{$h$};\end{tikzpicture}}\end{array}$&&$\begin{array}{c}\includeTikz{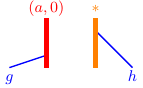}{\begin{tikzpicture}[xscale=.5,yscale=.25,every node/.style={scale=.5}]\inflationbrhs{$g$}{$(a,0)$}{}{$*$}{$h$};\end{tikzpicture}}\end{array}$\\
\greycline{2-2}\greycline{4-4}\greycline{6-6}
&$T\otimes_{\vvec{\ZZ{p}}}F_r$&&$\begin{array}{c}\includeTikz{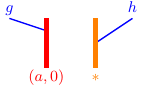}{\begin{tikzpicture}[xscale=.5,yscale=.25,every node/.style={scale=.5}]\inflationarhs{$g$}{$(a,0)$}{}{$*$}{$h$};\end{tikzpicture}}\end{array}$&&$\begin{array}{c}\includeTikz{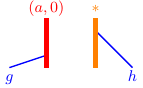}{\begin{tikzpicture}[xscale=.5,yscale=.25,every node/.style={scale=.5}]\inflationbrhs{$g$}{$(a,0)$}{}{$*$}{$h$};\end{tikzpicture}}\end{array}$\\
\greycline{2-2}\greycline{4-4}\greycline{6-6}
&$R\otimes_{\vvec{\ZZ{p}}}R$&&$\begin{array}{c}\includeTikz{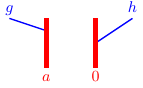}{\begin{tikzpicture}[xscale=.5,yscale=.25,every node/.style={scale=.5}]\inflationarhss{$g$}{$a$}{}{$0$}{$h$};\end{tikzpicture}}\end{array}$&&$\begin{array}{c}\includeTikz{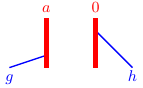}{\begin{tikzpicture}[xscale=.5,yscale=.25,every node/.style={scale=.5}]\inflationbrhss{$g$}{$a$}{}{$0$}{$h$};\end{tikzpicture}}\end{array}$\\
\greycline{2-2}\greycline{4-4}\greycline{6-6}
&$R\otimes_{\vvec{\ZZ{p}}}F_0$&&$\frac{1}{p}\sum_k\omega^{k \mu}\begin{array}{c}\includeTikz{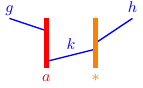}{\begin{tikzpicture}[xscale=.5,yscale=.25,every node/.style={scale=.5}]\inflationarhs{$g$}{$a$}{$k$}{$*$}{$h$};\end{tikzpicture}}\end{array}$&&$\frac{1}{p}\sum_k\omega^{k \nu}\begin{array}{c}\includeTikz{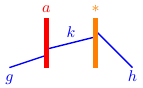}{\begin{tikzpicture}[xscale=.5,yscale=.25,every node/.style={scale=.5}]\inflationbrhs{$g$}{$a$}{$k$}{$*$}{$h$};\end{tikzpicture}}\end{array}$\\
\greycline{2-2}\greycline{4-4}\greycline{6-6}
&$R\otimes_{\vvec{\ZZ{p}}}X_l$&&$\begin{array}{c}\includeTikz{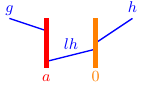}{\begin{tikzpicture}[xscale=.5,yscale=.25,every node/.style={scale=.5}]\inflationarhs{$g$}{$a$}{$lh$}{$0$}{$h$};\end{tikzpicture}}\end{array}$&&$\begin{array}{c}\includeTikz{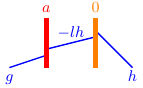}{\begin{tikzpicture}[xscale=.5,yscale=.25,every node/.style={scale=.5}]\inflationbrhs{$g$}{$a$}{$-lh$}{$0$}{$h$};\end{tikzpicture}}\end{array}$\\
\greycline{2-2}\greycline{4-4}\greycline{6-6}
&$X_k\otimes_{\vvec{\ZZ{p}}}R$&&$\begin{array}{c}\includeTikz{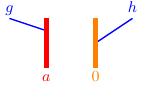}{\begin{tikzpicture}[xscale=.5,yscale=.25,every node/.style={scale=.5}]\inflationarhs{$g$}{$a$}{}{$0$}{$h$};\end{tikzpicture}}\end{array}$&&$\begin{array}{c}\includeTikz{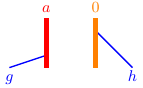}{\begin{tikzpicture}[xscale=.5,yscale=.25,every node/.style={scale=.5}]\inflationbrhs{$g$}{$a$}{}{$0$}{$h$};\end{tikzpicture}}\end{array}$\\
\greycline{2-2}\greycline{4-4}\greycline{6-6}
\multirow{-20}{*}{$R$}&$F_q\otimes_{\vvec{\ZZ{p}}}F_0$&&$\frac{1}{p}\sum_k\omega^{qak}\begin{array}{c}\includeTikz{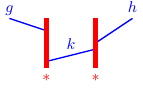}{\begin{tikzpicture}[xscale=.5,yscale=.25,every node/.style={scale=.5}]\inflationarhss{$g$}{$*$}{$k$}{$*$}{$h$};\end{tikzpicture}}\end{array}$&&$\frac{1}{p}\sum_k\omega^{qak}\begin{array}{c}\includeTikz{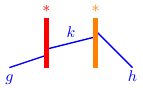}{\begin{tikzpicture}[xscale=.5,yscale=.25,every node/.style={scale=.5}]\inflationbrhs{$g$}{$*$}{$k$}{$*$}{$h$};\end{tikzpicture}}\end{array}$\\
\toprule[1pt]
		&$L\otimes_{\vvec{\ZZ{p}}}R$&\multirow{20}{*}{
			$\begin{array}{c}
			\includeTikz{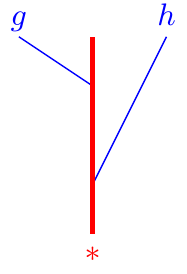}{\begin{tikzpicture}
				\inflationalhs{$g$}{$*$}{$h$};
				\end{tikzpicture}}
			\end{array}$}&$\begin{array}{c}\includeTikz{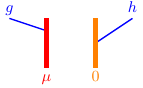}{\begin{tikzpicture}[xscale=.5,yscale=.25,every node/.style={scale=.5}]\inflationarhs{$g$}{$\mu$}{}{$0$}{$h$};\end{tikzpicture}}\end{array}$&\multirow{20}{*}{
			$\begin{array}{c}
			\includeTikz{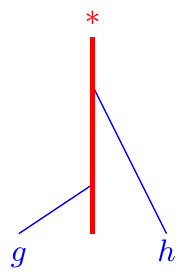}{\begin{tikzpicture}
				\inflationblhs{$g$}{$*$}{$h$};
				\end{tikzpicture}}
			\end{array}$}&$\begin{array}{c}\includeTikz{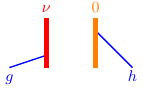}{\begin{tikzpicture}[xscale=.5,yscale=.25,every node/.style={scale=.5}]\inflationbrhs{$g$}{$\nu$}{}{$0$}{$h$};\end{tikzpicture}}\end{array}$\\
		\greycline{2-2}\greycline{4-4}\greycline{6-6}
		&$L\otimes_{\vvec{\ZZ{p}}}F_0$&&$\begin{array}{c}\includeTikz{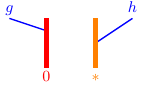}{\begin{tikzpicture}[xscale=.5,yscale=.25,every node/.style={scale=.5}]\inflationarhs{$g$}{$0$}{}{$*$}{$h$};\end{tikzpicture}}\end{array}$&&$\begin{array}{c}\includeTikz{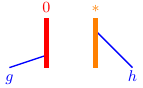}{\begin{tikzpicture}[xscale=.5,yscale=.25,every node/.style={scale=.5}]\inflationbrhs{$g$}{$0$}{}{$*$}{$h$};\end{tikzpicture}}\end{array}$\\
		\greycline{2-2}\greycline{4-4}\greycline{6-6}
		&$L\otimes_{\vvec{\ZZ{p}}}F_r$&&$\begin{array}{c}\includeTikz{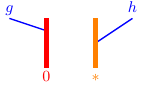}{\begin{tikzpicture}[xscale=.5,yscale=.25,every node/.style={scale=.5}]\inflationarhs{$g$}{$0$}{}{$*$}{$h$};\end{tikzpicture}}\end{array}$&&$\begin{array}{c}\includeTikz{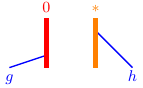}{\begin{tikzpicture}[xscale=.5,yscale=.25,every node/.style={scale=.5}]\inflationbrhs{$g$}{$0$}{}{$*$}{$h$};\end{tikzpicture}}\end{array}$\\
		\greycline{2-2}\greycline{4-4}\greycline{6-6}
		&$F_0\otimes_{\vvec{\ZZ{p}}}R$&&$\begin{array}{c}\includeTikz{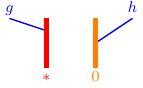}{\begin{tikzpicture}[xscale=.5,yscale=.25,every node/.style={scale=.5}]\inflationarhs{$g$}{$*$}{}{$0$}{$h$};\end{tikzpicture}}\end{array}$&&$\begin{array}{c}\includeTikz{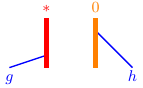}{\begin{tikzpicture}[xscale=.5,yscale=.25,every node/.style={scale=.5}]\inflationbrhs{$g$}{$*$}{}{$0$}{$h$};\end{tikzpicture}}\end{array}$\\
		\greycline{2-2}\greycline{4-4}\greycline{6-6}
		&$F_0\otimes_{\vvec{\ZZ{p}}}F_0$&&$\frac{1}{p}\sum_k\omega^{k \mu}\begin{array}{c}\includeTikz{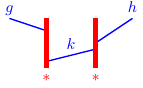}{\begin{tikzpicture}[xscale=.5,yscale=.25,every node/.style={scale=.5}]\inflationarhss{$g$}{$*$}{$k$}{$*$}{$h$};\end{tikzpicture}}\end{array}$&&$\frac{1}{p}\sum_k\omega^{k \nu}\begin{array}{c}\includeTikz{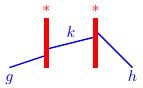}{\begin{tikzpicture}[xscale=.5,yscale=.25,every node/.style={scale=.5}]\inflationbrhss{$g$}{$*$}{$k$}{$*$}{$h$};\end{tikzpicture}}\end{array}$\\
		\greycline{2-2}\greycline{4-4}\greycline{6-6}
		&$F_0\otimes_{\vvec{\ZZ{p}}}X_l$&&$\begin{array}{c}\includeTikz{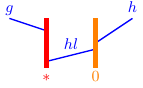}{\begin{tikzpicture}[xscale=.5,yscale=.25,every node/.style={scale=.5}]\inflationarhs{$g$}{$*$}{$hl$}{$0$}{$h$};\end{tikzpicture}}\end{array}$&&$\begin{array}{c}\includeTikz{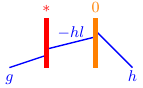}{\begin{tikzpicture}[xscale=.5,yscale=.25,every node/.style={scale=.5}]\inflationbrhs{$g$}{$*$}{$-hl$}{$0$}{$h$};\end{tikzpicture}}\end{array}$\\
		\greycline{2-2}\greycline{4-4}\greycline{6-6}
		&$X_k\otimes_{\vvec{\ZZ{p}}}F_0$&&$\begin{array}{c}\includeTikz{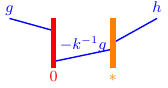}{\begin{tikzpicture}[xscale=.6,yscale=.25,every node/.style={scale=.5}]\inflationarhs{$g$}{$0$}{$-k^{-1}g$}{$*$}{$h$};\end{tikzpicture}}\end{array}$&&$\begin{array}{c}\includeTikz{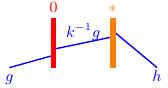}{\begin{tikzpicture}[xscale=.6,yscale=.25,every node/.style={scale=.5}]\inflationbrhs{$g$}{$0$}{$k^{-1}g$}{$*$}{$h$};\end{tikzpicture}}\end{array}$\\
		\greycline{2-2}\greycline{4-4}\greycline{6-6}
		\multirow{-20}{*}{$F_0$}&$F_q\otimes_{\vvec{\ZZ{p}}}R$&&$\begin{array}{c}\includeTikz{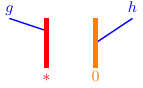}{\begin{tikzpicture}[xscale=.5,yscale=.25,every node/.style={scale=.5}]\inflationarhs{$g$}{$*$}{}{$0$}{$h$};\end{tikzpicture}}\end{array}$&&$\begin{array}{c}\includeTikz{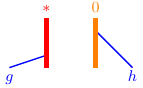}{\begin{tikzpicture}[xscale=.5,yscale=.25,every node/.style={scale=.5}]\inflationbrhs{$g$}{$*$}{}{$0$}{$h$};\end{tikzpicture}}\end{array}$\\
\toprule[1pt]
		&$X_k\otimes_{\vvec{\ZZ{p}}}X_l,\,m=kl$&\multirow{2}{*}{
			$\begin{array}{c}
			\includeTikz{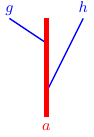}{\begin{tikzpicture}[scale=.5,every node/.style={scale=.5}]
				\inflationalhs{$g$}{$a$}{$h$};
				\end{tikzpicture}}
			\end{array}$}&$\begin{array}{c}\includeTikz{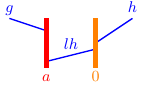}{\begin{tikzpicture}[xscale=.5,yscale=.25,every node/.style={scale=.5}]\inflationarhs{$g$}{$a$}{$lh$}{$0$}{$h$};\end{tikzpicture}}\end{array}$&\multirow{2}{*}{
			$\begin{array}{c}
			\includeTikz{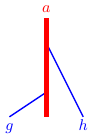}{\begin{tikzpicture}[scale=.5,every node/.style={scale=.5}]
				\inflationblhs{$g$}{$a$}{$h$};
				\end{tikzpicture}}
			\end{array}$}&$\begin{array}{c}\includeTikz{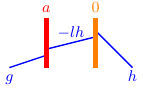}{\begin{tikzpicture}[xscale=.5,yscale=.25,every node/.style={scale=.5}]\inflationbrhs{$g$}{$a$}{$-lh$}{$0$}{$h$};\end{tikzpicture}}\end{array}$\\
		\greycline{2-2}\greycline{4-4}\greycline{6-6}
		\multirow{-4}{*}{$X_m$}&$F_q\otimes_{\vvec{\ZZ{p}}}F_r,\,m=q^{-1}r$&&$\frac{1}{p}\sum_k\omega^{qka}\begin{array}{c}\includeTikz{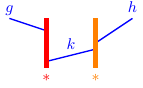}{\begin{tikzpicture}[xscale=.5,yscale=.25,every node/.style={scale=.5}]\inflationarhs{$g$}{$*$}{$k$}{$*$}{$h$};\end{tikzpicture}}\end{array}$&&$\frac{1}{p}\sum_k\omega^{(qa-rh)k}\begin{array}{c}\includeTikz{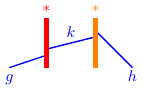}{\begin{tikzpicture}[xscale=.5,yscale=.25,every node/.style={scale=.5}]\inflationbrhs{$g$}{$*$}{$k$}{$*$}{$h$};\end{tikzpicture}}\end{array}$\\
\toprule[1pt]
		&$X_k\otimes_{\vvec{\ZZ{p}}}F_r,\,n=k^{-1}r$&\multirow{2}{*}{
			$\begin{array}{c}
			\includeTikz{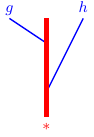}{\begin{tikzpicture}[scale=.5,every node/.style={scale=.5}]
				\inflationalhs{$g$}{$*$}{$h$};
				\end{tikzpicture}}
			\end{array}$}&$\omega^{-ngh}\begin{array}{c}\includeTikz{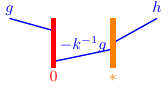}{\begin{tikzpicture}[xscale=.6,yscale=.25,every node/.style={scale=.5}]\inflationarhs{$g$}{$0$}{$-k^{-1}g$}{$*$}{$h$};\end{tikzpicture}}\end{array}$&\multirow{2}{*}{
			$\begin{array}{c}
			\includeTikz{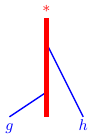}{\begin{tikzpicture}[scale=.5,every node/.style={scale=.5}]
				\inflationblhs{$g$}{$*$}{$h$};
				\end{tikzpicture}}
			\end{array}$}&$\begin{array}{c}\includeTikz{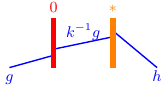}{\begin{tikzpicture}[xscale=.6,yscale=.25,every node/.style={scale=.5}]\inflationbrhs{$g$}{$0$}{$k^{-1}g$}{$*$}{$h$};\end{tikzpicture}}\end{array}$\\
		\greycline{2-2}\greycline{4-4}\greycline{6-6}
		\multirow{-4}{*}{$F_n$}&$F_q\otimes_{\vvec{\ZZ{p}}}X_l,\,n=ql$&&$\begin{array}{c}\includeTikz{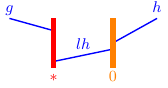}{\begin{tikzpicture}[xscale=.6,yscale=.25,every node/.style={scale=.5}]\inflationarhs{$g$}{$*$}{$lh$}{$0$}{$h$};\end{tikzpicture}}\end{array}$&&$\begin{array}{c}\includeTikz{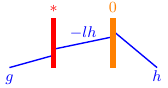}{\begin{tikzpicture}[xscale=.6,yscale=.25,every node/.style={scale=.5}]\inflationbrhs{$g$}{$*$}{$-lh$}{$0$}{$h$};\end{tikzpicture}}\end{array}$\\
\toprule[1pt]
	\end{tabular}
	\caption{Inflations (part b). All $\mu,\,\nu$ occurring label components of the tensor decomposition} \label{tab:inflation_2}
\end{table}

\setcounter{table}{\value{curtable}}
\renewcommand{\thetable}{\Roman{table}}


\clearpage

\addtocounter{table}{1}

\setcounter{curtable}{0}
\addtocounter{curtable}{\value{table}}

\setcounter{table}{0}
\renewcommand{\thetable}{\Roman{curtable}(\alph{table})}

\section{Horizontal fusion outcomes} \label{sec:horizontal_fusion_table}
\vspace*{12mm}
\begin{table}[h]
	\begin{minipage}[t][.8\textheight][c]{\textwidth}
		\renewcommand{\arraystretch}{1.5}
\begin{tabular}{!{\vrule width 1pt}>{\columncolor[gray]{.9}[\tabcolsep]}c!{\vrule width 1pt}c!{\color[gray]{.8}\vrule}c!{\color[gray]{.8}\vrule}c!{\color[gray]{.8}\vrule}c!{\color[gray]{.8}\vrule}c!{\color[gray]{.8}\vrule}c!{\vrule width 1pt}}
			\toprule[1pt]
			\rowcolor[gray]{.9}[\tabcolsep]&$ \defect{T}{T}{c}{d}{} $&$ \defect{T}{L}{c}{}{} $&$ \defect{T}{R}{c}{}{} $&$ \defect{T}{F_{0}}{}{}{} $&$ \defect{T}{X_{n}}{c}{}{} $&$ \defect{T}{F_{t}}{}{}{} $\\
			\toprule[1pt]
			$ \defect{T}{T}{a}{b}{} $&$ \delta_{\nu -\mu}^{b+c} \defect{T}{T}{a}{d}{\mu,\nu} $&$ \defect{T}{T}{a}{c}{\mu} $&$ \delta_{\nu -\mu}^{b+c} \defect{T}{R}{a}{}{\mu,\nu} $&$ \defect{T}{R}{a}{}{\mu} $&$ \defect{T}{T}{a}{n^{-1}(b+c+\mu )}{\mu} $&$ \defect{T}{R}{a}{}{\mu} $\\\greyhline
			$ \defect{T}{L}{a}{}{} $&$ \delta_{\nu -\mu}^{a+c} \defect{T}{L}{d}{}{\mu,\nu} $&$ \defect{T}{L}{c}{}{\mu} $&$ \delta_{\nu -\mu}^{a+c} \defect{T}{F_{0}}{}{}{\mu,\nu} $&$ \defect{T}{F_{0}}{}{}{\mu} $&$ \defect{T}{L}{n^{-1}(a+c+\mu )}{}{\mu} $&$ \defect{T}{F_{0}}{}{}{\mu} $\\\greyhline
			$ \defect{T}{R}{a}{}{} $&$ \defect{T}{T}{a}{d}{\mu} $&$ \defect{T}{T}{a}{c}{\mu,\nu} $&$ \defect{T}{R}{a}{}{\mu} $&$ \defect{T}{R}{a}{}{\mu,\nu} $&$ \defect{T}{R}{a}{}{\mu} $&$ \oplus_{\beta} \defect{T}{T}{a}{\beta}{\mu} $\\\greyhline
			$ \defect{T}{F_{0}}{}{}{} $&$ \defect{T}{L}{d}{}{\mu} $&$ \defect{T}{L}{c}{}{\mu,\nu} $&$ \defect{T}{F_{0}}{}{}{\mu} $&$ \defect{T}{F_{0}}{}{}{\mu,\nu} $&$ \defect{T}{F_{0}}{}{}{\mu} $&$ \oplus_{\alpha} \defect{T}{L}{\alpha}{}{\mu} $\\\greyhline
			$ \defect{T}{X_{l}}{a}{}{} $&$ \defect{T}{T}{a+l (c+\mu )}{d}{\mu} $&$ \defect{T}{L}{c}{}{\mu} $&$ \defect{T}{R}{a+l (c+\mu )}{}{\mu} $&$ \defect{T}{F_{0}}{}{}{\mu} $&$ \defect{T}{X_{l n}}{a+l (c+\mu )}{}{\mu} $&$ \defect{T}{F_{l^{-1}t}}{}{}{\mu} $\\
\end{tabular}
			\end{minipage}
			\vspace*{16mm}
			\caption{Defect fusion table (part a). $\mu$ ($\nu$) indexes degeneracy in the bottom (top) domain wall fusion.}
			\vspace*{-30mm}
			\label{tab:horizontal_table_1}
			\end{table}

\begin{table} 
	\renewcommand{\arraystretch}{1.5}
	\begin{tabular}{!{\vrule width 1pt}>{\columncolor[gray]{.9}[\tabcolsep]}c!{\vrule width 1pt}c!{\color[gray]{.8}\vrule}c!{\color[gray]{.8}\vrule}c!{\color[gray]{.8}\vrule}c!{\color[gray]{.8}\vrule}c!{\color[gray]{.8}\vrule}c!{\vrule width 1pt}}
		\toprule[1pt]
		\rowcolor[gray]{.9}[\tabcolsep]&$ \defect{L}{T}{c}{}{} $&$ \defect{L}{L}{c}{z}{} $&$ \defect{L}{R}{}{}{} $&$ \defect{L}{F_{0}}{z}{}{} $&$ \defect{L}{X_{n}}{}{}{} $&$ \defect{L}{F_{t}}{z}{}{} $\\\toprule[1pt]
		$ \defect{T}{T}{a}{b}{} $&$ \defect{T}{T}{a}{c}{\nu} $&$ \defect{T}{T}{a}{c}{} $&$ \defect{T}{R}{a}{}{\nu} $&$ \defect{T}{R}{a}{}{} $&$ \oplus_{\beta} \defect{T}{T}{a}{\beta}{} $&$ \defect{T}{R}{a}{}{} $\\\greyhline
		$ \defect{T}{L}{a}{}{} $&$ \defect{T}{L}{c}{}{\nu} $&$ \defect{T}{L}{c}{}{} $&$ \defect{T}{F_{0}}{}{}{\nu} $&$ \defect{T}{F_{0}}{}{}{} $&$ \oplus_{\alpha} \defect{T}{L}{\alpha}{}{} $&$ \defect{T}{F_{0}}{}{}{} $\\\greyhline
		$ \defect{T}{R}{a}{}{} $&$ p\cdot \defect{T}{T}{a}{c}{} $&$ \defect{T}{T}{a}{c}{\nu} $&$ p\cdot \defect{T}{R}{a}{}{} $&$ \defect{T}{R}{a}{}{\nu} $&$ p\cdot \defect{T}{R}{a}{}{} $&$ \oplus_{\beta} \defect{T}{T}{a}{\beta}{} $\\\greyhline
		$ \defect{T}{F_{0}}{}{}{} $&$ p\cdot \defect{T}{L}{c}{}{} $&$ \defect{T}{L}{c}{}{\nu} $&$ p\cdot \defect{T}{F_{0}}{}{}{} $&$ \defect{T}{F_{0}}{}{}{\nu} $&$ p\cdot \defect{T}{F_{0}}{}{}{} $&$ \oplus_{\alpha} \defect{T}{L}{\alpha}{}{} $\\\greyhline
		$ \defect{T}{X_{l}}{a}{}{} $&$ \oplus_{\alpha} \defect{T}{T}{\alpha}{c}{} $&$ \defect{T}{L}{c}{}{} $&$ \oplus_{\alpha} \defect{T}{R}{\alpha}{}{} $&$ \defect{T}{F_{0}}{}{}{} $&$ \oplus_{\alpha} \defect{T}{X_{l n}}{\alpha}{}{} $&$ \defect{T}{F_{l^{-1}t}}{}{}{} $\\\greyhline
		$ \defect{T}{F_{r}}{}{}{} $&$ p\cdot \defect{T}{L}{c}{}{} $&$ \oplus_{\alpha} \defect{T}{T}{\alpha}{c}{} $&$ p\cdot \defect{T}{F_{0}}{}{}{} $&$ \oplus_{\alpha} \defect{T}{R}{\alpha}{}{} $&$ p\cdot \defect{T}{F_{n r}}{}{}{} $&$ \oplus_{\alpha} \defect{T}{X_{r^{-1}t}}{\alpha}{}{} $\\
\end{tabular}
	\caption{Defect fusion table (part b). $\mu$ ($\nu$) indexes degeneracy in the bottom (top) domain wall fusion.}
	\label{tab:horizontal_table_2}
\end{table}

\begin{table} 
	\renewcommand{\arraystretch}{1.5}
	\begin{tabular}{!{\vrule width 1pt}>{\columncolor[gray]{.9}[\tabcolsep]}c!{\vrule width 1pt}c!{\color[gray]{.8}\vrule}c!{\color[gray]{.8}\vrule}c!{\color[gray]{.8}\vrule}c!{\color[gray]{.8}\vrule}c!{\color[gray]{.8}\vrule}c!{\vrule width 1pt}}
		\toprule[1pt]
		\rowcolor[gray]{.9}[\tabcolsep]&$ \defect{R}{T}{c}{}{} $&$ \defect{R}{L}{}{}{} $&$ \defect{R}{R}{c}{z}{} $&$ \defect{R}{F_{0}}{z}{}{} $&$ \defect{R}{X_{n}}{}{}{} $&$ \defect{R}{F_{t}}{z}{}{} $\\\toprule[1pt]
		$ \defect{T}{T}{a}{b}{} $&$ \delta_{\nu -\mu}^{b+c} \defect{R}{T}{a}{}{\mu,\nu} $&$ \defect{R}{T}{a}{}{\mu} $&$ \delta_{\nu -\mu}^{b+c} \defect{R}{R}{a}{z}{\mu,\nu} $&$ \defect{R}{R}{a}{z}{\mu} $&$ \defect{R}{T}{a}{}{\mu} $&$ \defect{R}{R}{a}{z+t (b+\mu )}{\mu} $\\\greyhline
		$ \defect{T}{L}{a}{}{} $&$ \delta_{\nu -\mu}^{a+c} \defect{R}{L}{}{}{\mu,\nu} $&$ \defect{R}{L}{}{}{\mu} $&$ \delta_{\nu -\mu}^{a+c} \defect{R}{F_{0}}{z}{}{\mu,\nu} $&$ \defect{R}{F_{0}}{z}{}{\mu} $&$ \defect{R}{L}{}{}{\mu} $&$ \defect{R}{F_{0}}{z+t (a+\mu )}{}{\mu} $\\\greyhline
		$ \defect{T}{R}{a}{}{} $&$ \defect{R}{T}{a}{}{\mu} $&$ \defect{R}{T}{a}{}{\mu,\nu} $&$ \defect{R}{R}{a}{z}{\mu} $&$ \defect{R}{R}{a}{z}{\mu,\nu} $&$ \oplus_{\beta} \defect{R}{R}{a}{\beta}{\mu} $&$ \defect{R}{T}{a}{}{\mu} $\\\greyhline
		$ \defect{T}{F_{0}}{}{}{} $&$ \defect{R}{L}{}{}{\mu} $&$ \defect{R}{L}{}{}{\mu,\nu} $&$ \defect{R}{F_{0}}{z}{}{\mu} $&$ \defect{R}{F_{0}}{z}{}{\mu,\nu} $&$ \oplus_{\alpha} \defect{R}{F_{0}}{\alpha}{}{\mu} $&$ \defect{R}{L}{}{}{\mu} $\\\greyhline
		$ \defect{T}{X_{l}}{a}{}{} $&$ \defect{R}{T}{a+l (c+\mu )}{}{\mu} $&$ \defect{R}{L}{}{}{\mu} $&$ \defect{R}{R}{a+l (c+\mu )}{z}{\mu} $&$ \defect{R}{F_{0}}{z}{}{\mu} $&$ \defect{R}{X_{l n}}{}{}{\mu} $&$ \defect{R}{F_{l^{-1}t}}{\mu  t+l^{-1}(a t)+z}{}{\mu} $\\
\end{tabular}
	\caption{Defect fusion table (part c). $\mu$ ($\nu$) indexes degeneracy in the bottom (top) domain wall fusion.}
	\label{tab:horizontal_table_3}
\end{table}

\begin{table} 
	\renewcommand{\arraystretch}{1.5}
	\begin{tabular}{!{\vrule width 1pt}>{\columncolor[gray]{.9}[\tabcolsep]}c!{\vrule width 1pt}c!{\color[gray]{.8}\vrule}c!{\color[gray]{.8}\vrule}c!{\color[gray]{.8}\vrule}c!{\color[gray]{.8}\vrule}c!{\color[gray]{.8}\vrule}c!{\vrule width 1pt}}
		\toprule[1pt]
		\rowcolor[gray]{.9}[\tabcolsep]&$ \defect{F_{0}}{T}{}{}{} $&$ \defect{F_{0}}{L}{z}{}{} $&$ \defect{F_{0}}{R}{z}{}{} $&$ \defect{F_{0}}{F_{0}}{z}{w}{} $&$ \defect{F_{0}}{X_{n}}{z}{}{} $&$ \defect{F_{0}}{F_{t}}{}{}{} $\\\toprule[1pt]
		$ \defect{T}{T}{a}{b}{} $&$ \defect{R}{T}{a}{}{\nu} $&$ \defect{R}{T}{a}{}{} $&$ \defect{R}{R}{a}{z}{\nu} $&$ \defect{R}{R}{a}{w}{} $&$ \defect{R}{T}{a}{}{} $&$ \oplus_{\beta} \defect{R}{R}{a}{\beta}{} $\\\greyhline
		$ \defect{T}{L}{a}{}{} $&$ \defect{R}{L}{}{}{\nu} $&$ \defect{R}{L}{}{}{} $&$ \defect{R}{F_{0}}{z}{}{\nu} $&$ \defect{R}{F_{0}}{w}{}{} $&$ \defect{R}{L}{}{}{} $&$ \oplus_{\alpha} \defect{R}{F_{0}}{\alpha}{}{} $\\\greyhline
		$ \defect{T}{R}{a}{}{} $&$ p\cdot \defect{R}{T}{a}{}{} $&$ \defect{R}{T}{a}{}{\nu} $&$ p\cdot \defect{R}{R}{a}{z}{} $&$ \defect{R}{R}{a}{w}{\nu} $&$ \oplus_{\beta} \defect{R}{R}{a}{\beta}{} $&$ p\cdot \defect{R}{T}{a}{}{} $\\\greyhline
		$ \defect{T}{F_{0}}{}{}{} $&$ p\cdot \defect{R}{L}{}{}{} $&$ \defect{R}{L}{}{}{\nu} $&$ p\cdot \defect{R}{F_{0}}{z}{}{} $&$ \defect{R}{F_{0}}{w}{}{\nu} $&$ \oplus_{\alpha} \defect{R}{F_{0}}{\alpha}{}{} $&$ p\cdot \defect{R}{L}{}{}{} $\\\greyhline
		$ \defect{T}{X_{l}}{a}{}{} $&$ \oplus_{\alpha} \defect{R}{T}{\alpha}{}{} $&$ \defect{R}{L}{}{}{} $&$ \oplus_{\alpha} \defect{R}{R}{\alpha}{z}{} $&$ \defect{R}{F_{0}}{w}{}{} $&$ \defect{R}{X_{l n}}{}{}{} $&$ \oplus_{\alpha} \defect{R}{F_{l^{-1}t}}{\alpha}{}{} $\\\greyhline
		$ \defect{T}{F_{r}}{}{}{} $&$ p\cdot \defect{R}{L}{}{}{} $&$ \oplus_{\alpha} \defect{R}{T}{\alpha}{}{} $&$ p\cdot \defect{R}{F_{0}}{z}{}{} $&$ \oplus_{\alpha} \defect{R}{R}{\alpha}{w}{} $&$ \oplus_{\alpha} \defect{R}{F_{n r}}{\alpha}{}{} $&$ p\cdot \defect{R}{X_{r^{-1}t}}{}{}{} $\\
\end{tabular}
	\caption{Defect fusion table (part d). $\mu$ ($\nu$) indexes degeneracy in the bottom (top) domain wall fusion.}
	\label{tab:horizontal_table_4}
\end{table}

\begin{table}  
	\resizebox{\textwidth}{!}{
		\renewcommand{\arraystretch}{1.5}
		\begin{tabular}{!{\vrule width 1pt}>{\columncolor[gray]{.9}[\tabcolsep]}c!{\vrule width 1pt}c!{\color[gray]{.8}\vrule}c!{\color[gray]{.8}\vrule}c!{\color[gray]{.8}\vrule}c!{\color[gray]{.8}\vrule}c!{\color[gray]{.8}\vrule}c!{\color[gray]{.8}\vrule}c!{\vrule width 1pt}}
			\toprule[1pt]
			\rowcolor[gray]{.9}[\tabcolsep]&$ \defect{X_{m}}{T}{c}{}{} $&$ \defect{X_{m}}{L}{}{}{} $&$ \defect{X_{m}}{R}{}{}{} $&$ \defect{X_{m}}{F_{0}}{z}{}{} $&$ \defect{X_{m}}{X_{m}}{c}{z}{} $&$ \defect{X_{m}}{X_{n}}{}{}{} $&$ \defect{X_{m}}{F_{t}}{z}{}{} $\\
			\toprule[1pt]
			$ \defect{T}{T}{a}{b}{} $&$ \defect{T}{T}{a}{m^{-1}(b+c-\nu )}{\nu} $&$ \oplus_{\beta} \defect{T}{T}{a}{\beta}{} $&$ \defect{T}{R}{a}{}{\nu} $&$ \defect{T}{R}{a}{}{} $&$ \defect{T}{T}{a}{m^{-1}(b+c)}{} $&$ \oplus_{\beta} \defect{T}{T}{a}{\beta}{} $&$ \defect{T}{R}{a}{}{} $\\\greyhline
			$ \defect{T}{L}{a}{}{} $&$ \defect{T}{L}{m^{-1}(a+c-\nu )}{}{\nu} $&$ \oplus_{\alpha} \defect{T}{L}{\alpha}{}{} $&$ \defect{T}{F_{0}}{}{}{\nu} $&$ \defect{T}{F_{0}}{}{}{} $&$ \defect{T}{L}{m^{-1}(a+c)}{}{} $&$ \oplus_{\alpha} \defect{T}{L}{\alpha}{}{} $&$ \defect{T}{F_{0}}{}{}{} $\\\greyhline
			$ \defect{T}{R}{a}{}{} $&$ \oplus_{\beta} \defect{T}{T}{a}{\beta}{} $&$ \oplus_{\beta} \defect{T}{T}{a}{\beta}{\nu} $&$ p\cdot \defect{T}{R}{a}{}{} $&$ \defect{T}{R}{a}{}{\nu} $&$ \defect{T}{R}{a}{}{} $&$ p\cdot \defect{T}{R}{a}{}{} $&$ \oplus_{\beta} \defect{T}{T}{a}{\beta}{} $\\\greyhline
			$ \defect{T}{F_{0}}{}{}{} $&$ \oplus_{\alpha} \defect{T}{L}{\alpha}{}{} $&$ \oplus_{\alpha} \defect{T}{L}{\alpha}{}{\nu} $&$ p\cdot \defect{T}{F_{0}}{}{}{} $&$ \defect{T}{F_{0}}{}{}{\nu} $&$ \defect{T}{F_{0}}{}{}{} $&$ p\cdot \defect{T}{F_{0}}{}{}{} $&$ \oplus_{\alpha} \defect{T}{L}{\alpha}{}{} $\\
\end{tabular}
	}
	\caption{Defect fusion table (part e). $\mu$ ($\nu$) indexes degeneracy in the bottom (top) domain wall fusion.}
	\label{tab:horizontal_table_5}
\end{table}

\begin{table} 
	\resizebox{\textwidth}{!}{
		\renewcommand{\arraystretch}{1.5}
		\begin{tabular}{!{\vrule width 1pt}>{\columncolor[gray]{.9}[\tabcolsep]}c!{\vrule width 1pt}c!{\color[gray]{.8}\vrule}c!{\color[gray]{.8}\vrule}c!{\color[gray]{.8}\vrule}c!{\color[gray]{.8}\vrule}c!{\color[gray]{.8}\vrule}c!{\color[gray]{.8}\vrule}c!{\vrule width 1pt}}
			\toprule[1pt]
			\rowcolor[gray]{.9}[\tabcolsep]&$ \defect{F_{s}}{T}{}{}{} $&$ \defect{F_{s}}{L}{z}{}{} $&$ \defect{F_{s}}{R}{z}{}{} $&$ \defect{F_{s}}{F_{0}}{}{}{} $&$ \defect{F_{s}}{X_{n}}{z}{}{} $&$ \defect{F_{s}}{F_{s}}{z}{w}{} $&$ \defect{F_{s}}{F_{t}}{}{}{} $\\
			\toprule[1pt]
			$ \defect{T}{T}{a}{b}{} $&$ \defect{R}{T}{a}{}{\nu} $&$ \defect{R}{T}{a}{}{} $&$ \defect{R}{R}{a}{b s-\nu  s+z}{\nu} $&$ \oplus_{\beta} \defect{R}{R}{a}{\beta}{} $&$ \defect{R}{T}{a}{}{} $&$ \defect{R}{R}{a}{b s+w}{} $&$ \oplus_{\beta} \defect{R}{R}{a}{\beta}{} $\\\greyhline
			$ \defect{T}{L}{a}{}{} $&$ \defect{R}{L}{}{}{\nu} $&$ \defect{R}{L}{}{}{} $&$ \defect{R}{F_{0}}{a s-\nu  s+z}{}{\nu} $&$ \oplus_{\alpha} \defect{R}{F_{0}}{\alpha}{}{} $&$ \defect{R}{L}{}{}{} $&$ \defect{R}{F_{0}}{a s+w}{}{} $&$ \oplus_{\alpha} \defect{R}{F_{0}}{\alpha}{}{} $\\\greyhline
			$ \defect{T}{R}{a}{}{} $&$ p\cdot \defect{R}{T}{a}{}{} $&$ \defect{R}{T}{a}{}{\nu} $&$ \oplus_{\beta} \defect{R}{R}{a}{\beta}{} $&$ \oplus_{\beta} \defect{R}{R}{a}{\beta}{\nu} $&$ \oplus_{\beta} \defect{R}{R}{a}{\beta}{} $&$ \defect{R}{T}{a}{}{} $&$ p\cdot \defect{R}{T}{a}{}{} $\\\greyhline
			$ \defect{T}{F_{0}}{}{}{} $&$ p\cdot \defect{R}{L}{}{}{} $&$ \defect{R}{L}{}{}{\nu} $&$ \oplus_{\alpha} \defect{R}{F_{0}}{\alpha}{}{} $&$ \oplus_{\alpha} \defect{R}{F_{0}}{\alpha}{}{\nu} $&$ \oplus_{\alpha} \defect{R}{F_{0}}{\alpha}{}{} $&$ \defect{R}{L}{}{}{} $&$ p\cdot \defect{R}{L}{}{}{} $\\
\end{tabular}
	}
	\caption{Defect fusion table (part f). $\mu$ ($\nu$) indexes degeneracy in the bottom (top) domain wall fusion.}
	\label{tab:horizontal_table_6}
\end{table}

\setcounter{table}{\value{curtable}}
\renewcommand{\thetable}{\Roman{table}}


\clearpage

\section{Vertical fusion outcomes} \label{sec:vertical_fusion_table}
\newlength{\tabwidth}
\setlength{\tabwidth}{.49\textheight}
\newlength{\tabheight}
\setlength{\tabheight}{.3\textwidth}
\vspace*{10mm}
\begin{table}[h] 
	\begin{minipage}[t][.85\textheight][c]{\textwidth}
	\renewcommand{\arraystretch}{2.9}
\rotatebox[]{90}{
	\resizebox*{\tabwidth}{\tabheight}{
		\begin{tabular}{!{\vrule width 1pt}>{\columncolor[gray]{.9}[\tabcolsep]}c!{\vrule width 1pt}c !{\color[gray]{.8}\vrule} c !{\color[gray]{.8}\vrule} c !{\color[gray]{.8}\vrule} c!{\vrule width 1pt}c !{\color[gray]{.8}\vrule} c!{\vrule width 1pt}}
			\toprule[1pt]
			\rowcolor[gray]{.9}[\tabcolsep]\cellcolor{red!20}$L$&$\defect{L}{T}{c}{}{}$&$\defect{L}{L}{c}{z}{}$&$\defect{L}{R}{}{}{}$&$\defect{L}{F_{0}}{z}{}{}$&$\defect{L}{X_{m}}{}{}{}$&$\defect{L}{F_{s}}{z}{}{}$\\
			\toprule[1pt]
			$\defect{T}{L}{a}{}{}$&$\oplus_{\alpha}\defect{T}{T}{\alpha}{a+c}{}$&$\defect{T}{L}{a+c}{}{}$&$\oplus_{\alpha}\defect{T}{R}{\alpha}{}{}$&$\defect{T}{F_{0}}{}{}{}$&$\oplus_{\alpha}\defect{T}{X_{m}}{\alpha}{}{}$&$\defect{T}{F_{s}}{}{}{}$\\
			\greycline{2-7}
			$\defect{L}{L}{a}{x}{}$&$\defect{L}{T}{a+c}{}{}$&$\defect{L}{L}{a+c}{x+z}{}$&$\defect{L}{R}{}{}{}$&$\defect{L}{F_{0}}{x+z}{}{}$&$\defect{L}{X_{m}}{}{}{}$&$\defect{L}{F_{s}}{x+z+sa}{}{}$\\
			\greycline{2-7}
			$\defect{R}{L}{}{}{}$&$\oplus_{\alpha}\defect{R}{T}{\alpha}{}{}$&$\defect{R}{L}{}{}{}$&$\oplus_{\alpha,\zeta}\defect{R}{R}{\alpha}{\zeta}{}$&$\oplus_{\zeta}\defect{R}{F_{0}}{\zeta}{}{}$&$p \cdot \defect{R}{X_{m}}{}{}{}$&$\oplus_{\zeta}\defect{R}{F_{s}}{\zeta}{}{}$\\
			\greycline{2-7}
			$\defect{F_{0}}{L}{x}{}{}$&$\defect{F_{0}}{T}{}{}{}$&$\defect{F_{0}}{L}{x+z}{}{}$&$\oplus_{\zeta}\defect{F_{0}}{R}{\zeta}{}{}$&$\oplus_{\eta}\defect{F_{0}}{F_{0}}{x+z}{\eta}{}$&$\oplus_{\zeta}\defect{F_{0}}{X_{m}}{\zeta}{}{}$&$\defect{F_{0}}{F_{s}}{}{}{}$\\
			\toprule[1pt]
			$\defect{X_{k}}{L}{}{}{}$&$\oplus_{\alpha}\defect{X_{k}}{T}{\alpha}{}{}$&$\defect{X_{k}}{L}{}{}{}$&$p \cdot \defect{X_{k}}{R}{}{}{}$&$\oplus_{\zeta}\defect{X_{k}}{F_{0}}{\zeta}{}{}$&$\begin{cases}\oplus_{\alpha,\beta}\defect{X_{k}}{X_{k}}{\alpha}{\beta}{}&k=m\\p\cdot\defect{X_{k}}{X_{m}}{}{}{}&k\neq m\end{cases}$&$\oplus_{\alpha}\defect{X_{k}}{F_{s}}{\alpha}{}{}$\\
			\greycline{2-7}
			$\defect{F_{q}}{L}{x}{}{}$&$\defect{F_{q}}{T}{}{}{}$&$\defect{F_{q}}{L}{x+z+qc}{}{}$&$\oplus_{\zeta}\defect{F_{q}}{R}{\zeta}{}{}$&$\defect{F_{q}}{F_{0}}{}{}{}$&$\oplus_{\zeta}\defect{F_{q}}{X_{m}}{\zeta}{}{}$&$\begin{cases}\oplus_{\alpha}\defect{F_{q}}{F_{q}}{x+z}{\alpha}{}&q=s\\\defect{F_{q}}{F_{s}}{}{}{}&q\neq s\end{cases}$\\
			\toprule[1pt]
		\end{tabular}
	}
}
\rotatebox[]{90}{
	\resizebox*{\tabwidth}{\tabheight}{
		\begin{tabular}{!{\vrule width 1pt}>{\columncolor[gray]{.9}[\tabcolsep]}c!{\vrule width 1pt}c !{\color[gray]{.8}\vrule} c !{\color[gray]{.8}\vrule} c !{\color[gray]{.8}\vrule} c!{\vrule width 1pt}c !{\color[gray]{.8}\vrule} c!{\vrule width 1pt}}
			\toprule[1pt]
			\rowcolor[gray]{.9}[\tabcolsep]\cellcolor{red!20}$F_0$&$\defect{F_{0}}{T}{}{}{}$&$\defect{F_{0}}{L}{z}{}{}$&$\defect{F_{0}}{R}{z}{}{}$&$\defect{F_{0}}{F_{0}}{z}{w}{}$&$\defect{F_{0}}{X_{m}}{z}{}{}$&$\defect{F_{0}}{F_{s}}{}{}{}$\\
			\toprule[1pt]
			$\defect{T}{F_{0}}{}{}{}$&$\oplus_{\alpha,\beta}\defect{T}{T}{\alpha}{\beta}{}$&$\oplus_{\alpha}\defect{T}{L}{\alpha}{}{}$&$\oplus_{\alpha}\defect{T}{R}{\alpha}{}{}$&$\defect{T}{F_{0}}{}{}{}$&$\oplus_{\alpha}\defect{T}{X_{m}}{\alpha}{}{}$&$p \cdot \defect{T}{F_{s}}{}{}{}$\\
			\greycline{2-7}
			$\defect{L}{F_{0}}{x}{}{}$&$\oplus_{\alpha}\defect{L}{T}{\alpha}{}{}$&$\oplus_{\alpha}\defect{L}{L}{\alpha}{x+z}{}$&$\defect{L}{R}{}{}{}$&$\defect{L}{F_{0}}{x+z}{}{}$&$\defect{L}{X_{m}}{}{}{}$&$\oplus_{\zeta}\defect{L}{F_{s}}{\zeta}{}{}$\\
			\greycline{2-7}
			$\defect{R}{F_{0}}{x}{}{}$&$\oplus_{\alpha}\defect{R}{T}{\alpha}{}{}$&$\defect{R}{L}{}{}{}$&$\oplus_{\alpha}\defect{R}{R}{\alpha}{x+z}{}$&$\defect{R}{F_{0}}{x+w}{}{}$&$\defect{R}{X_{m}}{}{}{}$&$\oplus_{\zeta}\defect{R}{F_{s}}{\zeta}{}{}$\\
			\greycline{2-7}
			$\defect{F_{0}}{F_{0}}{x}{y}{}$&$\defect{F_{0}}{T}{}{}{}$&$\defect{F_{0}}{L}{x+z}{}{}$&$\defect{F_{0}}{R}{y+z}{}{}$&$\defect{F_{0}}{F_{0}}{x+z}{y+w}{}$&$\defect{F_{0}}{X_{m}}{x+z+m^{-1}y}{}{}$&$\defect{F_{0}}{F_{s}}{}{}{}$\\
			\toprule[1pt]
			$\defect{X_{k}}{F_{0}}{x}{}{}$&$\oplus_{\alpha}\defect{X_{k}}{T}{\alpha}{}{}$&$\defect{X_{k}}{L}{}{}{}$&$\defect{X_{k}}{R}{}{}{}$&$\defect{X_{k}}{F_{0}}{x+z+k^{-1}w}{}{}$&$\begin{cases}\oplus_{\alpha}\defect{X_{k}}{X_{k}}{\alpha}{k(x+z)}{}&k=m\\\defect{X_{k}}{X_{m}}{}{}{}&k\neq m\end{cases}$&$\oplus_{\zeta}\defect{X_{k}}{F_{s}}{\zeta}{}{}$\\
			\greycline{2-7}
			$\defect{F_{q}}{F_{0}}{}{}{}$&$p\cdot\defect{F_{q}}{T}{}{}{}$&$\oplus_{\zeta}\defect{F_{q}}{L}{\zeta}{}{}$&$\oplus_{\zeta}\defect{F_{q}}{R}{\zeta}{}{}$&$\defect{F_{q}}{F_{0}}{}{}{}$&$\oplus_{\zeta}\defect{F_{q}}{X_{m}}{\zeta}{}{}$&$\begin{cases}\oplus_{\zeta,\eta}\defect{F_{q}}{F_{q}}{\zeta}{\eta}{}&q=s\\ p \cdot \defect{F_{q}}{F_{s}}{}{}{}&q\neq s\end{cases}$\\
			\toprule[1pt]
		\end{tabular}
	}
}
\rotatebox[]{90}{
	\resizebox*{\tabwidth}{\tabheight}{
		\begin{tabular}{!{\vrule width 1pt}>{\columncolor[gray]{.9}[\tabcolsep]}c!{\vrule width 1pt}c !{\color[gray]{.8}\vrule} c !{\color[gray]{.8}\vrule} c !{\color[gray]{.8}\vrule} c!{\vrule width 1pt}c !{\color[gray]{.8}\vrule} c !{\color[gray]{.8}\vrule} c!{\vrule width 1pt}}
			\toprule[1pt]
			\rowcolor[gray]{.9}[\tabcolsep]\cellcolor{red!20}$F_r$&$\defect{F_{r}}{T}{}{}{}$&$\defect{F_{r}}{L}{z}{}{}$&$\defect{F_{r}}{R}{z}{}{}$&$\defect{F_{r}}{F_{0}}{}{}{}$&$\defect{F_{r}}{X_{m}}{z}{}{}$&$\defect{F_{r}}{F_{r}}{z}{w}{}$&$\defect{F_{r}}{F_{s}}{}{}{}$\\
			\toprule[1pt]
			$\defect{T}{F_{r}}{}{}{}$&$\oplus_{\alpha,\beta}\defect{T}{T}{\alpha}{\beta}{}$&$\oplus_{\alpha}\defect{T}{L}{\alpha}{}{}$&$\oplus_{\alpha}\defect{T}{R}{\alpha}{}{}$&$p\cdot\defect{T}{F_{0}}{}{}{}$&$\oplus_{\alpha}\defect{T}{X_{m}}{\alpha}{}{}$&$\defect{T}{F_{r}}{}{}{}$&$p\cdot\defect{T}{F_{s}}{}{}{}$\\
			\greycline{2-8}
			$\defect{L}{F_{r}}{x}{}{}$&$\oplus_{\alpha}\defect{L}{T}{\alpha}{}{}$&$\oplus_{\alpha}\defect{L}{L}{\alpha}{x+z-r\alpha}{}$&$\defect{L}{R}{}{}{}$&$\oplus_{\zeta}\defect{L}{F_{0}}{\zeta}{}{}$&$\defect{L}{X_{m}}{}{}{}$&$\defect{L}{F_{r}}{x+z}{}{}$&$\oplus_{\zeta}\defect{L}{F_{s}}{\zeta}{}{}$\\
			\greycline{2-8}
			$\defect{R}{F_{r}}{x}{}{}$&$\oplus_{\alpha}\defect{R}{T}{}{}{}$&$\defect{R}{L}{}{}{}$&$\oplus_{\alpha}\defect{R}{R}{\alpha}{x+z-r\alpha}{}$&$\oplus_{\zeta}\defect{R}{F_{0}}{\zeta}{}{}$&$\defect{R}{X_{m}}{}{}{}$&$\defect{R}{F_{r}}{x+w}{}{}$&$\oplus_{\zeta}\defect{R}{F_{s}}{\zeta}{}{}$\\
			\greycline{2-8}
			$\defect{F_{0}}{F_{r}}{}{}{}$&$p\cdot\defect{F_{0}}{T}{}{}{}$&$\oplus_{\zeta}\defect{F_{0}}{L}{\zeta}{}{}$&$\oplus_{\zeta}\defect{F_{0}}{R}{\zeta}{}{}$&$\oplus_{\zeta,\eta}\defect{F_{0}}{F_{0}}{\zeta}{\eta}{}$&$\sum_\alpha \defect{F_{0}}{X_{m}}{\alpha}{}{}$&$\defect{F_{0}}{F_{r}}{}{}{}$&$p\cdot\defect{F_{0}}{F_{s}}{}{}{}$\\
			\toprule[1pt]
			$\defect{X_{k}}{F_{r}}{x}{}{}$&$\oplus_{\alpha}\defect{X_{k}}{T}{\alpha}{}{}$&$\defect{X_{k}}{L}{}{}{}$&$\defect{X_{k}}{R}{}{}{}$&$\oplus_{\zeta}\defect{X_{k}}{F_{0}}{\zeta}{}{}$&$\begin{cases}\oplus_{\alpha}\defect{X_{k}}{X_{k}}{\alpha}{x+z-r\alpha}{}&k=m\\\defect{X_{k}}{X_{m}}{}{}{}&k\neq m\end{cases}$&$\defect{X_{k}}{F_{r}}{w+x+kz}{}{}$&$\oplus_{\zeta}\defect{X_{k}}{F_{s}}{\zeta}{}{}$\\
			\greycline{2-8}
			$\defect{F_{r}}{F_{r}}{x}{y}{}$&$\defect{F_{r}}{T}{}{}{}$&$\defect{F_{q}}{L}{x+z}{}{}$&$\defect{F_{q}}{R}{y+z}{}{}$&$\defect{F_{q}}{F_{0}}{}{}{}$&$\defect{F_{q}}{X_{m}}{y+z+mx}{}{}$&$\defect{F_{r}}{F_{r}}{x+z}{y+w}{}$&$\defect{F_{r}}{F_{s}}{}{}{}$\\
			\greycline{2-8}
			$\defect{F_{q}}{F_{r}}{}{}{}$&$p\cdot\defect{F_{q}}{T}{}{}{}$&$\oplus_{\zeta}\defect{F_{q}}{L}{\zeta}{}{}$&$\oplus_{\zeta}\defect{F_{q}}{R}{\zeta}{}{}$&$p\cdot\defect{F_{q}}{F_{0}}{}{}{}$&$\oplus_{\zeta}\defect{F_{q}}{X_{m}}{\zeta}{}{}$&$\defect{F_{q}}{F_{r}}{}{}{}$&$\begin{cases}\oplus_{\zeta,\eta}\defect{F_{q}}{F_{q}}{\zeta}{\eta}{}&q=s\\p\cdot\defect{F_{q}}{F_{s}}{}{}{}&q\neq s\end{cases}$\\
			\toprule[1pt]
		\end{tabular}
	}
}
\rotatebox[]{90}{
\resizebox*{\tabwidth}{\tabheight}{
	\begin{tabular}{!{\vrule width 1pt}>{\columncolor[gray]{.9}[\tabcolsep]}c!{\vrule width 1pt}c !{\color[gray]{.8}\vrule} c !{\color[gray]{.8}\vrule} c !{\color[gray]{.8}\vrule} c!{\vrule width 1pt}c !{\color[gray]{.8}\vrule} c!{\vrule width 1pt}}
		\toprule[1pt]
		\rowcolor[gray]{.9}[\tabcolsep]\cellcolor{red!20}$T$&$\defect{T}{T}{c}{d}{}$&$\defect{T}{L}{c}{}{}$&$\defect{T}{R}{c}{}{}$&$\defect{T}{F_{0}}{}{}{}$&$\defect{T}{X_{m}}{c}{}{}$&$\defect{T}{F_{s}}{}{}{}$\\
		\toprule[1pt]
		$\defect{T}{T}{a}{b}{}$&$\defect{T}{T}{a+c}{b+d}{}$&$\defect{T}{L}{b+c}{}{}$&$\defect{T}{R}{a+c}{}{}$&$\defect{T}{F_{0}}{}{}{}$&$\defect{T}{X_{m}}{a+c+mb}{}{}$&$\defect{T}{F_{s}}{}{}{}$\\
				\greycline{2-7}
		$\defect{L}{T}{a}{}{}$&$\defect{L}{T}{a+d}{}{}$&$\oplus_{\beta}\defect{L}{L}{a+c}{\beta}{}$&$\defect{L}{R}{}{}{}$&$\oplus_{\alpha}\defect{L}{F_{0}}{\alpha}{}{}$&$\defect{L}{X_{m}}{}{}{}$&$\oplus_{\alpha}\defect{L}{F_{s}}{\alpha}{}{}$\\
				\greycline{2-7}
		$\defect{R}{T}{a}{}{}$&$\defect{R}{T}{a+c}{}{}$&$\defect{R}{L}{}{}{}$&$\oplus_{\beta}\defect{R}{R}{a+c}{\beta}{}$&$\oplus_\alpha \defect{R}{F_{0}}{\alpha}{}{}$&$\defect{R}{X_{m}}{}{}{}$&$\oplus_{\alpha}\defect{R}{F_{s}}{\alpha}{}{}$\\
				\greycline{2-7}
		$\defect{F_{0}}{T}{}{}{}$&$\defect{F_{0}}{T}{}{}{}$&$\oplus_{\alpha}\defect{F_{0}}{L}{\alpha}{}{}$&$\oplus_{\alpha}\defect{F_{0}}{R}{\alpha}{}{}$&$\oplus_{\alpha,\beta}\defect{F_{0}}{F_{0}}{\alpha}{\beta}{}$&$\oplus_{\alpha}\defect{F_{0}}{X_{m}}{\alpha}{}{}$&$p\cdot\defect{F_{0}}{F_{s}}{}{}{}$\\
		\toprule[1pt]
		$\defect{X_{k}}{T}{a}{}{}$&$\defect{X_{k}}{T}{a+c+kd}{}{}$&$\defect{X_{k}}{L}{}{}{}$&$\defect{X_{k}}{R}{}{}{}$&$\oplus_{\alpha}\defect{X_{k}}{F_{0}}{\alpha}{}{}$&$\begin{cases}\oplus_{\beta}\defect{X_{k}}{X_{k}}{a+c}{\beta}{}&k=m\\\defect{X_{k}}{X_{m}}{}{}{}&k\neq m\end{cases}$&$\oplus_{\alpha}\defect{X_{k}}{F_{s}}{\alpha}{}{}$\\
				\greycline{2-7}
		$\defect{F_{q}}{T}{}{}{}$&$\defect{F_{q}}{T}{}{}{}$&$\oplus_{\alpha}\defect{F_{q}}{L}{\alpha}{}{}$&$\oplus_{\alpha}\defect{F_{q}}{R}{\alpha}{}{}$&$p\cdot\defect{F_{q}}{F_{0}}{}{}{}$&$\oplus_{\alpha}\defect{F_{q}}{X_{m}}{\alpha}{}{}$&$\begin{cases}\oplus_{\alpha,\beta}\defect{F_{q}}{F_{q}}{\alpha}{\beta}{}&q=s\\p\cdot\defect{F_{q}}{F_{s}}{}{}{}&q\neq s\end{cases}$\\
		\toprule[1pt]
	\end{tabular}
}
}
\rotatebox[]{90}{
\resizebox*{\tabwidth}{\tabheight}{
	\begin{tabular}{!{\vrule width 1pt}>{\columncolor[gray]{.9}[\tabcolsep]}c!{\vrule width 1pt}c !{\color[gray]{.8}\vrule} c !{\color[gray]{.8}\vrule} c !{\color[gray]{.8}\vrule} c!{\vrule width 1pt}c !{\color[gray]{.8}\vrule} c!{\vrule width 1pt}}
		\toprule[1pt]
		\rowcolor[gray]{.9}[\tabcolsep]\cellcolor{red!20}$R$&$\defect{R}{T}{c}{}{}$&$\defect{R}{L}{}{}{}$&$\defect{R}{R}{c}{z}{}$&$\defect{R}{F_{0}}{z}{}{}$&$\defect{R}{X_{m}}{}{}{}$&$\defect{R}{F_{s}}{z}{}{}$\\
		\toprule[1pt]
		$\defect{T}{R}{a}{}{}$&$\oplus_{\beta}\defect{T}{T}{a+c}{\beta}{}$&$\oplus_{\alpha}\defect{T}{L}{\alpha}{}{}$&$\defect{T}{R}{a+c}{}{}$&$\defect{T}{F_{0}}{}{}{}$&$\oplus_{\alpha}\defect{T}{X_{m}}{\alpha}{}{}$&$\defect{T}{F_{s}}{}{}{}$\\
		\greycline{2-7}
		$\defect{L}{R}{}{}{}$&$\oplus_{\alpha}\defect{L}{T}{\alpha}{}{}$&$\oplus_{\alpha,\zeta}\defect{L}{L}{\alpha}{\zeta}{}$&$\defect{L}{R}{}{}{}$&$\oplus_{\zeta}\defect{L}{F_{0}}{\zeta}{}{}$&$p\cdot\defect{L}{X_{m}}{}{}{}$&$\oplus_{\zeta}\defect{L}{F_{s}}{\zeta}{}{}$\\
		\greycline{2-7}
		$\defect{R}{R}{a}{x}{}$&$\defect{R}{T}{a+c}{}{}$&$\defect{R}{L}{}{}{}$&$\defect{R}{R}{a+c}{x+z}{}$&$\defect{R}{F_{0}}{x+z}{}{}$&$\defect{R}{X_{m}}{}{}{}$&$\defect{R}{F_{s}}{x+z+sa}{}{}$\\
		\greycline{2-7}
		$\defect{F_{0}}{R}{x}{}{}$&$\defect{F_{0}}{T}{}{}{}$&$\oplus_{\zeta}\defect{F_{0}}{L}{\zeta}{}{}$&$\defect{F_{0}}{R}{x+z}{}{}$&$\oplus_{\zeta}\defect{F_{0}}{F_{0}}{\zeta}{x+z}{}$&$\oplus_{\zeta}\defect{F_{0}}{X_{m}}{\zeta}{}{}$&$\defect{F_{0}}{F_{s}}{}{}{}$\\
		\toprule[1pt]
		$\defect{X_{k}}{R}{}{}{}$&$\oplus_{\alpha}\defect{X_{k}}{T}{\alpha}{}{}$&$p\cdot\defect{X_{k}}{L}{}{}{}$&$\defect{X_{k}}{R}{}{}{}$&$\oplus_{\zeta}\defect{X_{k}}{F_{0}}{\zeta}{}{}$&$\begin{cases}\oplus_{\alpha,\zeta}\defect{X_{k}}{X_{k}}{\alpha}{\zeta}{}&k=m\\ p\cdot \defect{X_{k}}{X_{m}}{}{}{}&k\neq m\end{cases}$&$\oplus_{\zeta}\defect{X_{k}}{F_{s}}{\zeta}{}{}$\\
		\greycline{2-7}
		$\defect{F_{q}}{R}{x}{}{}$&$\defect{F_{q}}{T}{}{}{}$&$\oplus_{\zeta}\defect{F_{q}}{L}{\zeta}{}{}$&$\defect{F_{q}}{R}{x+z+qc}{}{}$&$\defect{F_{q}}{F_{0}}{}{}{}$&$\oplus_{\zeta}\defect{F_{q}}{X_{m}}{\zeta}{}{}$&$\begin{cases}\oplus_{\zeta}\defect{F_{q}}{F_{q}}{\zeta}{x+z}{}&q=s\\\defect{F_{q}}{F_{s}}{}{}{}&q\neq s\end{cases}$\\
		\toprule[1pt]
	\end{tabular}
}
}
\rotatebox[]{90}{
\resizebox*{\tabwidth}{\tabheight}{
	\begin{tabular}{!{\vrule width 1pt}>{\columncolor[gray]{.9}[\tabcolsep]}c!{\vrule width 1pt}c !{\color[gray]{.8}\vrule} c !{\color[gray]{.8}\vrule} c !{\color[gray]{.8}\vrule} c !{\vrule width 1pt} c!{\color[gray]{.8}\vrule}c !{\color[gray]{.8}\vrule} c!{\vrule width 1pt}}
		\toprule[1pt]
		\rowcolor[gray]{.9}[\tabcolsep]\cellcolor{red!20}$X_{l}$&$\defect{X_{l}}{T}{c}{}{}$&$\defect{X_{l}}{L}{}{}{}$&$\defect{X_{l}}{R}{}{}{}$&$\defect{X_{l}}{F_{0}}{z}{}{}$&$\defect{X_{l}}{X_{l}}{c}{z}{}$&$\defect{X_{l}}{X_{m}}{}{}{}$&$\defect{X_{l}}{F_{s}}{z}{}{}$\\
		\toprule[1pt]
		$\defect{T}{X_{l}}{a}{}{}$&$\oplus_{\beta}\defect{T}{T}{a+c-l\beta}{\beta}{}$&$\oplus_{\alpha}\defect{T}{L}{\alpha}{}{}$&$\oplus_{\alpha}\defect{T}{R}{\alpha}{}{}$&$\defect{T}{F_{0}}{}{}{}$&$\defect{T}{X_{m}}{a+c}{}{}$&$\oplus_{\alpha}\defect{T}{X_{m}}{\alpha}{}{}$&$\defect{T}{F_{s}}{}{}{}$\\
		\greycline{2-8}
		$\defect{L}{X_{l}}{}{}{}$&$\oplus_{\alpha}\defect{L}{T}{\alpha}{}{}$&$\oplus_{\alpha,\zeta}\defect{L}{L}{\alpha}{\zeta}{}$&$p\cdot\defect{L}{R}{}{}{}$&$\oplus_{\zeta}\defect{L}{F_{0}}{\zeta}{}{}$&$\defect{L}{X_{m}}{}{}{}$&$p\cdot\defect{L}{X_{m}}{}{}{}$&$\oplus_{\zeta}\defect{L}{F_{s}}{\zeta}{}{}$\\
		\greycline{2-8}
		$\defect{R}{X_{l}}{}{}{}$&$\oplus_{\alpha}\defect{R}{T}{\alpha}{}{}$&$p\cdot\defect{R}{L}{}{}{}$&$\oplus_{\alpha,\zeta}\defect{R}{R}{\alpha}{\zeta}{}$&$\oplus_{\zeta}\defect{R}{F_{0}}{\zeta}{}{}$&$\defect{R}{X_{m}}{}{}{}$&$p\cdot\defect{R}{X_{m}}{}{}{}$&$\oplus_{\zeta}\defect{R}{F_{s}}{\zeta}{}{}$\\
		\greycline{2-8}
		$\defect{F_{0}}{X_{l}}{x}{}{}$&$\defect{F_{0}}{T}{}{}{}$&$\oplus_{\zeta}\defect{F_{0}}{L}{\zeta}{}{}$&$\oplus_{\zeta}\defect{F_{0}}{R}{\zeta}{}{}$&$\oplus_{\eta}\defect{F_{0}}{F_{0}}{x+z-l^{-1}\eta}{\eta}{}$&$\defect{F_{0}}{X_{l}}{x+l^{-1}z}{}{}$&$\oplus_{\zeta}\defect{F_{0}}{X_{m}}{\zeta}{}{}$&$\defect{F_{0}}{F_{s}}{}{}{}$\\
		\toprule[1pt]
		$\defect{X_{l}}{X_{l}}{a}{x}{}$&$\defect{X_{k}}{T}{a+c}{}{}$&$\defect{X_{k}}{L}{}{}{}$&$\defect{X_{k}}{R}{}{}{}$&$\defect{X_{k}}{F_{0}}{l^{-1}x+z}{}{}$&$\defect{X_{l}}{X_{l}}{a+c}{x+z}{}$&$\defect{X_{l}}{X_{m}}{}{}{}$&$\defect{X_{k}}{F_{s}}{x+z+as}{}{}$\\
		\greycline{2-8}
		$\defect{X_{k}}{X_{l}}{}{}{}$&$\oplus_{\alpha}\defect{X_{k}}{T}{\alpha}{}{}$&$p\cdot\defect{X_{k}}{L}{}{}{}$&$p\cdot\defect{X_{k}}{R}{}{}{}$&$\oplus_{\zeta}\defect{X_{k}}{F_{0}}{\zeta}{}{}$&$\defect{X_{k}}{X_{l}}{}{}{}$&$\begin{cases}\oplus_{\alpha,\zeta}\defect{X_{k}}{X_{k}}{\alpha}{\zeta}{}&k=m\\p\cdot\defect{X_{k}}{X_{m}}{}{}{}&k\neq m\end{cases}$&$\oplus_{\zeta}\defect{X_{k}}{F_{s}}{\zeta}{}{}$\\
		\greycline{2-8}
		$\defect{F_{q}}{X_{l}}{x}{}{}$&$\defect{F_{q}}{T}{}{}{}$&$\oplus_{\zeta}\defect{F_{q}}{L}{\zeta}{}{}$&$\oplus_{\zeta}\defect{F_{q}}{R}{\zeta}{}{}$&$\defect{F_{q}}{F_{0}}{}{}{}$&$\defect{F_{q}}{X_{m}}{x+z+qc}{}{}$&$\oplus_{\zeta}\defect{F_{q}}{X_{m}}{\zeta}{}{}$&$\begin{cases}\oplus_{\zeta}\defect{F_{q}}{F_{q}}{\zeta}{x+z-l\zeta}{}&q=s\\\defect{F_{q}}{F_{s}}{}{}{}&q\neq s\end{cases}$\\
		\toprule[1pt]
	\end{tabular}
}
}
\end{minipage}
\vspace*{14mm}
	\caption{Vertical fusion tables}\label{tab:vertical_fusion_tables}
	\renewcommand{\arraystretch}{1}
\vspace*{-20mm}
\end{table}

\end{document}